\documentclass[12pt, a4paper]{article}

\usepackage[english]{babel}
\usepackage[right=2.8cm,left=2.8cm,top=3.2cm,bottom=3.0cm]{geometry}

\usepackage{authblk}

\usepackage{setspace}
\usepackage{multirow}
\usepackage{graphicx}
\usepackage[version=4]{mhchem}
\usepackage{array}
\usepackage{bm}
\usepackage{fullpage}
\usepackage{pgfplots,tikz}
\usepackage{pgfplotstable}
\usepackage{hyperref}

\usepackage{mathtools}

\usepackage{amsmath}
\usepackage{amssymb}
\usepackage{amsthm}
\usepackage[latin1]{inputenc}
\usepackage{subcaption}
\usepackage{graphicx}
\usepackage{ulem}
\usepackage{graphics}
\usepackage{xcolor}
\usepackage{float}
\usepackage{array}
\usepackage{url}
\usepackage{blindtext}
\definecolor{lightgray}{gray}{0.9}
\usepackage{listings}
\usepackage[toc,page]{appendix}
\usepackage{textcomp}

\usepackage{epstopdf}
\usepackage{enumitem}
\usepackage{comment}
\usepackage{mathrsfs}
\usepackage{cite}
\usepackage[final]{pdfpages}
\usepackage{gensymb}

\usepackage{caption}
\usepackage{color}

\newcommand*\diff{\mathop{}\!\mathrm{d}}

\newcommand{\odiff}[2]{\frac{\diff #1}{\diff #2}}

\newcommand{\pdiff}[2]{\frac{\partial #1}{\partial #2}}

\usepackage{amsthm}
\theoremstyle{plain}

\theoremstyle{definition}

\usepackage{color} 
\definecolor{mygreen}{RGB}{28,172,0} 
\definecolor{mylilas}{RGB}{170,55,241}

\newcommand{\MA}{\leavevmode Monge-Amp\`ere}
\newcommand{\MAe}{\leavevmode Monge-Amp\`ere equation}
\newcommand{\MAE}{\leavevmode Monge-Amp\`ere Equation}

\usepackage[linesnumbered]{algorithm2e}
\usepackage{chngcntr}
\counterwithout{equation}{section}

\definecolor{deepgreen}{rgb}{0,0.5,0}
\lstdefinelanguage{Python}{
	keywords={typeof, null, catch, switch, in, int, str, float, self, as, import},
	keywordstyle=\color{blue}\bfseries,
	ndkeywords={boolean, throw, import},
	ndkeywords={return, class, if ,elif, endif, while, do, else, True, False , catch, def},
	ndkeywordstyle=\color{blue}\bfseries,
	identifierstyle=\color{black},
	sensitive=true,
	comment=[l]{\#},
	morecomment=[s]{/*}{*/},
	commentstyle=\color{deepgreen}\ttfamily,
	stringstyle=\color{red}\ttfamily,
	tabsize=4,
}

\usepackage{pgfplots}
\usetikzlibrary{arrows.meta}
\usetikzlibrary{pgfplots.patchplots}
\usetikzlibrary{calc}
\usetikzlibrary{positioning}
\usetikzlibrary{decorations.pathreplacing}
\usetikzlibrary{decorations.markings}

\usepackage{calc}
\makeatletter
\newcounter{tmp@cnt}
\newcommand*\@labelpunc{.}
\newcommand*\combine[1][2]{%
	\refstepcounter{enumi}
	\setcounter{tmp@cnt}{\value{enumi}}
	\addtocounter{enumi}{#1-1}
	\item[\thetmp@cnt--\theenumi\@labelpunc]}
\newcommand*\labeltype[2][]{\gdef\@labelpunc{#1}\renewcommand\thetmp@cnt{#2{tmp@cnt}}}
\makeatother

\newlength\tindent
\usepackage{afterpage} 
\usepackage[toc,page]{appendix}
\usepackage{empheq}
\usepackage[nottoc]{tocbibind}

\providecommand{\keywords}[1]
{
    \small	
    \textbf{\textit{Keywords---}} #1
}

\title{Numerical Methods for the Hyperbolic \MAE{} Based on the Method of Characteristics}
\date{}
\author[1,*]{M.W.M.C. Bertens}
\author[1]{E.M.T. Vugts}
\author[1]{{}M.J.H. Anthonissen}
\author[1]{\\J.H.M. ten Thije Boonkkamp}
\author[1,2]{W.L. IJzerman}
\affil[1]{CASA, Department of Mathematics and Computer Science, Eindhoven University of Technology, PO Box 513, 5600 MB Eindhoven, The Netherlands}
\affil[2]{Signify Research, High Tech Campus 7, 5656 AE Eindhoven, The Netherlands}
\affil[*]{Corresponding author: m.w.m.c.bertens@tue.nl}
\begin{document}

    
    \maketitle

    \abstract{
    \noindent{}We present three alternative derivations of the method of characteristics (MOC) for a second order nonlinear hyperbolic partial differential equation. The MOC gives rise to two mutually coupled systems of ordinary differential equations. As a special case we consider the \MA{} (MA) equation, for which we solve the system of ODE's using explicit one-step methods (Euler, Runge-Kutta) and spline interpolation. Numerical examples demonstrate the performance of the methods. 
    }


    \section{Introduction}
    The general \MAe{} for a variable $u$ in two independent variables $x,y$ is of the form
    \begin{align}
    A(u_{xx} u_{yy} - u_{xy}^2) + B u_{xx} + C u_{xy} + D u_{yy} + E = 0,
    \end{align}
    where $A, B, C, D$ and $E$ are functions, possibly dependent on $x, y, u, u_x$ and $u_y$. The linearity in the Hessian $u_{xx} u_{yy} - u_{xy}^2$ is the defining feature of the \MAe{}.
    Applications of the \MAe{} are found, a.o., in fluid dynamics to compute the velocity of an incompressible fluid from the pressure using the streamline formulation~\cite{fluid_dynamics_MAE}, in mathematical finance to determine optimal portfolio strategies~\cite{Caboussat2014MathFinance} and in Riemannian geometry to compute the surface of a manifold given the Gauss curvature~\cite{Chen_2019}.
     
    Our interest lies in designing freeform optical surfaces, i.e., mirrors or lenses without any symmetries, to transfer a given light source distribution to a desired target distribution for some optical systems. Combining the optical map, i.e., the relation between a point in the  source domain and a point in the target domain, with conservation of energy gives rise to two variants of the \MAe{}, viz. the elliptic and the hyperbolic equation~\cite{CorienThesis}.
    The elliptic equation is well established in the literature and used in optical design \cite{CorienThesis, Yadav_2019, Romijn_2019, Romijn_2021}. On the other hand, the hyperbolic equation is more exotic and the literature is scarce. The papers \cite{Brickell_1976, Westcott_1976} are the most notable references in the context of numerical results in illumination optics, and the papers \cite{GlobalSmoothSolutions, Tunitsky_2017} are the most important results regarding existence and uniqueness results. 
    It is conjectured that designing optical systems using the hyperbolic equation allows for the construction of more compact optics.
    
\begin{figure}[H]
        \centering
        \captionsetup{width=0.9\linewidth}
        \includegraphics[width = 0.55\linewidth]{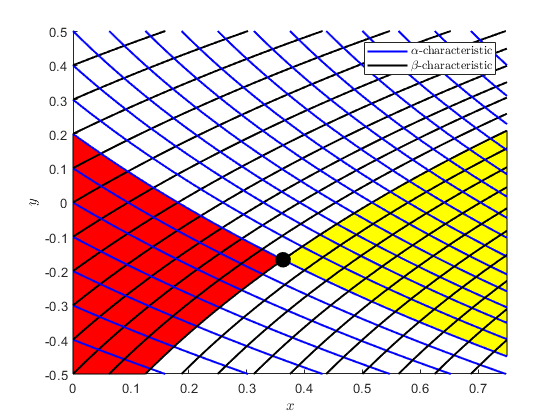}
        \caption{ Schematic representation of domain of dependence (red) and domain of influence (yellow). }
        \label{fig:domainOfDependenceAndInfluence}
    \end{figure}

    The hyperbolic \MAe{} has proven to be more difficult to solve than its elliptic counterpart. This is due to the existence of two mutually coupled families of characteristics. 
    The two characteristics through an interior point $(x_{0},y_{0})$ facing back to the boundary enclose the domain of dependence of this point. The solution $u(x_{0},y_{0})$ depends on all function values $u(x,y)$ with $(x,y)$ in this domain. Conversely, the characteristics emanating from $(x_{0},y_{0})$ bound the region of influence of $(x_{0},y_{0})$, which is the region where the solution is determined by $u(x_{0},y_{0})$.   
    Figure~\ref{fig:domainOfDependenceAndInfluence} shows one such example where the blue and black lines indicate (a few of the) characteristics, where $(x_0, y_0) = (0.363, -0.167)$ which is denoted by the black dot and where the red and yellow parts indicate the domain of dependence and the domain of influence, respectively.
    Hence, boundary data determine the solution in the interior domain, and vice versa. Therefore, we distinguish between entering and leaving characteristics and fix a so-called initial strip, determined by a chosen parameterization, on which we prescribe Cauchy conditions. The remaining boundary conditions then follow from the course of the characteristics, by considering the domain of dependence, i.e., by considering where along the boundary, the characteristics enter or leave the domain. Depending on the number of characteristics entering and leaving, we can either solely prescribe $u$, prescribe $u$ and the normal derivative of $u$, or we should not specify any boundary conditions at all. Therefore, it is of utmost importance that a numerical solution procedure is able to accurate approximate the location of the characteristics and identify the correct boundary conditions. In other words, a numerical method that violates either one of these conditions will completely destroy the solution.
    
    To start with, we introduce a general framework for second order nonlinear hyperbolic PDEs and subsequently restrict ourselves to the special case of the standard hyperbolic \MAe{}, viz.
    \begin{align}
    u_{xx} u_{yy} - u_{xy}^2 + f^2 = 0,
    \end{align}
    where $f$ is a given continuously differentiable function dependent on $x$ and $y$. The method of characteristics gives rise to two mutually coupled ODE systems, which we integrate with standard explicit one-step methods in the $x$-direction. Unfortunately, the direction of one characteristic depends on the other and vice versa, and characteristics do not necessarily pass from a numerical grid line to the next. Therefore missing information on a characteristic and in the grid points is obtained by spline interpolation. In order to control the interpolation error, we control the displacement of the characteristics in the $y$-direction by tuning the step size in the $x$-direction.
    This consequently determines the numerical domain of dependence, which is defined similar to the regular domain of dependence, but is formed by the numerical approximations of the characteristics instead.
    Furthermore, for a scheme for a hyperbolic PDE to be numerically stable, the (physical) domain of dependence should be enclosed by the numerical domain of dependence \cite[p. 366]{JanBook}. By appropriately interpolating only within the area enclosed by the numerical characteristics and by step size control, the analytical domain of dependence lies within the domain of dependence of the numerical scheme and the developed schemes are stable in practice.
    
    Estimates for the rate of convergence of the numerical methods are made, and tested for a variety of examples. We measure the convergence indirectly because direct measures are generally infeasible as they rely on analytical solutions, which are often unavailable. We do so by reformulating the \MAe{} as an integral equation and measuring its residual via Gauss-Legendre quadrature rules. We present examples which confirm $u$ to be a saddle surface, one example where the number of required boundary conditions varies along the boundary, and one example with discontinuous third derivatives, for which no analytical solution is known.
    
    We have organized our paper as follows. The theoretical framework for a hyperbolic second order PDE is introduced in Section~\ref{sec:DerivationODEs}. In Section~\ref{sec:HMAE} we apply this to the hyperbolic \MAe{}, and discuss the boundary conditions. Subsequently, we introduce the numerical methods in Section~\ref{sec:numericalMethods}. In Section~\ref{sec:numericalResults}, various numerical results are given and analyzed, and finally a brief discussion and concluding remarks are given in Section~\ref{sec:Conclusion}.
    
    \section{Method of characteristics for a second {order} nonlinear hyperbolic PDE}
    \label{sec:DerivationODEs}
    We start by introducing the method of characteristics for a general nonlinear second order PDE in two variables. To this end we assume, unless explicitly stated otherwise, that all functions are continuous and have continuous derivatives of all orders involved. Let the PDE of interest be given by
    \begin{align}
    \label{eqn:nonlinearBaseEqn}
    F(x,y,u,p,q,r,s,t) = 0, \hspace{5pt} (x,y) \in \Omega,
    \end{align}
    where $u = u(x,y)$, $p = u_x$, $q=u_y$, $r=u_{xx}$, $s=u_{xy}$, $t=u_{yy}$ and $\Omega \subseteq \mathbb{R}^2$ the domain of interest.

    \subsection{An introduction to the method of characteristics}
     In this section we give a brief introduction to the method of characteristics.
     
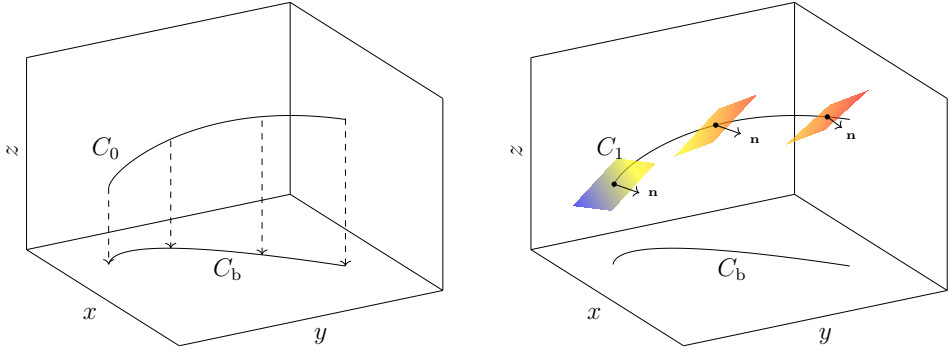
\begin{figure}[H]
         \centering
         \captionsetup{width=0.9\linewidth}
         \scalebox{0.8}{
         \begin{tikzpicture}
             \begin{axis}
             [
             ticks=none,
             xlabel={$x$},
             ylabel={$y$},
             zlabel={$z$},
             xmin=-.5,
             xmax=1.5,
             ymin=0.5,
             ymax=3.5,
             zmin=0,
             zmax=2.5,
             view={60}{30},
             ]
             
             \addplot3[
             domain=0:1,
             samples = 60,
             samples y=0,
             ]
             ({1.8 * x^2 - x},
             {1 + x^2 + x},
             {1 + x * sin(deg(2*x))});
             
             \addplot3[
             domain=0:1,
             samples = 60,
             samples y=0,
             ]
             ({1.8 * x^2 - x},
             {1 + x^2 + x},
             {0});

              \draw[->, dashed] (axis cs:0, 1, 1) -- (axis cs:0, 1, 0) ;
              \draw[->, dashed] (axis cs:-0.05, 1.75, 1.4) -- (axis cs:-0.05, 1.75, 0) ;
              \draw[->, dashed] (axis cs:0.3520, 2.44, 1.75) -- (axis cs:0.3520, 2.44, 0) ;
              \draw[->, dashed] (axis cs:0.8, 3, 1.9212) -- (axis cs:0.8, 3, 0) ;
              \end{axis}
                          
             \node[label={$C_\text{b}$}] at (3.5, 1, 0.5) {};
             \node[label={$C_0$}] at (1.5, 3, 0.5) {};
         \end{tikzpicture}
         \hspace{20pt}
         \begin{tikzpicture}
             \begin{axis}
             [
             ticks=none,
             xlabel={$x$},
             ylabel={$y$},
             zlabel={$z$},
             xmin=-.5,
             xmax=1.5,
             ymin=0.5,
             ymax=3.5,
             zmin=0,
             zmax=2.5,
             view={60}{30},
             ]
             
             \addplot3[
             domain=0:1,
             samples = 60,
             samples y=0,
             ]
             ({1.8 * x^2 - x},
             {1 + x^2 + x},
             {1 + x * sin(deg(2*x))});
             
             \addplot3[
             domain=0:1,
             samples = 60,
             samples y=0,
             ]
             ({1.8 * x^2 - x},
             {1 + x^2 + x},
             {0});

             \addplot3[
             surf,
             shader = interp,
             opacity = 0.65,
             domain = -0.2955:0.2045,
             domain y = 0.8025:1.3025,
             ] {0.5 * x + 1 * y - 0.0247583};
             \addplot3[mark=*,black, mark size=1pt] coordinates {(-0.0455, 1.0525, 1.00499)};
             \draw[->] (axis cs:-0.0455, 1.0525, 1.00499) -- (axis cs:-0.0455+0.1, 1.0525+0.2, 1.00499-0.1) node [below, right, xshift=0.0cm, yshift =-0.25]  {\tiny $\mathbf{n}$};

             \addplot3[
             surf,
             shader = interp,
             opacity = 0.65,
             domain = -0.1395:0.3605,
             domain y = 1.8225:2.3225,
             ] {1 * x + 1 * y - 0.556687};
             \addplot3[mark=*,black, mark size=1pt] coordinates {(0.1105, 2.0725, 1.62631)};
             \draw[->] (axis cs:0.1105, 2.0725, 1.62631) -- (axis cs:0.1105+0.1, 2.0725+0.2, 1.62631-0.1) node [right, xshift=0.0cm, yshift=-0.1cm] {\tiny $\mathbf{n}$};

             \addplot3[
             surf,
             shader = interp,
             opacity = 0.65,
             domain = 0.4245:0.9245,
             domain y = 2.6025:3.1025,
             ] {1 * x + 4/5 * y - 1.05751};
             \addplot3[mark=*,black, mark size=1pt] coordinates {(0.6745, 2.8525, 1.89899)};
             \draw[->] (axis cs:0.6745, 2.8525, 1.89899) -- (axis cs:0.6745+0.1, 2.8525+4/50, 1.89899-0.1) node [below, xshift=0.15cm, yshift=0.05cm] {\tiny $\mathbf{n}$};
             
             \end{axis}
             
             \node[label={$C_\text{b}$}] at (3.5, 1, 0.5) {};
             \node[label={$C_1$}] at (1.5, 3, 0.5) {};
         \end{tikzpicture}
         }
         \caption{Schematic representation of the curves $C_\text{b}$, $C_0$ and the $C_1$-strip. Three tangent planes and their normal vectors $\mathbf{n}$ are drawn.}
         \label{fig:CbC0C1}
     \end{figure}

     Let $C_\text{b}$ be a curve in the $(x,y)$-plane, parameterized by $\lambda \in I$, for an interval $I \subset \mathbb{R}$, i.e., $C_\text{b} = \{(X(\lambda), Y(\lambda)) | \lambda \in I\}$ with $X,Y:I\rightarrow\mathbb{R}$. 
     Let $C_0$ be a corresponding curve in $(x,y,z)$-space, which we also parameterize by $\lambda \in I$, i.e., $C_0 = \{(X(\lambda), Y(\lambda), U(\lambda)) | \lambda \in I\}$ where $U:I\rightarrow\mathbb{R}$. The projection of $C_0$ on the $(x,y)$-plane yields the curve $C_\text{b}$, see Figure~\ref{fig:CbC0C1}. We call $C_\text{b}$ the base curve of $C_0$, or simply the base curve.
     
     A base curve $C_\text{b}$ is said to be differentiable if the corresponding map $\lambda \mapsto (X(\lambda), Y(\lambda))$ is differentiable for every $\lambda \in I$.
     A curve is regular if it is differentiable and the tangent vector has non-zero length for all $\lambda \in I$. We generally assume $C_\text{b}$ to be regular, implying that $\odiff{}{\lambda} (X(\lambda),Y(\lambda))^\text{T} \neq \mathbf{0}$, or equivalently $X_\lambda^2 + Y_\lambda^2 \neq 0$ for all $\lambda \in I$, where a subscript denotes differentiation.
     
     Let $\mathbf{v} = (X_\lambda, Y_\lambda, U_\lambda)^\text{T}$ be the tangent vector to $(X(\lambda), Y(\lambda), U(\lambda)) \in C_0$ with $\lambda \in I$. A plane through the point $(X(\lambda), Y(\lambda), U(\lambda))$ with normal vector $\mathbf{n}$ is tangent to the curve $C_0$ if $\mathbf{v} \cdot \mathbf{n} = 0$. To identify those planes let $\mathbf{n} = (P, Q, -1)^\text{T}$ with $P,Q:I\rightarrow\mathbb{R}$. Note that if the third component $n_3 \neq 0$,  $\mathbf{n}$ can always be reduced to such form by scaling the components. If for each point on $C_0$ we fix the tangent plane, then the collection of $C_0$ together with said tangent planes forms a so-called $C_1$-strip, i.e., 
     \begin{align}
     \label{eqn:C_1_def}
     C_1 = \{(X(\lambda), Y(\lambda), U(\lambda), P(\lambda), Q(\lambda)) | \lambda \in I\},
     \end{align}
     sometimes referred to as a strip of first order. Figure~\ref{fig:CbC0C1} shows the $C_1$-strip for three tangent planes with corresponding normals. We use the notation $C_1$ interchangeably to denote either the strip's type, or the strip itself as given by~\eqref{eqn:C_1_def}.
     From $\mathbf{v} \cdot \mathbf{n} = 0$ it follows that
     \begin{align}
     \label{eqn:1stOrderStripCond}
     P X_\lambda + Q Y_\lambda - U_\lambda = 0,
     \end{align}
     which is the strip condition of first order, in short, the strip condition.
     Note that so far the $C_1$-strip and the strip condition have no connection to the PDE~\eqref{eqn:nonlinearBaseEqn}.
     
     Let $u(x,y)$ be a solution of~\eqref{eqn:nonlinearBaseEqn}, then $z = u(x,y)$ is called an integral surface of $\eqref{eqn:nonlinearBaseEqn}$.
     An integral surface $z = u(x,y)$ naturally induces a $C_1$-strip. 
     Given a base curve $C_\text{b}$, let $u(\lambda) := u(X(\lambda), Y(\lambda))$, $p(\lambda) := p(X(\lambda), Y(\lambda))$ and $q(\lambda) := q(X(\lambda), Y(\lambda))$.
%
%
%
%
     The normal of the integral surface $u(x,y) - z = 0$ is given by $\mathbf{n} = (u_x, u_y, -1)^\text{T}$ in $(x,y,z)$-space. Hence the strip $C_1 = \{(x(\lambda), y(\lambda), u(\lambda), p(\lambda), q(\lambda)) | \lambda \in I\}$ is obtained. From the chain rule we conclude
     \begin{align}
     \label{eqn:stripConditionU}
     u_\lambda = u_x X_\lambda + u_y Y_\lambda =  p X_\lambda + q Y_\lambda.
     \end{align}
    which is identical to the strip condition~\eqref{eqn:1stOrderStripCond} with $P = p = u_x$, $Q = q = u_y$ and $U = u$, the solution of~\eqref{eqn:nonlinearBaseEqn}.
    
    One can naturally generalize first order strips to higher order strips. A $C_2$-strip consists of the $C_1$-strip together with the tangent planes of the curves \hfill \newline $(X(\lambda), Y(\lambda), P(\lambda))$ and $(X(\lambda), Y(\lambda), Q(\lambda))$.
    Higher order strip conditions are also found naturally in the following way: with $(X(\lambda), Y(\lambda), P(\lambda))$ we can associate two functions $R,S:I\rightarrow\mathbb{R}$ such that the normal vector of a tangent plane is $(R, S, -1)$. The tangent vector of $(X(\lambda), Y(\lambda), P(\lambda))$ equals $(X_\lambda, Y_\lambda, P_\lambda)$. The same reasoning as before applies and we find
    \begin{align}
    P_\lambda = R X_\lambda + S Y_\lambda.
    \end{align}
    Analogously, for $(X(\lambda), Y(\lambda), Q(\lambda))$ let the normal vector of a tangent plane be $(\tilde{S}, T, -1)$, it then follows that
    \begin{align}
    Q_\lambda = \tilde{S} X_\lambda + T Y_\lambda.
    \end{align}
    Note that the functions $S$ and $\tilde{S}$ are not necessarily equal as $(R,S,-1)^\text{T}$ should be perpendicular to the curve $(X_\lambda, Y_\lambda, P_\lambda)$ and $(\tilde{S}, T,-1)^\text{T}$ should be perpendicular to the curve $(X_\lambda, Y_\lambda, Q_\lambda)$. 
     As before, an integral surface $z=u(x,y)$ induces a $C_2$-strip where we identify $Q$, $R$, $S$, $\tilde{S}$ and $T$ with the values $r, s, t$ via $R = r = u_{xx}(X(\lambda), Y(\lambda))$, $S = s = u_{xy}(X(\lambda), Y(\lambda))$, $\tilde{S} = s = u_{yx}(X(\lambda), Y(\lambda))$ and $T = t = u_{yy}(X(\lambda), Y(\lambda))$. We assume $u$ is twice continuously differentiable, and therefore $S=u_{xy} = u_{yx} = \tilde{S}$. Henceforth strip conditions of second order for $u$ become
    \begin{subequations}
    \label{eqn:stripConditionPQ}    
    \begin{align}
    p_\lambda & = r X_\lambda + s Y_\lambda, \\ 
    q_\lambda & = s X_\lambda + t Y_\lambda. 
    \end{align}
    \end{subequations} 
    For completeness we give the strip conditions of third order, viz.,
    \begin{subequations}
    \label{eqn:C2StripConditions}
    \begin{alignat}{3}
    r_\lambda & = u_{xxx} X_\lambda + u_{xxy} Y_\lambda && = r_x X_\lambda + r_y Y_\lambda && = r_x X_\lambda + s_x Y_\lambda, \\ 
    s_\lambda & = u_{xyx} X_\lambda + u_{xyy} Y_\lambda && = s_x X_\lambda + s_y Y_\lambda && = s_x X_\lambda + t_x Y_\lambda = r_y X_\lambda + s_y Y_\lambda, \\
    t_\lambda & = u_{yyx} X_\lambda + u_{yyy} Y_\lambda && = t_x X_\lambda + t_y Y_\lambda && = s_y X_\lambda + t_y Y_\lambda.
    \end{alignat}
    \end{subequations}

    The process of finding higher order strips is called extending. To clarify, the curve $C_0 = \{(X(\lambda), Y(\lambda), U(\lambda)) | \lambda \in I\}$ (a strip of zeroth order) is extended to a strip $C_1$, given by~\eqref{eqn:C_1_def}. Similarily $C_1$ is extended to a strip of second order, given by
    \begin{align}
    \label{eqn:C_2_def}
    C_2 = \{(X(\lambda), Y(\lambda), U(\lambda), P(\lambda), Q(\lambda), R(\lambda), S(\lambda), \tilde{S}(\lambda), T(\lambda)) | \lambda \in I\}.
    \end{align}
    
    We define $C_2$ to be an integral strip if there exists a $C_1$-strip which can be extended to the $C_2$-strip uniquely, solely using the PDE~\eqref{eqn:nonlinearBaseEqn} and the strip conditions~\eqref{eqn:stripConditionPQ}. In this case $C_1$ is called a \textit{free strip}.
    If $C_1$ is not a free strip, additional requirements should be prescribed in order for $C_1$ to be extendable to an integral strip $C_2$. In this case we call $C_1$ a \textit{characteristic strip} which implies that not all second order derivatives of $u$ can be determined uniquely from $C_1$, the PDE~\eqref{eqn:nonlinearBaseEqn} and the strip conditions.
    To put into context, let $C_\text{b} = \{(x(\lambda), y(\lambda)) | \lambda \in I\}$ be a base curve, $z = u(x,y)$ be an integral surface of \eqref{eqn:nonlinearBaseEqn} and let $C_0$ be a corresponding zeroth order strip. Furthermore, we supplement $p = u_x$ and $q = u_y$ to obtain a first order strip $C_1$. If by using the PDE~\eqref{eqn:nonlinearBaseEqn} and the strip conditions~\eqref{eqn:stripConditionPQ} we are able to determine $r = u_{xx}, s=u_{xy}$ and $t = u_{yy}$ uniquely, then the strip $C_2 = \{(x(\lambda), y(\lambda), u(\lambda), p(\lambda), q(\lambda), r(\lambda), s(\lambda), t(\lambda))|\lambda \in I\}$ is called an integral strip, and $C_1$ is a free strip, otherwise $C_1$ is called a characteristic strip.
    
    Note that along a free strip, but not along a characteristic strip, the derivatives $u_{xx}, u_{xy}$ and $u_{yy}$ can all be determined along the strip, either by being interior derivatives with respect to $C_1$, or by combining the PDE~\eqref{eqn:nonlinearBaseEqn} with the remaining interior derivatives. To illustrate, given $u$, $p = u_x$ and $q=u_y$, on a vertical line segment, i.e., $X_\lambda = 0$, by differentiation with respect to $y$ one can obtain $u_{xy}$ and $u_{yy}$, and $u_{xx}$ follows from the PDE, as will be shown in Section~\ref{sec:BoundaryConditionsHMA}.
    
    For completeness, if $C_1$ is a characteristic strip, its carrier $C_0$ will be called a characteristic curve in $(x,y,z)$-space, and the base curve $C_\text{b}$, will be called a characteristic base curve. Generally we refer to a characteristic strip, characteristic curve and characteristic base curve simply as `the characteristic'.
    
    Note that thus far we considered an entire curve/strip to be either free or characteristic. Formally this should be evaluated pointwise, which introduces the notion of a characteristic base point, a characteristic point and a characteristic element for a $2$-, $3$- and $5$-dimensional point on $C_\text{b}$, $C_0$ and $C_1$, respectively. This distinction is often not necessary due to the fact that every strip, which has one point in common with the integral surface and all its tangent planes equal to that of the integral surface, lies entirely on said surface. To see this consider the strip $\gamma$ parameterized by $\lambda \in I$, given by $\gamma(\lambda) = (X(\lambda), Y(\lambda), U(\lambda), u_x(\lambda), u_y(\lambda))$ with strip condition $U_\lambda = u_x X_\lambda + u_y Y_\lambda$. Let $(x_0, y_0, u_0)$ lie on the integral surface $z = u(x,y)$, so $u_0 = u(x_0, y_0)$. Furthermore, let $\gamma$ pass through $(x_0, y_0, u_0)$, i.e., there exists a $\lambda_0$ such that $(x_0, y_0, u_0) = (X(\lambda_0), Y(\lambda_0), U(\lambda_0))$ and $u_0 = u(X_0, Y_0) = u(X(\lambda_0), Y(\lambda_0))$. Let $d(\lambda) = U(\lambda) - u(X(\lambda), Y(\lambda))$ be the pointwise signed vertical distance between the integral surface $z=u(x,y)$ and the strip $\gamma$ at the point $(X(\lambda), Y(\lambda), u(X(\lambda), Y(\lambda)))$. If $d \equiv 0$ then clearly $\gamma$ lies on the integral surface $z=u(x,y)$. Obviously it holds that $d(\lambda_0) = 0$. Furthermore, the change in the signed distance $d$ for $\lambda \in I$ can be found by
    \begin{align}
    \odiff{d}{\lambda} = \odiff{U}{\lambda} - u_x \odiff{X}{\lambda} - u_y \odiff{Y}{\lambda} = \odiff{U}{\lambda} - (u_x X_\lambda + u_y Y_\lambda) = \odiff{U}{\lambda} - U_\lambda = 0,
    \end{align}
    where we applied the strip condition. 
    Because $\odiff{d}{\lambda} = 0$ for all $\lambda \in I$ and $d(\lambda_0) = 0$, $d \equiv 0$, and hence the strip lies entirely on the integral surface $u$.

    \subsection{The characteristic condition}
    We will derive and discuss the conditions under which a strip $C_1$ is a characteristic strip in this section. These conditions will be called the characteristic conditions. We will impose conditions on $C_1$, based on our starting equation~\eqref{eqn:nonlinearBaseEqn}, such that (at least one of the) second and higher order derivatives cannot be determined uniquely. 
    We will discuss three different approaches to obtaining the characteristic conditions.
            
    \subsubsection{The characteristic condition by the implicit function \hfill\break\mbox{theorem}}
    Fundamentally, we are looking for conditions on the solvability for the second order derivatives $r, s$ and $t$. One way to derive the characteristic condition is to apply the implicit function theorem. To this end, let the $C_1$-strip be parameterized by $\lambda$ as before. Define
    \begin{align}
    \mathbf{f}(x,y,u,p,q,x_\lambda, y_\lambda, p_\lambda, q_\lambda \mid r,s,t) = 
    \begin{pmatrix}
    F(x,y,u,p,q,r,s,t) \\
    x_\lambda r + y_\lambda s - p_\lambda \\ 
    x_\lambda s + y_\lambda t - q_\lambda
    \end{pmatrix}.
    \end{align}
    The components of the vector-valued function $\mathbf{f}$ are formed by our PDE~\eqref{eqn:nonlinearBaseEqn} and the two strip conditions~\eqref{eqn:stripConditionPQ} for $p_\lambda$ and $q_\lambda$.
    The implicit function theorem \cite[p. 731]{Adams} states that if there exists a $\lambda_0$ such that
    \begin{align}
    \begin{split}
    \mathbf{f}(x(\lambda_0),y(\lambda_0),u(\lambda_0),p(\lambda_0),q(\lambda_0),x_\lambda(\lambda_0), y_\lambda(\lambda_0), \hspace{30pt} \\ p_\lambda(\lambda_0), q_\lambda(\lambda_0) \mid r(\lambda_0),s(\lambda_0),t(\lambda_0)) = \mathbf{0},
    \end{split}
    \end{align} and the Jacobi matrix
    \begin{align}
    \label{eqn:A_def}
    \mathbf{A} := \pdiff{\mathbf{f}}{(r,s,t)} = \begin{pmatrix}
    F_r & F_s & F_t \\
    x_\lambda & y_\lambda & 0 \\
    0 & x_\lambda & y_\lambda 
    \end{pmatrix},
    \end{align}
    is nonsingular, then there is an open set $\Lambda \subset \mathbb{R}$ containing $\lambda_0$ and a unique continuously differentiable function $\mathbf{g} : \Lambda \rightarrow \mathbb{R}^3$ with $\mathbf{g}(\lambda_0) = (r(\lambda_0), s(\lambda_0), t(\lambda_0))^\text{T}$ such that 
    \begin{align}
    \mathbf{f}(x(\lambda),y(\lambda),u(\lambda),p(\lambda),q(\lambda),x_\lambda(\lambda), y_\lambda(\lambda), p_\lambda(\lambda), q_\lambda(\lambda)|\mathbf{g}(\lambda)) = \mathbf{0}, \text{ for all } \lambda \in \Lambda.
    \end{align}
    If $D := \det\left(\mathbf{A}\right) \neq 0$, then $(r,s,t)$ can be found uniquely along $C_1$, i.e., we have a free strip. Alternatively, if $D = 0$, then $(r,s,t)$ cannot be determined uniquely, hence we have a characteristic strip. The case $D = 0$ is therefore called the characteristic condition and it can be written as
    \begin{align}
    \label{eqn:Q}
    D = F_r y_\lambda^2 - F_s x_\lambda y_\lambda + F_t x_\lambda^2 = 0.
    \end{align}
            
    \subsubsection{The characteristic condition by a coordinate \hfill\break\mbox{transformation}}
    \label{sec:coordinate_transformation}
    The characteristic condition can also be derived by means of a coordinate transformation \cite[p. 419]{HilbertCourant} This can be achieved due to the following equivalent definition of a characteristic.
    If the differential equation $F = 0$ represents an interior differential equation along a strip $C_1$, then $C_1$ is a characteristic strip. The term interior differential operator here means that along $C_1$ the second order differential operator $F$ can be expressed solely in terms of derivatives of $x$, $y$ and $u$ with respect to the parameter describing the base curve of $C_1$.
    
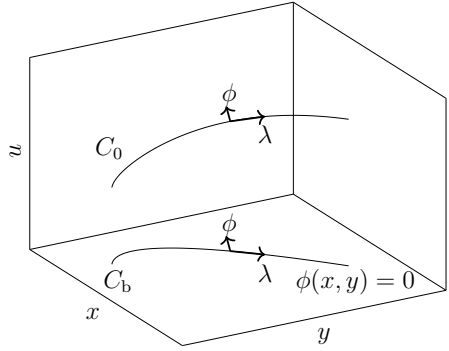
\begin{figure}[H]
        \centering
        \captionsetup{width=0.9\linewidth}
        \scalebox{0.8}{
        \begin{tikzpicture}
            \begin{axis}
            [
            ticks=none,
            xlabel={$x$},
            ylabel={$y$},
            zlabel={$u$},
            xmin=-.5,
            xmax=1.5,
            ymin=0.5,
            ymax=3.5,
            zmin=0,
            zmax=2.5,
            view={60}{30},
            ]

            \addplot3[
            domain=0:1,
            samples = 60,
            samples y=0,
            ]
            ({1.8 * x^2 - x},
            {1 + x^2 + x},
            {1 + x * sin(deg(2*x))});
            
            \draw[->, line width=0.3mm] (axis cs:0.182, 2.19, 1.6898) -- (axis cs:-0.068, 2.36, 1.6898) node [above, xshift=0.05cm, yshift=-0.15cm] {$\phi$};
            \draw[->, line width=0.3mm] (axis cs:0.182, 2.19, 1.6898) -- (axis cs:0.352, 2.44, 1.7997) node [below] {$\lambda$};

            \addplot3[
            domain=0:1,
            samples = 60,
            samples y=0,
            ]
            ({1.8 * x^2 - x},
            {1 + x^2 + x},
            {0});

            \draw[->, line width=0.3mm] (axis cs:0.182, 2.19, 0.0) -- (axis cs:-0.068, 2.36, 0.0) node [above, xshift=0.05cm, yshift=-0.15cm] {$\phi$};
            \draw[->, line width=0.3mm] (axis cs:0.182, 2.19, 0.0) -- (axis cs:0.352, 2.44, 0.0) node [below] {$\lambda$};
            
            \end{axis}
            
            \node (Cb) [label={$C_\text{b}$}] at (1.45, 0.6, 0) {};
            \node [label={$\phi(x,y) = 0$}] at (5.55, 0.75, 0.5){};
            \node[label={$C_0$}] at (1.5, 3, 0.5) {};
        \end{tikzpicture}
        }
        \caption{ Schematic representation of a curve $C_0$, with its base curve $C_\text{b}$ and the accompanying parameters along and perpendicular to the curves.}
        \label{fig:interiorDif}
    \end{figure}

    Let $C_1$ be the strip of interest with corresponding base curve $C_\text{b}$ as shown in Figure~\ref{fig:interiorDif}. 
    We introduce the coordinate transformation 
    \begin{align}
    \label{eqn:coordinateTrans}
    (x,y) \rightarrow (\phi(x,y), \lambda(x,y)),
    \end{align}
    where $\lambda$ is the parameter along the curve $C_\text{b}$ and $\phi$ leads away from $C_\text{b}$.
    \noindent Recall that we require $C_\text{b}$ to be a regular curve, i.e., $x_\lambda^2 + y_\lambda^2 \neq 0$. 
    Adopting the new coordinates the derivatives of $u$ are:
    \begin{subequations}
    \begin{align}
        u_x & = u_\phi \phi_x + u_\lambda \lambda_x, \\
        u_y & = u_\phi \phi_y + u_\lambda \lambda_y, \\
        \label{eqn:transform_xx}u_{xx} & = u_{\phi \phi} (\phi_x)^2 + 2 u_{\phi \lambda} \phi_x \lambda_x + u_{\lambda \lambda} (\lambda_x)^2 + u_\phi \phi_{xx} + u_\lambda \lambda_{xx}, \\
        u_{xy} & = u_{\phi \phi} \phi_x \phi_y + u_{\phi \lambda} (\phi_x \lambda_y  + \phi_y \lambda_x) + u_{\lambda \lambda} \lambda_x \lambda_y + u_\phi \phi_{xy} + u_\lambda \lambda_{xy}, \\
        \label{eqn:transform_yy}u_{yy} & = u_{\phi \phi} (\phi_y)^2 + 2 u_{\phi \lambda} \phi_y \lambda_y + u_{\lambda \lambda} (\lambda_y)^2 + u_\phi \phi_{yy} + u_\lambda \lambda_{yy}.
    \end{align}
    \end{subequations}
    Using these relations, one can construct a function $G$ such that
    \begin{align}
        \label{eqn:FGrelationCoordinateTransformation}
    F(x,y,u,u_x,u_y,u_{xx},u_{xy},u_{yy}) = G(\phi, \lambda, u, u_\phi, u_\lambda, u_{\phi \phi}, u_{\phi \lambda}, u_{\lambda \lambda}) = 0.
    \end{align}
    Recall that a $C_1$-strip is characteristic if not all higher order derivatives can be determined uniquely along the strip. As $p$ and $q$ are known on a $C_1$-strip, both $u_\phi$ and $u_\lambda$ can be obtained provided the Jacobian of the coordinate transformation~\eqref{eqn:coordinateTrans} $\phi_x \lambda_y - \phi_y \lambda_x \neq 0 $. Because $\lambda$ is the parameter along the strip, $(u_\phi)_\lambda$ and $(u_\lambda)_\lambda$ can naturally be found by differentiating along the strip. Therefore, for $C_1$ to be characteristic, $u_{\phi \phi}$ should be undetermined. 
    If $G_{u_{\phi\phi}} = 0$ then $u_{\phi \phi}$ cannot be determined from $G = 0$. Hence differentiating~\eqref{eqn:FGrelationCoordinateTransformation} and applying \eqref{eqn:transform_xx}-\eqref{eqn:transform_yy} we find
    \begin{align}
    \begin{split}
        G_{u_{\phi\phi}}(\phi, \lambda, u, u_\phi, u_\lambda, u_{\phi \phi}, u_{\phi \lambda}, u_{\lambda \lambda})
        & = F_{u_{xx}} \pdiff{u_{xx}}{u_{\phi\phi}} + F_{u_{xy}} \pdiff{u_{xy}}{u_{\phi\phi}} + F_{u_{yy}} \pdiff{u_{yy}}{u_{\phi\phi}} \\
        & = F_{r} \phi_x^2 + F_{s} \phi_x \phi_y + F_{t} \phi_y^2 \\
        & = 0,
    \end{split}
    \end{align}
    which is the characteristic condition. Assuming $\phi_y \neq 0$, this can further be rewritten as
    \begin{align}
    \label{eqn:charConditionPhi}
    F_{r} \left(\frac{\phi_x}{\phi_y}\right)^2 \hspace{-2pt} + F_{s} \frac{\phi_x}{\phi_y} + F_{t} = 0.
    \end{align}
    To see that this is equivalent to~\eqref{eqn:Q} consider the following:
    $\phi$ is constant along $C_\text{b}$, therefore for fixed $x,y$ it should not depend on $\lambda$, i.e., 
    \begin{align}
        \phi_\lambda = \phi_x x_\lambda + \phi_y y_\lambda = 0,
    \end{align}
    which is equivalent to 
    \begin{align}
        \frac{\phi_x}{\phi_y} = - \frac{y_\lambda}{x_\lambda}.
    \end{align}
    Substituting this in~\eqref{eqn:charConditionPhi} yields the previously found characteristic condition~\eqref{eqn:Q}.

\subsubsection{The characteristic condition for second order strips}
\label{sec:TheCharacteristicConditionByLinearAlgebra}
In the previous sections we established relations such that $C_1$ is uniquely defined, while $C_2$-strips are not. Although the characteristic condition $D = 0$ is fundamental, it yields no practical means to determine the evolution of the solution along a $C_1$-strip. In this section we show that the characteristic condition obtained also holds for second order strips, which does provide insights on how to determine the evolution of a $C_1$-strip. This evolution will be further discussed in Section~\ref{sec:CompatibilityConditions} and Section~\ref{sec:EvolutionAlongTheCurves}.

Not all second order derivatives can be determined along a characteristic $C_1$-strip. Likewise, not all third order derivatives can be determined uniquely either.
To determine relations for the derivatives of $r,s,t$ we apply strip conditions~\eqref{eqn:C2StripConditions}, together with additional relations, which are found by differentiating the PDE \eqref{eqn:nonlinearBaseEqn} with respect to $x$ and $y$, viz.,
\begin{subequations}
    \begin{align}
    \odiff{F}{x} & = F_x + F_u p + F_p r + F_q s + F_r r_x + F_s s_x + F_t t_x = 0, \\
    \odiff{F}{y} & = F_y + F_u q + F_p s + F_q t + F_r r_y + F_s s_y + F_t t_y = 0.
    \end{align}
\end{subequations}
Combining these with strip conditions~\eqref{eqn:C2StripConditions} with $X(\lambda) = x(\lambda)$ etc., yields two systems of equations 
\begin{align}
\label{eqn:StripSystem}
\mathbf{A} \begin{pmatrix}
r_x \\ s_x \\ t_x
\end{pmatrix}= \begin{pmatrix}
- F^x \\ r_\lambda \\ s_\lambda
\end{pmatrix}
,
\hspace{10pt}
\mathbf{A} \begin{pmatrix}
r_y \\ s_y \\ t_y
\end{pmatrix}= \begin{pmatrix}
- F^y \\ s_\lambda \\ t_\lambda
\end{pmatrix}
, 
\end{align}
where the equations are formed by collecting the $x$- and $y$-derivatives respectively, and where $F^x :=  F_x + F_u p + F_p r + F_q s$, $F^y = F_y + F_u q + F_p s + F_q t$ and $\mathbf{A}$ is given by~\eqref{eqn:A_def}. Because $F$ becomes an interior operator along a characteristic strip, $r_\lambda, s_\lambda$ and $t_\lambda$ can be determined along a $C_2$-strip while not all of $r_x, s_x, t_x, r_y, s_y$ and $t_y$ can. Therefore $\mathbf{A}$ should be singular. Hence we obtain once more the characteristic condition~\eqref{eqn:Q}.

    \subsection{Compatibility conditions}
    \label{sec:CompatibilityConditions}
    Because $D = \det(\mathbf{A}) = 0$, the systems in~\eqref{eqn:StripSystem} may not have solutions. In this section we will derive compatibility conditions such that solutions do exist. Consider the rank of $\mathbf{A}$.
    The rank of a matrix equals the order of the largest non-vanishing minor, which is known as the determinantal rank. The matrix $\mathbf{A}$ has 9 minors of order 2, the minors formed by the lower right and lower left $2\times2$ submatrices are
    \begin{align}
    M_{1,1} = \begin{vmatrix}
    y_\lambda & 0 \\
    x_\lambda & y_\lambda
    \end{vmatrix}
    =
    y_\lambda^2, \hspace{20pt}
    M_{1,3} = \begin{vmatrix}
    x_\lambda & y_\lambda \\
    0 & x_\lambda
    \end{vmatrix} = x_\lambda^2,
    \end{align}
    where $M_{ij}$ denotes the minor formed by deleting the $i$th row and $j$th column.
    Because $x_\lambda^2 + y_\lambda^2 \neq 0$, $M_{1,1}$ and $ M_{1,3} $ cannot be 0 simultaneously, hence $D = 0$ implies $\text{rank}(\mathbf{A}) = 2$. From \eqref{eqn:StripSystem} we conclude that solutions $(r_x, s_x, t_x)^\text{T}$ and $(r_y, s_y, t_y)^\text{T}$ do not always exist if $\mathbf{A}$ is singular. Therefore we require the left and right hand sides of \eqref{eqn:StripSystem} to be compatible, i.e., the vectors $(- F^x, r_\lambda, s_\lambda)^\text{T}$ and $(- F^y, s_\lambda, t_\lambda)^\text{T}$ should be in the column-space of any two columns of $\mathbf{A}$. Hence the matrix of any two column vectors of $\mathbf{A}$ with either $(- F^x, r_\lambda, s_\lambda)^\text{T}$ or $(- F^y, s_\lambda, t_\lambda)^\text{T}$ should be singular. To this end we introduce the matrices
    \begin{subequations}
        \label{eqn:augmentedMA}
        \begin{align}
        \mathbf{A}^x = 
        \left(
        \begin{array}{ccc|c}
        F_r & F_s & F_t & -F^x\\
        x_\lambda & y_\lambda & 0 & r_\lambda\\
        0 & x_\lambda & y_\lambda & s_\lambda
        \end{array}
        \right),
        \\
        \mathbf{A}^y = 
        \left(
        \begin{array}{ccc|c}
        F_r & F_s & F_t & -F^y\\
        x_\lambda & y_\lambda & 0 & s_\lambda\\
        0 & x_\lambda & y_\lambda & t_\lambda
        \end{array}
        \right).
        \end{align}
    \end{subequations}
    Let $D^x_{klm}$ denote the determinant formed by selecting the columns $k, l$ and $m$ of $\mathbf{A}^x$ and similarly, we introduce $D^y_{klm}$. We find by including the first, third and fourth column
    \begin{subequations}
        \label{eqn:Axy}
        \begin{align}
        D^x_{134} = 
        \begin{vmatrix}
        F_r & F_t & -F^x \\ x_\lambda & 0 & r_\lambda \\ 0 & y_\lambda & s_\lambda
        \end{vmatrix} & = - (F^x x_\lambda y_\lambda + F_r y_\lambda r_\lambda + F_t x_\lambda s_\lambda), 
        \\
        D^y_{134} = 
        \begin{vmatrix}
        F_r & F_t & -F^y \\ x_\lambda & 0 & s_\lambda \\ 0 & y_\lambda & t_\lambda
        \end{vmatrix} & = - (F^y x_\lambda y_\lambda + F_r y_\lambda s_\lambda + F_t x_\lambda t_\lambda),
        \end{align}
    \end{subequations}
    which both should equal zero due to $\mathbf{A}$ being singular.
    Trivially, since $D = 0$ similar calculations yield $D^x_{124} =  D^x_{134} \, x_\lambda / y_\lambda$, $D^x_{234} = D^x_{134}  \, y_\lambda / x_\lambda $, $D^y_{124} = D^y_{134}  \, x_\lambda / y_\lambda$, $D^y_{234} = D^y_{134}  \, y_\lambda / x_\lambda $ for $x_\lambda, y_\lambda \neq 0$.
    
    Condition~\eqref{eqn:Q} turns out to be enough to determine the evolution of $x, y, u, p$ and $q$ along the characteristics. Conditions~\eqref{eqn:Axy} are needed for the evolutions of other variables which are introduced in the next section.

    \subsection{Evolution along the characteristics}
    \label{sec:EvolutionAlongTheCurves}
    To derive the evolution of the solution along the characteristics we rewrite \eqref{eqn:Q} as a second order polynomial equation, viz.
    \begin{align}
    \label{eqn:Q2}
    D(\mu) = F_r \mu^2 - F_s \mu + F_t = 0,
    \end{align}
    where we introduced $\mu = y_\lambda / x_\lambda$ assuming $x_\lambda \neq 0$. In case $x_\lambda = 0$, we can use $\tilde{\mu} =  x_\lambda / y_\lambda$ instead. Solving \eqref{eqn:Q2} for $\mu$ yields two roots $a$ and $b$ viz.
    \begin{align}
    \label{eqn:abBaseDef}
    a = \frac{F_s + \sqrt{\Delta}}{2 F_r}, \hspace{15pt} b = \frac{F_s - \sqrt{\Delta}}{2 F_r},
    \end{align}
    where $\Delta = F_s^2 - 4 F_r F_t$ is the discriminant. The discriminant allows us to classify the differential operator $F$. If at a point $\mathbf{x}_0=(x_0,y_0) \in \Omega$ the discriminant \mbox{$\Delta(\mathbf{x}_0) > 0$}, then the PDE~\eqref{eqn:nonlinearBaseEqn} is hyperbolic in that point \cite[p. 420]{HilbertCourant}. Naturally there exists a (small) neighborhood of $\mathbf{x}_0$ for which the PDE is hyperbolic.
    Similarly we call the PDE parabolic in $\mathbf{x}_0$ if $\Delta(\mathbf{x}_0) = 0$ and elliptic if $\Delta(\mathbf{x}_0) < 0$.
    If for all $\mathbf{x}_0 \in \Omega$ we have $\Delta(\mathbf{x}_0)>0$, then $F$ is called hyperbolic, or hyperbolic in the entire domain.
    In the following we restrict ourselves to the hyperbolic case. By definition, we have two separate families of characteristic curves defined by $a(x(\lambda), y(\lambda)) x_\lambda = y_\lambda$ or $b(x(\lambda), y(\lambda)) x_\lambda = y_\lambda$, respectively, passing through the point $(x(\lambda), y(\lambda))$.

    We can express $F_r, F_s$ and $F_t$ in terms of $a, b$ and $\Delta$ using \eqref{eqn:abBaseDef} as
    \begin{align}
    \label{eqn:abEqs1}
    F_t = \frac{ab}{a-b} \sqrt{\Delta}, \hspace{15pt} F_s = \frac{a + b}{a - b} \sqrt{\Delta}, \hspace{15pt} F_r = \frac{1}{a - b} \sqrt{\Delta}.
    \end{align}
    Alternatively, we can express $a$ and $b$ in terms of $F_r, F_s$ and $F_t$ as 
    \begin{align}
    \label{eqn:abEqs2}
     a + b = \frac{F_s}{F_r}, \hspace{15pt} a b = \frac{F_t}{F_r}.
    \end{align}
    Using the definition of $\mu$ we find $y_\lambda / x_\lambda = a$ or $y_\lambda / x_\lambda = b$, implying we have two distinct families of characteristics, one induced by $a$, and the other induced by $b$. To distinguish the characteristics, we write $x = x(\alpha)$, $y = y(\alpha)$ and $x = x(\beta)$, $y = y(\beta)$ for the characteristic induced by $a$ and $b$, respectively. Henceforth $\alpha$ and $\beta$ effectively take over the role of $\lambda$. As such instead of $x_\lambda$ we write $\odiff{x}{\alpha} = x_\alpha$ for the derivative of $x$ with respect to $\alpha$, and similarly for the other variables and for differentiation with respect to $\beta$. 
    
    The matrix $\mathbf{A}$, given in \eqref{eqn:A_def} actually represents two distinct matrices, $\mathbf{A}^\alpha$ and $\mathbf{A}^\beta$, because the derivatives w.r.t. $\lambda$ can be associated with both $\alpha$ and $\beta$. Because further derivations for either characteristic is done analogously for the other, we will only treat the characteristic induced by $a$, the $\alpha$-characteristic, and postulate the results for the $\beta$-characteristic. Note that for fixed $(x_0,y_0) \in \Omega$, two characteristics pass through $(x_0, y_0)$, i.e., both the $\alpha$- and $\beta$-characteristic. 
    Because the matrix $\mathbf{A}^\alpha$ has rank 2, the rows are linearly dependent and therefore $\kappa^\alpha_1$ and $\kappa^\alpha_2$ exist such that
    \begin{align}
    \label{eqn:FrFsFtKappaRelation}
    \begin{pmatrix}
    F_r \\ F_s \\ F_t 
    \end{pmatrix} = 
    \kappa^\alpha_1
    \begin{pmatrix}
    x_\alpha \\ y_\alpha \\ 0
    \end{pmatrix}
    + 
    \kappa^\alpha_2
    \begin{pmatrix}
    0 \\ x_\alpha \\ y_\alpha
    \end{pmatrix}.
    \end{align}
    The first row gives $x_\alpha = F_r /\kappa^\alpha_1$. By definition we have $y_\alpha = a x_\alpha$ and hence $y_\alpha = a F_r/\kappa^\alpha_1$. The third row then yields $\kappa^\alpha_2 = F_t/y_\alpha = \kappa^\alpha_1 F_t / (a F_r)$ which yields $ \kappa^\alpha_2 = b  \kappa^\alpha_1$ by \eqref{eqn:abEqs2}. For the sake of brevity we write $\kappa^\alpha = \kappa^\alpha_1$. Then~\eqref{eqn:FrFsFtKappaRelation} reduces to 
    \begin{subequations}
    \label{eqn:kappaRelations}
    \begin{align}
    F_r & = \kappa^\alpha x_\alpha, \\
    F_s & = \kappa^\alpha y_\alpha + b \kappa^\alpha x_\alpha, \\
    F_t & = b \kappa^\alpha y_\alpha.
    \end{align}
    \end{subequations}
    The evolution of $x, y, u, p$ and $q$ along the characteristics can be determined from~\eqref{eqn:kappaRelations}, the strip conditions \eqref{eqn:stripConditionU} and \eqref{eqn:stripConditionPQ}, respectively, giving
    \begin{subequations}
    \label{eqn:evolutionXYUPQ}
    \begin{align}
    x_\alpha & = \frac{F_r}{\kappa^\alpha}, \\
    y_\alpha & = a \frac{F_r}{\kappa^\alpha}, \\
    u_\alpha & = (p + a q) \frac{F_r}{\kappa^\alpha}, \\
    p_\alpha & = (r + a s) \frac{F_r}{\kappa^\alpha}, \\
    q_\alpha & = (s + a t) \frac{F_r}{\kappa^\alpha}, 
    \end{align}
    \end{subequations}
    where the choice of $\kappa^\alpha$ determines the parametric scaling of the base curve. The evolution of $r,s$ and $t$ can be obtained using the compatibility conditions~\eqref{eqn:Axy}. To that purpose we rewrite~\eqref{eqn:Axy} as the underdetermined linear system
    \begin{align}
    \label{eqn:CharcondRST}
    \begin{pmatrix}
    F_r y_\alpha & F_t x_\alpha & 0 \\
    0 & F_r y_\alpha & F_t x_\alpha
    \end{pmatrix}
    \begin{pmatrix}
    r_\alpha \\ s_\alpha \\ t_\alpha
    \end{pmatrix}
    = - x_\alpha y_\alpha \begin{pmatrix}
    F^x \\ F^y
    \end{pmatrix}.
    \end{align}
    By the rank-nullity theorem \cite[p. 175]{LinAlgebra} the general solution of \eqref{eqn:CharcondRST} reads
    \begin{align}
    \begin{pmatrix}
    r_\alpha \\ s_\alpha \\ t_\alpha
    \end{pmatrix}
    = 
    - \begin{pmatrix}
    \frac{F^x x_\alpha}{F_r} \\ 0 \\ \frac{F^y y_\alpha}{F_t}
    \end{pmatrix}
    + \theta^\alpha
    \begin{pmatrix}
    \frac{F_t x_\alpha}{F_r y_\alpha} \\ -1 \\ \frac{F_r y_\alpha}{F_t x_\alpha}
    \end{pmatrix},
    \end{align}
    where the first terms is the particular solution with $s_\alpha = 0$ and the second terms is an element of the null space of the matrix for arbitrary $\theta^\alpha$. Rewriting this, using \eqref{eqn:kappaRelations}, yields
    \begin{align}
    \begin{pmatrix}
    r_\alpha \\ s_\alpha \\ t_\alpha
    \end{pmatrix}
    = 
    - \frac{1}{\kappa^\alpha} \begin{pmatrix}
    F^x \\ 0 \\ \frac{F^y}{b}
    \end{pmatrix}
    + \theta^\alpha
    \begin{pmatrix}
    b \\ -1 \\ \frac{1}{b}
    \end{pmatrix}.
    \end{align}
    Because the ODE system~\eqref{eqn:evolutionXYUPQ} depends on $a$ and $b$, what remains is to determine the evolution of $a$ and $b$ along the characteristics. These are straightforwardly calculated by taking the derivative of~\eqref{eqn:abBaseDef} w.r.t. $\alpha$. Treating the $\beta$-characteristic analogously to the $\alpha$-characteristic, we similarly obtain $\kappa^\beta$ and $\theta^\beta$, and the ODE systems read
    \begin{align}
    \label{eqn:fullGeneralODE}
     \begin{aligned}[c]
     x_\alpha & = \frac{F_r}{\kappa^\alpha}, \\
     y_\alpha & = a \frac{F_r}{\kappa^\alpha}, \\
     u_\alpha & = (p + a q) \frac{F_r}{\kappa^\alpha}, \\
     p_\alpha & = (r + a s) \frac{F_r}{\kappa^\alpha}, \\
     q_\alpha & = (s + a t) \frac{F_r}{\kappa^\alpha}, \\
     r_\alpha & = \theta^\alpha b - \frac{1}{\kappa^\alpha} F^x, \\
     s_\alpha & = - \theta^\alpha, \\
     t_\alpha & = \frac{1}{b} \left( \theta^\alpha - \frac{1}{\kappa^\alpha} F^y \right), \\
     a_\alpha & = \left(\frac{F_s + \sqrt{\Delta}}{2 F_r}\right)_\alpha, \\
     b_\alpha & = \left(\frac{F_s - \sqrt{\Delta}}{2 F_r}\right)_\alpha,
     \end{aligned}
     \qquad
     \begin{aligned}[c]
      x_\beta & = \frac{F_r}{\kappa^\beta}, \\
      y_\beta & = b \frac{F_r}{\kappa^\beta}, \\
      u_\beta & = (p + b q) \frac{F_r}{\kappa^\beta}, \\
      p_\beta & = (r + b s) \frac{F_r}{\kappa^\beta}, \\
      q_\beta & = (s + b t) \frac{F_r}{\kappa^\beta}, \\
      r_\beta & = \theta^\beta a - \frac{1}{\kappa^\beta} F^x, \\
      s_\beta & = - \theta^\beta, \\
      t_\beta & = \frac{1}{a} \left( \theta^\beta - \frac{1}{a \kappa^\beta} F^y \right), \\
      a_\beta & = \left(\frac{F_s + \sqrt{\Delta}}{2 F_r}\right)_\beta,\\
      b_\beta & = \left(\frac{F_s - \sqrt{\Delta}}{2 F_r}\right)_\beta,\\
      \end{aligned}
    \end{align}
    where the further evaluation of $a_\alpha, b_\alpha, a_\beta$ and $b_\beta$ yield no meaningful insight. Note the coupling between the two ODE systems, for example the evolution of $r_\alpha$ depends on $b$, which forms the direction of the other characteristic via $b x_\beta = y_\beta$. 
    All expressions for the evolution of $a$ and $b$ can be expanded, but neither can be fully expressed in $x, y, u$ and the derivatives of $u$. Furthermore, as we will see for the \MAe{}, the evolution for either $a$ or $b$, does depend both on $a$ and $b$. Hence,~\eqref{eqn:fullGeneralODE} is a mutually coupled system of ODEs.
%
%
    
    \section{The hyperbolic \MAe{}}
    \label{sec:HMAE}
    In the remaining of the paper we will consider the hyperbolic \MAe{}. Recall that the hyperbolic \MAe{} is given by
    \begin{align}
        \label{eqn:MASimple}
        F(x,y,u,p,q,r,s,t) = r t - s^2+ f^2 = 0 \hspace{5pt} \text{for} \hspace{5pt} (x,y) \in \Omega,
    \end{align}
    for $\Omega \subseteq \mathbb{R}^2$, the unknown function $u = u(x,y) \in C^3(\Omega)$ and the known function $f = f(x,y) \in C^1(\Omega)$, with $f \neq 0$ on $\bar{\Omega}$. The derivatives of $F$ are
    \begin{equation}
        \label{eqn:introductionVariables_F_}
        \begin{aligned}
            F_r  = t, \hspace{15pt}
            F_s = - 2 s, \hspace{15pt}
            F_t = r.
        \end{aligned}
    \end{equation}
    Furthermore, the characteristic equation is given by
    \begin{align}
        D(\mu) = F_r \mu^2 - F_s \mu + F_t = t \mu^2 + 2 s \mu + r = 0,
    \end{align}
    and the corresponding discriminant is
    \begin{equation}
    \label{eqn:introductionVariables_discriminant}
    \begin{aligned}
         \Delta & = F_s^2 - 4 F_r F_t = 4 s^2 - 4 r t = 4 f^2,
     \end{aligned}
     \end{equation}
     which is positive, and hence~\eqref{eqn:MASimple} is hyperbolic. The two real and distinct roots of~\eqref{eqn:introductionVariables_discriminant} are given by 
    \begin{equation}
    \label{eqn:introductionVariables_ab}
    \begin{aligned}     
              a = \frac{F_s + \sqrt{\Delta}}{2 F_r} = \frac{-s + f}{t}, &&\hspace{20pt}
              b = \frac{F_s - \sqrt{\Delta}}{2 F_r} = \frac{-s -f}{t}.
   \end{aligned}
   \end{equation}
   for $t \neq 0$. In case $t = 0$, we express $a$ and $b$ as
   \begin{equation}
   \begin{aligned}
    a = \frac{r}{-s - f}, &&\hspace{20pt} b = \frac{r}{-s + f},
   \end{aligned}
  \end{equation}
   instead. Furthermore, the auxiliary functions $F^x$ and $F^y$ are given by
    \begin{equation}
    \label{eqn:introductionVariables_FXFY}
    \begin{aligned} 
            F^x & = F_x + F_u p + F_p r + F_q s = F_x = 2 f f_x, \\
            F^y & = F_y + F_u q + F_p s + F_q t = F_y = 2 f f_y.
        \end{aligned}
    \end{equation}
    From~\eqref{eqn:MASimple} and~\eqref{eqn:introductionVariables_ab}, it follows that
    \begin{align}
    \label{eqn:equivalenceRelationsAB}
    a + b = \frac{-2 s}{t}, && a-b = \frac{2f}{t}, && a b = \frac{s^2 - f^2}{t^2} = \frac{r}{t}. 
    \end{align}
    From these relations we can express the second derivatives in terms of $a,b$ and $f$ as follows
    \begin{align}
        \label{eqn:equivalenceRelations}
        r = \frac{2 ab}{a-b} f, && s = -\frac{a+b}{a-b} f,  && t = \frac{2f}{a - b}.
    \end{align}
    The systems of ODEs \eqref{eqn:fullGeneralODE} then read
    \begin{equation}
        \label{eqn:evolutionFull}
        \begin{aligned}[c]
            x_\alpha & = \frac{t}{\kappa^\alpha}, \\
            y_\alpha & = a \frac{t}{\kappa^\alpha}, \\
            u_\alpha & = (p + a q) \frac{t}{\kappa^\alpha}, \\
            p_\alpha & = (r + a s) \frac{t}{\kappa^\alpha}, \\
            q_\alpha & = (s + a t) \frac{t}{\kappa^\alpha}, \\
            r_\alpha & = \theta^\alpha b - \frac{1}{\kappa^\alpha} \left(f^2\right)_x , \\
            s_\alpha & = - \theta^\alpha, \\
            t_\alpha & = \frac{1}{b} \Big(\theta^\alpha - \frac{1}{\kappa^\alpha} \left(f^2\right)_y \Big), \\
        \end{aligned}
        \qquad
        \begin{aligned}[c]
            x_\beta & = \frac{t}{\kappa^\beta}, \\
            y_\beta & = b \frac{t}{\kappa^\beta}, \\
            u_\beta & = (p + b q) \frac{t}{\kappa^\beta}, \\
            p_\beta & = (r + b s) \frac{t}{\kappa^\beta}, \\
            q_\beta & = (s + b t) \frac{t}{\kappa^\beta}, \\
            r_\beta & = \theta^\beta a - \frac{1}{\kappa^\beta} \left(f^2\right)_x , \\
            s_\beta & = - \theta^\beta, \\
            t_\beta & = \frac{1}{a} \Big(\theta^\beta - \frac{1}{\kappa^\beta} \left(f^2\right)_y \Big) . \\
        \end{aligned}
    \end{equation}
    Note that the ODE systems in~\eqref{eqn:evolutionFull} contain four parameters, viz. $\kappa^{\alpha}, \kappa^{\beta}$, which are determined by an appropriate scaling, and $\theta^{\alpha}$ and $\theta^{\beta}$, which are free parameters. Consequently, the derivatives of $r$, $s$ and $t$ cannot be rewritten such that they no longer depend on $\theta^\alpha$ as this would uniquely determine $r,s,t$ along the characteristic strip which contradicts the definition of a characteristic strip.
    By differentiating the expressions for $a$ and $b$ and using~\eqref{eqn:evolutionFull} we find the evolution along the $\alpha$-characteristic, viz.
    \begin{subequations}
    \begin{align}
    \begin{split}
    a_\alpha & 
    = \frac{1}{t} \left(1 - \frac{a}{b}\right) \theta^\alpha + \frac{1}{\kappa^\alpha} \left(f_x +a f_y\right) + \frac{a}{b} \frac{1}{\kappa^\alpha t} \left(f^2\right)_y\\
    &
    = \frac{1}{t} \left(1 - \frac{a}{b}\right) \theta^\alpha + \frac{1}{\kappa^\alpha} \left(f_x + \frac{a^2}{b} f_y\right),
    \end{split}
    \\
    \begin{split}
    b_\alpha &
    = \frac{\theta^\alpha - f_x \frac{t}{\kappa^\alpha} - f_y a \frac{t}{\kappa^\alpha}}{t} - \frac{b}{t} \Big(\theta^\alpha - \frac{2 f f_y}{\kappa^\alpha} \Big) \frac{1}{b} \\
    & 
    = -\frac{1}{\kappa^\alpha}(f_x + b f_y).
    \end{split}
    \end{align}
    \end{subequations}
    The expression for $a_\alpha$ could be rewritten, for example by using~\eqref{eqn:equivalenceRelations}, but due to $a\neq b$, it will include the unknown $\theta^\alpha$. More fundamentally, we cannot determine $a_\alpha$ explicitly as then both $a$ and $b$ can be uniquely determined along the $\alpha$-characteristic, from which $r,s,t$ would follow by \eqref{eqn:equivalenceRelations}, which contradicts the definition of a characteristic.
    
    Note that we are free to choose $\kappa^\alpha$ and $\kappa^\beta$ due to the freedom in parameterization of the base curve. In the following we conveniently choose $\kappa^\alpha = \kappa^\beta = F_r = t$. 
    Using relations \eqref{eqn:equivalenceRelations} the ODE system \eqref{eqn:evolutionFull} reduces to 
    \begin{equation}
        \label{eqn:evolutionReduced}
        \begin{aligned}[c]
            x_\alpha & = 1, \\
            y_\alpha & = a , \\
            u_\alpha & = p + a q, \\
            p_\alpha & = - a f, \\
            q_\alpha & = f, \\
            r_\alpha & = \theta^\alpha b - (a-b) f_x, \\
            s_\alpha & = - \theta^\alpha, \\
            t_\alpha & = \frac{1}{b} \Big(\theta^\alpha - (a - b) f_y \Big), \\
            a_\alpha & = \frac{a - b}{2 f} \Big[ \left( 1 - \frac{a}{b} \right) \theta^\alpha + f_x + \frac{a^2}{b} f_y\Big],\\
            b_\alpha & = \frac{b-a}{2f} (f_x + b f_y) ,
        \end{aligned}
        \qquad
        \begin{aligned}[c]
            x_\beta & = 1, \\
            y_\beta & = b, \\
            u_\beta & = p + b q, \\
            p_\beta & = b f, \\
            q_\beta & = -f, \\
            r_\beta & = \theta^\beta a - (a-b) f_x, \\
            s_\beta & = - \theta^\beta, \\
            t_\beta & = \frac{1}{a} \Big(\theta^\beta - (a-b) f_y \Big) , \\
            a_\beta & = \frac{a-b}{2 f} (f_x + a f_y), \\
            b_\beta & = -\frac{a - b}{2 f} \Big[ \left( \frac{b}{a} - 1 \right) \theta^\beta + f_x + \frac{b^2}{a} f_y \Big].
        \end{aligned}
    \end{equation}
    There is a lot of redundancy in these equations which directly follows from \eqref{eqn:introductionVariables_ab}, \eqref{eqn:equivalenceRelationsAB} and \eqref{eqn:equivalenceRelations}. We therefore reduce~\eqref{eqn:evolutionReduced} by omitting the equations involving $\theta_\alpha$ and $\theta_\beta$. What remains are the ODE systems
        \begin{equation}
            \label{eqn:evolutionFullyReduced}
            \begin{aligned}[c]
                x_\alpha & = 1, \\
                y_\alpha & = a , \\
                u_\alpha & = p + a q,\\
                p_\alpha & = - a f, \\
                q_\alpha & = f, \\
                b_\alpha & = \frac{b-a}{2f} (f_x + b f_y),
            \end{aligned}
            \qquad
            \begin{aligned}[c]
                x_\beta & = 1, \\
                y_\beta & = b, \\
                u_\beta & = p + b q,\\
                p_\beta & = b f, \\
                q_\beta & = -f, \\
                a_\beta & = \frac{a-b}{2 f} (f_x + a f_y),
            \end{aligned}
        \end{equation}
    which we integrate numerically. Recall that: ``every strip, which has one point in common with the integral surface and all its tangent planes equal to that of the integral surface, lies entirely on said surface.'', which furthermore justifies the reduction of~\eqref{eqn:evolutionReduced} to~\eqref{eqn:evolutionFullyReduced}, as only $u, p, q$ and the characteristics, so $u, p, q, a$ and $b$, need to be known.
    
    \subsection{Boundary conditions} 
    \label{sec:BoundaryConditionsHMA}
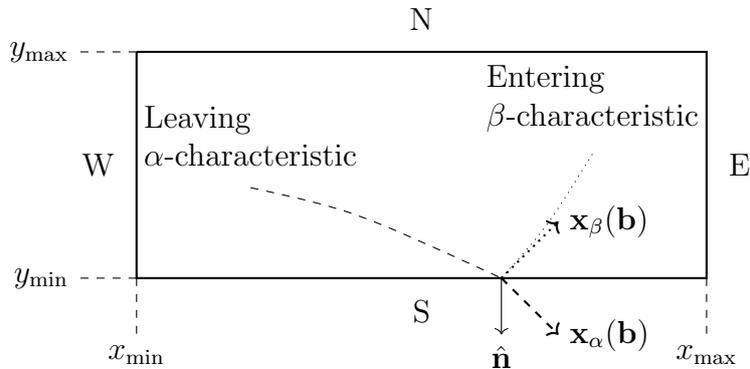
\begin{figure}[H]
    \centering
    \captionsetup{width=0.9\linewidth}
    \begin{tikzpicture}[scale=1.5]
    \coordinate (xi) at (3.2,0);
    \coordinate (alphaStart) at (1,0.8);
    \coordinate (betaEnd) at (4,1.1);

    \draw[dashed] (alphaStart) .. controls (1.9, 0.6) .. (xi); 
    \draw[dotted] (xi) .. controls (3.6, 0.45) .. (betaEnd); 
    \node[label={[align=left]Leaving\\$\alpha$-characteristic}] at (alphaStart) {};
    \node[label={[align=left]Entering\\$\beta$-characteristic}] at (betaEnd) {};
    
    \draw [-,thick] (0,2) |- (0,0) |- (5,0) |- (5,2) |- (0,2);
    
    \node[above] at (2.5, 2.1) {N};
    \node[left ] at (-.1, 1) {W};
    \node[below] at (2.5, -0.1) {S};
    \node[right] at (5.1, 1) {E};
    
    \draw [dashed] (0,0) |- (-1/2,0) node (ymin) [left] {$y_\text{min}$};
    \draw [dashed] (0,2) |- (-1/2,2) node (ymax) [left] {$y_\text{max}$};
    \draw [dashed] (0,0) -- (0,-1/2) node (xmin) [below] {$x_\text{min}$};
    \draw [dashed] (5,0) -- (5,-1/2) node (xmax) [below] {$x_\text{max}$};
    \draw [->] (xi) -- ($(xi)-(0,1/2)$) node (xilabelNode) [below] {$\hat{\mathbf{n}}$};
    
    \draw[dashed,->, thick] (xi) -- ($(xi)+(1/2,-1/2)$) node (boundaryLabel) [right] {$\mathbf{x}_\alpha(\mathbf{b})$}; 
    \draw[dotted,->, thick] (xi) -- ($(xi)+(1/2,1/2)$) node (boundaryLabel) [right] {$\mathbf{x}_\beta(\mathbf{b})$}; 
    
    \end{tikzpicture}
    \caption{ Schematic overview of a rectangular domain where the $\alpha$-characteristic leaves, and the $\beta$-characteristic enters the domain.}
    \label{fig:1boundaryChars}
    \end{figure}

    We solve~\eqref{eqn:MASimple} on a rectangular domain $\Omega = [x_{\text{min}}, x_{\text{max}}] \times [y_{\text{min}}, y_{\text{max}}]$, for $x_{\text{min}}, x_{\text{max}}$, $y_{\text{min}}, y_{\text{max}} \in \mathbb{R}$, with $x_{\text{min}} < x_{\text{max}}, y_{\text{min}} < y_{\text{max}}$ and we call $\{x_{\text{min}}\} \times [y_{\text{min}}, y_{\text{max}}]$ the initial base curve. We extend the initial base curve to an initial $C_1$-strip by supplementing it with $u_\text{W}, p_\text{W}: [y_{\text{min}}, y_{\text{max}}] \rightarrow\mathbb{R}$, where we prescribe $u(x_{\text{min}}, y) = u_\text{W}(y)$ and $ p(x_{\text{min}}, y) = p_\text{W}(y)$ for some $u_\text{W}$, $p_\text{W}$. The subscript `W' is used as the values are prescribed on the Western part of $\partial \Omega$, see Figure~\ref{fig:1boundaryChars}.
    
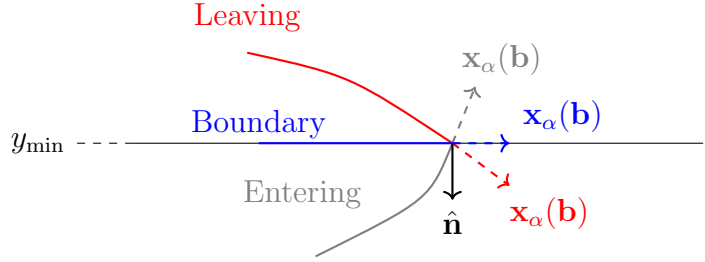
\begin{figure}[H]
        \centering
        \captionsetup{width=0.9\linewidth}
        \begin{tikzpicture}[scale=1.5]

        \draw [-] (0,0) |- (5,0);
        
        \coordinate (xi) at (2.8,0);
        \coordinate (alphaStart) at (1,0.8);
        \coordinate (betaEnd) at (1.6,-1.0);
        
        \draw[red, thick] (alphaStart) .. controls (1.9, 0.6) .. (xi); 
        \draw[dashed,->, red, thick] (xi) -- ($(xi)+(1/2,-3/8)$) node (leavingLabel) [right, below, xshift=15pt] {$\mathbf{x}_\alpha(\mathbf{b})$}; 
                
        \draw[gray, thick] (xi) .. controls (2.6, -0.5) .. (betaEnd); 
        \draw[dashed,->, gray, thick] (xi) -- ($(xi)+(1/5,1/2)$) node (enteringLabel) [right, above, xshift=10pt] {$\mathbf{x}_\alpha(\mathbf{b})$}; 
        
        \draw[blue, thick] (xi) |- ($(xi)-(1.7,0)$);
        \draw[dashed,->, blue, thick] (xi) -- ($(xi)+(1/2,0)$) node (boundaryLabel) [right, above, xshift=20pt] {$\mathbf{x}_\alpha(\mathbf{b})$}; 
        
        \node[label={[red, align=left]Leaving}] at (alphaStart) {};
        \node[xshift=-5pt, yshift=10pt, label={[gray, align=left]Entering}] at (betaEnd) {};
        \node[label={[blue, align=left]Boundary}] at ($(xi)-(1.7,0.15)$) {};
        
        \draw [dashed] (0,0) |- (-1/2,0) node (ymin) [left] {$y_\text{min}$};
        \draw [thick, ->] (xi) -- ($(xi)-(0,1/2)$) node (xilabelNode) [below] {$\hat{\mathbf{n}}$};
        
        \end{tikzpicture}
        \caption{ Schematic classification of characteristics. }
        \label{fig:classificationCharacteristics}
    \end{figure}

    If the initial strip is a free strip, then prescribing $u_\text{W}$, $p_\text{W}$ uniquely determines the $C_2$-strip as an extension of the initial base curve. To verify this, we check whether the characteristic condition~\eqref{eqn:Q2} holds. Therefore we parameterize the initial base curve as $x(\lambda) = 0$, $y(\lambda) = \lambda$, $\lambda \in [0,1]$, then $x_\lambda = 0$, $y_\lambda = 1$ and the characteristic condition yields $t y_\lambda^2 + 2 s x_\lambda y_\lambda + r x_\lambda^2 = t = 0$.
    Hence if $t\neq 0$ on the initial strip, i.e., if $u_\text{W}^{\prime\prime}(y)\neq 0$, the initial strip is a free strip. Henceforth we assume $u_\text{W}^{\prime\prime}(y) \neq0$, which then implies the initial base curve uniquely extends to a $C_2$-strip, and we can uniquely determine $q, a, b, r, s, t$ on the initial strip via
    \begin{subequations}
        \label{eqn:initialStripCalculations}
    \begin{align}
    q(x_{\text{min}}, y) & = u_\text{W}^{\prime}(y), \\
    t(x_{\text{min}}, y) & = u_\text{W}^{\prime\prime}(y), \\
    s(x_{\text{min}}, y) & = p_\text{W}^\prime(y), \\
    a(x_{\text{min}}, y) & = \frac{-s(x_{\text{min}}, y) + f_\text{W}(y)}{t(x_{\text{min}}, y)}, \\
    b(x_{\text{min}}, y) & = \frac{-s(x_{\text{min}}, y) - f_\text{W}(y)}{t(x_{\text{min}}, y)}, \\
    r(x_{\text{min}}, y) & = \frac{2 a(x_{\text{min}}, y) b(x_{\text{min}}, y)}{a(x_{\text{min}}, y) -b(x_{\text{min}}, y)}f_\text{W}(y),
    \end{align}
    \end{subequations}
    where $f_\text{W}(y) := f(x_{\text{min}}, y)$.
    For the lower and upper boundary, i.e., for $y \in \{y_{\text{min}}, y_{\text{max}}\}$, the required boundary conditions are more delicate. 
    To understand this let $\mathbf{x}_\alpha = (x_\alpha, y_\alpha)$ and $\mathbf{x}_\beta = (x_\beta, y_\beta)$ denote the tangent vectors of the characteristics. Let $\mathbf{b} \in \partial \Omega$ and $\mathbf{\hat{n}}$ the outward unit normal vector on the boundary. We classify the $\alpha$-characteristics, and similarly the $\beta$-characteristic, based on whether they are entering or leaving the domain as follows
    \begin{itemize}
        \item Leaving characteristic if $\mathbf{x}_\alpha(\mathbf{b}) \boldsymbol{\cdot} \hat{\mathbf{n}} > 0$,
        \item Entering characteristic if $\mathbf{x}_\alpha(\mathbf{b}) \boldsymbol{\cdot} \hat{\mathbf{n}} < 0$,
        \item Boundary characteristic if $\mathbf{x}_\alpha(\mathbf{b}) \boldsymbol{\cdot} \hat{\mathbf{n}} = 0$,
    \end{itemize}
    which is schematically shown in Figure~\ref{fig:classificationCharacteristics} for $y = y_\text{min}$.
    
    We assume $a(x,y)$ and $b(x,y)$ to be well defined for all $(x,y) \in \Omega$ and the \MAe{} to be hyperbolic, which implies there are two characteristics passing through each point $(x,y) \in \Omega$. By classifying the characteristics as entering or leaving, one can determine if and how many additional boundary conditions need to be prescribed at each boundary point. Boundary characteristics should be treated as leaving characteristics.
    We distinguish three cases: one characteristic entering, two characteristics entering and zero characteristics entering in the point $\mathbf{b}$ on the boundary. 
        
    \textit{Case 1. One characteristic leaving and one entering the domain, either}
    \begin{align}
        \mathbf{x}_\alpha \boldsymbol{\cdot} \hat{\mathbf{n}} > 0, \hspace{5pt} \mathbf{x}_\beta \boldsymbol{\cdot} \hat{\mathbf{n}} < 0 \hspace{10pt} \text{or}  \hspace{10pt}     \mathbf{x}_\alpha \boldsymbol{\cdot} \hat{\mathbf{n}} < 0, \hspace{5pt}  \mathbf{x}_\beta \boldsymbol{\cdot} \hat{\mathbf{n}} > 0.
    \end{align}
    An example of this case is shown in Figure \ref{fig:1boundaryChars} at $y = y_\text{min}$, with $N, W, S, E$ denoting the Western, Northern, Southern and Eastern boundary segments of the domain. The curves denote one $\alpha$- (dashed) and one $\beta$-characteristic (dotted). 
    At $\mathbf{b}$, the values $u$, $p$, $q$ and $b$ can be computed from the ODE system for the $\alpha$-characteristic. However, $a$ cannot be determined because the evolution of $a$ along the $\alpha$-characteristic is unknown. Therefore we should impose one boundary condition, which is the initial condition for the entering $\beta$-characteristic, such that $a$ can be computed. We can either prescribe $a$ directly, or prescribe either $r, s$ or $t$ and compute $a$ from inverting \eqref{eqn:equivalenceRelations}, viz.
    \begin{align}
    \label{eqn:aAsFunctionsOfb}
    & a = \frac{b r}{r - 2 b f}, && a = \frac{s-f}{s+f} b, && a = b + \frac{2 f}{t} &
    \end{align} 
    The remaining two unknowns of $r, s$ and $t$ then follow from~\eqref{eqn:equivalenceRelations}.
    
\begin{figure}[H]
        \centering
        \captionsetup{width=0.8\linewidth}
        \begin{tikzpicture}[scale=1.5]
        \coordinate (xi) at (1.5,0);
        \coordinate (alphaNaturalStart) at (1,0.6);
        \coordinate (alphaNaturalEnd) at (2.8,1.2);
        \coordinate (betaNaturalStart) at (1,0);
        \coordinate (betaNaturalEnd) at (2.5,1.2);
        \coordinate (collisionPoint) at ($(xi)-(0.3, 0)$);
        \coordinate (naturalBetaLabel) at ($(collisionPoint) - (0.6, 0.7)$);
        \coordinate (alphaEnd) at (3.45,1.25);
        \coordinate (betaEnd) at (2.6,1.1);


        \draw[dashed] (xi) .. controls (2.7, 0.5) .. (alphaEnd); 
        \draw[dotted] (xi) -- (betaEnd); 
        \node[yshift=0.4cm, xshift=0.9cm, label={[align=left]\small Entering\\ $\alpha$-characteristic}] at ($(betaEnd)-(-0.1,0.17)$) {};
        \node[yshift=-0.5cm, xshift=-2cm, label={[align=left]\small Entering\\ $\beta$-characteristic}] at ($(betaEnd)-(-0.1,0.17)$) {};
    
        \draw [-,thick] (0,2) |- (0,0) |- (5,0) |- (5,2) |- (0,2);
        
        \node[above] at (2.5, 2.1) {N};
        \node[left ] at (-.1, 1) {W};
        \node[below] at (2.5, -0.1) {S};
        \node[right] at (5.1, 1) {E};
        
        \draw [dashed] (0,0) |- (-1/2,0) node (ymin) [left] {$y_\text{min}$};
        \draw [dashed] (0,2) |- (-1/2,2) node (ymax) [left] {$y_\text{max}$};
        \draw [dashed] (0,0) -- (0,-1/2) node (xmin) [below] {$x_\text{min}$};
        \draw [dashed] (5,0) -- (5,-1/2) node (xmax) [below] {$x_\text{max}$};
        \draw [->] (xi) -- ($(xi)-(0,1/2)$) node (xilabelNode) [below] {$\hat{\mathbf{n}}$};

        \draw[dashed,->, thick] (xi) -- ($(xi)+(1/2,1/4)$) node (boundaryLabel) [right, yshift=-4pt, xshift=-2pt] {$\mathbf{x}_\alpha(\mathbf{b})$}; 
        \draw[dotted,->, thick] (xi) -- ( $(xi)+(2/6, 12/30)$ ) node (boundaryLabel) [left, yshift = -5pt, xshift=-4pt] {$\mathbf{x}_\beta(\mathbf{b})$}; 

        \end{tikzpicture}
        \caption{ Schematic overview of a rectangular domain where both an $\alpha$- and $\beta$-characteristic enter the domain at $y=y_\text{min}$.}
        \label{fig:0boundaryChars}
    \end{figure}
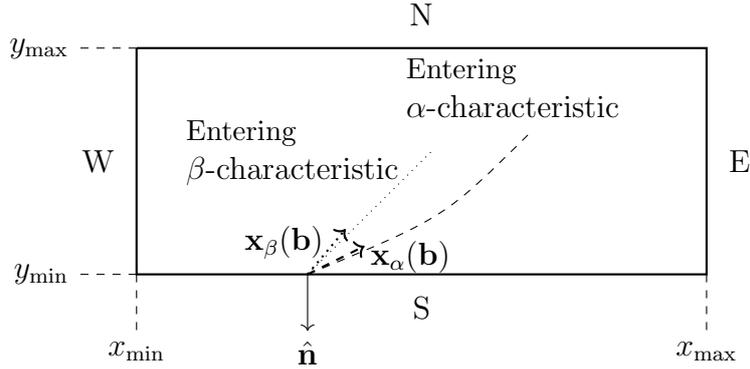

    \textit{Case 2. Two characteristics entering the domain, i.e.,}
    \begin{align}
    \mathbf{x}_\alpha \boldsymbol{\cdot} \hat{\mathbf{n}} < 0, \hspace{5pt} \mathbf{x}_\beta \boldsymbol{\cdot} \hat{\mathbf{n}} < 0.
    \end{align} 
    In this situation the values of $u, p, q, a$ and $b$ cannot be determined (see Figure~\ref{fig:0boundaryChars}). Therefore we prescribe $u$ and its derivative normal to the line segment, which is $q$ if the line segments is horizontal (as in Figure~\ref{fig:0boundaryChars}), or $p$ if it is vertical. 
    The calculation of the relevant variables at $y=y_\text{min}$ and $y = y_\text{max}$ follow analogously to~\eqref{eqn:initialStripCalculations}. As example we consider the line segment at $y = y_\text{min}$. Let $u_\text{S}, q_\text{S}: [x_{\text{min}}, x_{\text{max}}] \rightarrow\mathbb{R}$ be given and let $u(x, y_\text{min}) = u_\text{S}(x)$ and $q(x, y_\text{min}) = q_\text{S}(x)$. We obtain $u, p, a, b, r, s$ and $t$ at $y = y_\text{b}$ via
    \begin{subequations}
    \begin{align}
        p(x, y_\text{min}) & = u^\prime_\text{S}(x), \\
        r(x, y_\text{min}) & = u^{\prime\prime}_\text{S}(x), \\
        s(x, y_\text{min}) & = q^\prime_\text{S}(x), \\
        a(x, y_\text{min}) & = - \frac{r(x, y_\text{min})}{s(x, y_\text{min}) + f_\text{S}(x)}, \\
        b(x, y_\text{min}) & = - \frac{r(x, y_\text{min})}{s(x, y_\text{min}) - f_\text{S}(x)}, \\
        t(x, y_\text{min}) & = \frac{2 f_\text{S}(x)}{a(x, y_\text{min}) - b(x, y_\text{min})}.
    \end{align}
    \end{subequations}
    We require $u^{\prime\prime}_\text{S}(x) \neq 0$ in this case, such that the hyperbolicity condition $a(x, y_\text{b}) \neq b(x, y_\text{b})$ is satisfied.
    Note that the initial strip is one example of \textit{Case 2}, where two characteristics enter. To see this, note that
    \begin{align}
    \begin{split}
    \mathbf{x}_\alpha \cdot \mathbf{\hat{n}} = - x_\alpha = -1 < 0,  \\
    \mathbf{x}_\beta \cdot \mathbf{\hat{n}} = - x_\beta = -1 < 0,
    \end{split}
    \end{align}
    thus classifying both as entering characteristics.

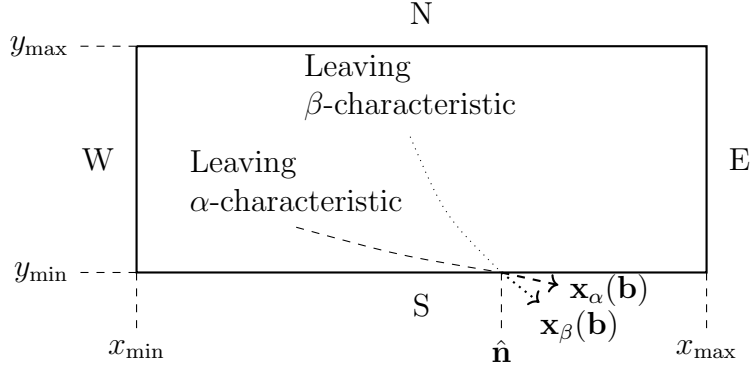
\begin{figure}[H]
        \centering
        \captionsetup{width=0.8\linewidth}
        \begin{tikzpicture}[scale=1.5]
        \coordinate (xi) at (3.2,0);
        \coordinate (alphaStart) at (1.4,0.4);
        \coordinate (betaStart) at (2.4,1.2);

        \draw[dashed] (alphaStart) .. controls (2.3, 0.15) .. (xi); 
        \draw[dotted] (betaStart) .. controls (2.7, 0.5) .. (xi); 
        \node[label={[align=left]Leaving\\$\alpha$-characteristic}, yshift=-2pt] at (alphaStart) {};
        \node[label={[align=left]Leaving\\$\beta$-characteristic}, yshift=-2pt] at (betaStart) {};
        
        \draw [-,thick] (0,2) |- (0,0) |- (5,0) |- (5,2) |- (0,2);
        
        \node[above] at (2.5, 2.1) {N};
        \node[left ] at (-.1, 1) {W};
        \node[below] at (2.5, -0.1) {S};
        \node[right] at (5.1, 1) {E};
        
        \draw [dashed] (0,0) |- (-1/2,0) node (ymin) [left] {$y_\text{min}$};
        \draw [dashed] (0,2) |- (-1/2,2) node (ymax) [left] {$y_\text{max}$};
        \draw [dashed] (0,0) -- (0,-1/2) node (xmin) [below] {$x_\text{min}$};
        \draw [dashed] (5,0) -- (5,-1/2) node (xmax) [below] {$x_\text{max}$};
        \draw [dashed] (xi) -- ($(xi)-(0,1/2)$) node (xilabelNode) [below] {$\hat{\mathbf{n}}$};
        
        \draw[dashed,->, thick] (xi) -- ($(xi)+(1/2, -1/9)$) node (boundaryLabel) [right, yshift=-2pt, xshift=0pt] {$\mathbf{x}_\alpha(\mathbf{b})$}; 
        \draw[dotted,->, thick] (xi) -- ( $(xi)+(2/6, -8/30)$ ) node (boundaryLabel) [right, yshift = -9pt, xshift=-4pt] {$\mathbf{x}_\beta(\mathbf{b})$}; 
        
        \end{tikzpicture}
        \caption{ Schematic overview of a rectangular domain where both an $\alpha$- and $\beta$-characteristic leave the domain $y=y_\text{min}$. }
        \label{fig:2boundaryChars}
    \end{figure}

    \textit{Case 3. Two characteristic leaving the domain, i.e.,}
    \begin{align}
    \mathbf{x}_\alpha \boldsymbol{\cdot} \hat{\mathbf{n}} > 0,  \hspace{5pt} \mathbf{x}_\beta \boldsymbol{\cdot} \hat{\mathbf{n}} > 0.
    \end{align} 
     Here, we should not prescribe anything at all as all values are known, or can be determined by integrating the ODE systems of the $\alpha$- and $\beta$-characteristics, see Figure~\ref{fig:2boundaryChars}. Note that this situation is identical to that of an interior point. 
            
    Note that at $x = x_{\text{max}}$, $x_\alpha = 1$ and $x_\beta = 1$ such that $\mathbf{x}_\alpha \boldsymbol{\cdot} \hat{\mathbf{n}} > 0$ and $\mathbf{x}_\beta \boldsymbol{\cdot} \hat{\mathbf{n}} > 0$, hence this boundary segment coincides with Case 3. The possible exception being the corner points $(x_{\text{max}}, y_\text{min})$ and $(x_{\text{max}}, y_\text{max})$ because the normal is not uniquely defined. In this case the point should be treated as in the Cases 1 or 2. 
    Classifying the boundary segment $x = x_{\text{max}}$ as Case 3 is a direct consequence of choice $\kappa^\alpha = \kappa^\beta = t$, such that $x_\alpha = 1$ and $x_\beta = 1$ for the two characteristic families.
            

    \section{Numerical methods}
    \label{sec:numericalMethods}   
    In order to solve~\eqref{eqn:evolutionFullyReduced} we rewrite it as
    \begin{subequations}
    \label{eqn:odeSystem_vwg}
    \begin{align}
        & \odiff{\mathbf{v}^\alpha}{\alpha} = \mathbf{g}^\alpha(\mathbf{v}^\alpha, a), \hspace{20pt} \odiff{\mathbf{v}^\beta}{\beta} = \mathbf{g}^\beta (\mathbf{v}^\beta, b), \\
    \mathbf{v}^\alpha = \begin{pmatrix} x \\ y \\ u \\ p \\ q \\ b \end{pmatrix}, \hspace{5pt}
    \mathbf{v}^\beta = & \begin{pmatrix} x \\ y \\ u \\ p \\ q \\ a \end{pmatrix}, \hspace{5pt}
    \mathbf{g}^\alpha = \begin{pmatrix} 1 \\ a \\ p + aq \\ - a f \\ f \\ \frac{b-a}{2 f}(f_x + b f_y) \end{pmatrix}, \hspace{5pt}
    \mathbf{g}^\beta = \begin{pmatrix} 1 \\ b \\ p + bq \\ b f \\ -f \\ \frac{a-b}{2 f}(f_x + a f_y) \end{pmatrix}.
    \end{align}
    \end{subequations}
    Equations~\eqref{eqn:odeSystem_vwg} are two mutually coupled systems because the evolution of $a$ and $b$ are determined on the other characteristic.
    By supplying initial conditions, the problem can be treated as a Cauchy problem which we solve by numerical integration.
   
    For our numerical grid we choose $N_x$ points in the $x$-direction and $N_y$ in the $y$-direction. Let the grid points $\mathbf{x}_{i,j}$ be given by $\mathbf{x}_{i,j} = (x_i, y_j)$ for $i = 1, \dots, N_x, j = 1,\dots N_y$. We choose the grid to be equidistant in the $y$-direction with spacing $h_y = (y_{\text{max}} - y_{\text{min}})/(N_y - 1)$. The grid spacing in the $x$-direction does not need to be equidistant, i.e., we write $(h_x)_i = x_{i+1} - x_{i}$. This adaptive stepsize will be detailed in Section~\ref{sec:dynamicStepSize}.
    We denote the numerical approximation of $u$ in a grid point as $u_{i,j} \approx u(\mathbf{x}_{i,j})$, and likewise for the other variables.
    
    When discussing numerical methods we generally consider one step at a time, i.e., we consider the evolution from the grid line $x = {x_i}$ to the line $x = x_{i+1}$. Therefore it is convenient to write $h_x = (h_x)_i$ when no ambiguity arises.


    \subsection{Numerical method based on forward Euler}
    \label{sec:ForwardEulerMethod}
    In this section we will introduce a numerical scheme based on the forward Euler method to calculate $\mathbf{v}_{i+1,j}^\alpha$ and $\mathbf{v}_{i+1,j}^\beta$ given $\mathbf{v}_{i,j}^\alpha$ and $\mathbf{v}_{i,j}^\beta$. 
\begin{figure}[H]
    \centering
    \captionsetup{width=0.85\linewidth}
    \scalebox{0.8}{
        \begin{tikzpicture}
        \coordinate (XAxisMin) at (-2,-2);
        \coordinate (XAxisMax) at (-1,-2);
        \coordinate (YAxisMin) at (-2,-2);
        \coordinate (YAxisMax) at (-2,-1);
        
        \draw [thin,-latex] (XAxisMin) -- (XAxisMax) node [right] {$x$};
        \draw [thin,-latex] (YAxisMin) -- (YAxisMax) node [above] {$y$};
        
        \foreach \x in {0, 2, 3.5}{
            \foreach \y in {0, 1,...,5}{
                \node[draw,circle,inner sep=2pt,fill] at (\x, \y) {};      
            }
        }
        
        \draw[line width=0.2mm] (4,0) -- (0,0) {};
        \draw[line width=0.2mm] (4,5) -- (0,5) {};
        
        \draw[line width=0.2mm] (0,0) -- (0,5) {};
        
        \draw[dashed, line width=0mm] (4,0) -- (-0.5,0) node [left] {$j=1$};
        \draw[dashed, line width=0mm] (4,1) -- (-0.5,1) node [left] {$j=2$};
        
        \draw[dashed, line width=0mm] (0,5) -- (0,-1)node (i1) [below] {$i=1$};
        \draw[dashed, line width=0mm] (2,5) -- (2,-1) node (i2) [below] {$i=2$};
        \draw[dashed, line width=0mm] (3.5,5) -- (3.5,-1) node (i3) [below] {$i=3$};
        
        \draw [decorate,decoration={brace,amplitude=5pt},xshift=-4pt,yshift=0pt]
        (0,3) -- (0,4.0) node [black,midway,xshift=-0.4cm] {\footnotesize $h_y$};
        
        \draw [decorate,decoration={brace,amplitude=5pt},xshift=0pt,yshift=-4pt]
        (2,0) -- (0,0) node [black,midway,yshift=-0.4cm] {\footnotesize $(h_x)_1$};            
        \draw [decorate,decoration={brace,amplitude=5pt},xshift=0pt,yshift=-4pt]
        (3.5, 0) -- (2,0) node [black,midway,yshift=-0.4cm] {\footnotesize $(h_x)_2$};        
        
        \draw[-{Latex[length=3mm]}, blue] (0,5) -- (2, 5.5) {};
        \draw[-{Latex[length=3mm]}, red, dashed] (0,5) -- (2, 4.7) {};
        
        \draw[-{Latex[length=3mm]}, blue] (0,4) -- (2, 4.5) {};
        \draw[-{Latex[length=3mm]}, red, dashed] (0,4) -- (2, 3.7) {};
        
        \draw[-{Latex[length=3mm]}, blue] (0,3) -- (2, 3.5) {};
        \draw[-{Latex[length=3mm]}, red, dashed] (0,3) -- (2, 2.7) {};
        
        \end{tikzpicture}
        \hspace{40pt}
        \begin{tikzpicture}
        \coordinate (Origin)   at (0,0);
        \coordinate (XAxisMin) at (-2,-2);
        \coordinate (XAxisMax) at (-1,-2);
        \coordinate (YAxisMin) at (-2,-2);
        \coordinate (YAxisMax) at (-2,-1);
        \draw [thin,-latex] (XAxisMin) -- (XAxisMax) node [right] {$x$};
        \draw [thin,-latex] (YAxisMin) -- (YAxisMax) node [above] {$y$};
        
        
        \node[draw,circle,inner sep=2pt,fill] at (0,0) {};
        \node[draw,circle,inner sep=2pt,fill] at (0,2) {};                    
        \node[draw,circle,inner sep=2pt,fill] at (0,4) {};
        \node[draw,circle,inner sep=2pt,fill] at (3,0) {};
        \node[draw,circle,inner sep=2pt,fill] at (3,2) {};                    
        \node[draw,circle,inner sep=2pt,fill] at (3,4) {};
        
        \node [below] at (0,2) {$\mathbf{x}_{i,j}$};
        \node [below] at (3,2) {$\mathbf{x}_{i+1,j}$};
        
        \draw[-{Latex[length=3mm]}, blue] (0,2) -- (3, 3.2) node [right, text=black] {$\widetilde{\mathbf{v}}^\alpha_{i+1}(j)$};
        \draw[-{Latex[length=3mm]}, red, dashed] (0,2) -- (3, 0.4) node [right, text=black] {$\widetilde{\mathbf{v}}^\beta_{i+1}(j)$};
        \node[above] at (1.5, 2.6) {$\alpha$};
        \node[below] at (1.5, 1.2) {$\beta$};
        
        \draw[dashed, line width=0mm] (0, 4.5) -- (0,-0.7)node [below] {$i$};
        \draw[dashed, line width=0mm] (3, 4.5) -- (3,-0.7) node [below] {$i+1$};
        
        \draw[dashed, line width=0mm] (3.5, 4) -- (-0.5,4) node [left] {$j+1$};
        \draw[dashed, line width=0mm] (3.5, 2) -- (-0.5,2) node [left] {$j$};
        \draw[dashed, line width=0mm] (3.5, 0) -- (-0.5,0) node [left] {$j-1$};
        
        \draw [decorate,decoration={brace,amplitude=10pt},xshift=0pt,yshift=-4pt]
        (3, 0) -- (0,0) node [black,midway,yshift=-0.5cm] {\footnotesize $h_x$}; 
        
        \end{tikzpicture}
        }
        \caption{Schematic representation of the numerical method using forward Euler.}	
        \label{fig:ExplicitEuler}
    \end{figure}
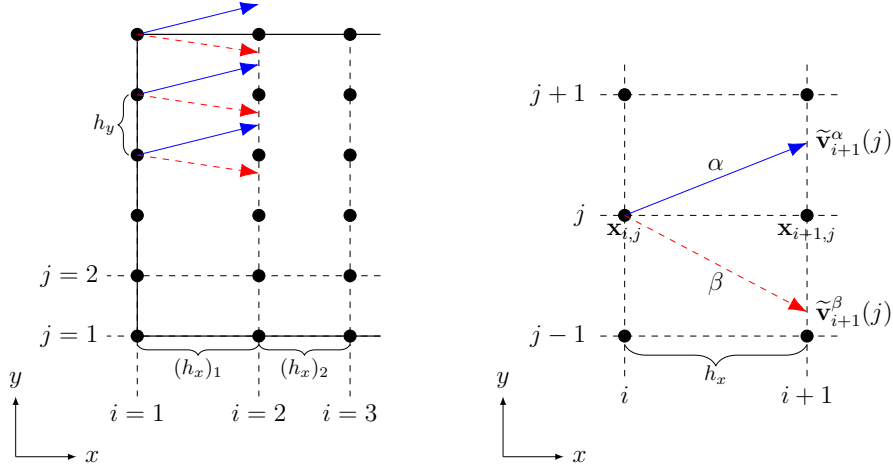

\begin{figure}[t!]
    \centering
    \captionsetup{width=0.85\linewidth}
    \scalebox{0.7}{
    \begin{tikzpicture}
    \coordinate (Origin)   at (0,0);
    \coordinate (XAxisMin) at (-2,-6.6);
    \coordinate (XAxisMax) at (-1,-6.6);
    \coordinate (YAxisMin) at (-2,-6.6);
    \coordinate (YAxisMax) at (-2,-5.6);
    \draw [thin,-latex] (XAxisMin) -- (XAxisMax) node [right] {$x$};
    \draw [thin,-latex] (YAxisMin) -- (YAxisMax) node [above] {$y$};
        
    \node[draw,circle,inner sep=2pt,fill] at (0,0) {};             
    \node[draw,circle,inner sep=2pt,fill] at (0,4) {};
    \node[draw,circle,inner sep=2pt,fill] at (0,-4) {};                    
    \node[draw,circle,inner sep=2pt,fill] at (3,0) {};                    
    \node[draw,circle,inner sep=2pt,fill] at (3,4) {};
    \node[draw,circle,inner sep=2pt,fill] at (3,-4) {};
    
    \node[above] at (1.5, 4.4) {$\alpha$};    
    \node[below] at (1.5, 3.2) {$\beta$};
    \node[above] at (1.5, 0.6) {$\alpha$};
    \node[below] at (1.5, -0.6) {$\beta$};
    \node[above] at (1.5, -3.2) {$\alpha$};
    \node[below] at (1.5, -4.4) {$\beta$};
        
    \draw[dashed, line width=0mm] (0, 5.5) -- (0,-5.5) node [below] {$i$};
    \draw[dashed, line width=0mm] (3, 5.5) -- (3,-5.5) node [below] {$i+1$};
    
    \draw[dashed, line width=0mm] (4.5, 4) -- (-0.5,4) node [left] {$j+1$};
    \draw[dashed, line width=0mm] (4.5, 0) -- (-0.5,0) node [left] {$j$};
    \draw[dashed, line width=0mm] (4.5, -4) -- (-0.5,-4) node [left] {$j-1$};
    
    \draw[->, blue] (0,4) -- (3, 4.8) node [draw,thick,minimum width=0.2cm,minimum height=0.2cm, fill] {};
    \draw[->, red] (0,4) -- (3, 2.4) node [draw,thick,minimum width=0.2cm,minimum height=0.2cm, fill] {};
    
    \draw[->, blue] (0,0) -- (3, 1.2) node [draw,thick,minimum width=0.2cm,minimum height=0.2cm, fill] {};
    \draw[->, red] (0,0) -- (3, -1.2) node [draw,thick,minimum width=0.2cm,minimum height=0.2cm, fill] {};
    
    \draw[->, blue] (0,-4) -- (3, -2.4) node [draw,thick,minimum width=0.2cm,minimum height=0.2cm, fill] {};
    \draw[->, red] (0,-4) -- (3, -4.8) node [draw,thick,minimum width=0.2cm,minimum height=0.2cm, fill] {};
    
    \draw [decorate,decoration={brace,amplitude=12pt,aspect=0.222222, raise=0pt},xshift=0pt,yshift=0pt, blue]
    (4, 4.8) -- (4, 1.2) {};
    
    \draw [decorate,decoration={brace,amplitude=12pt,aspect=0.333333, raise=0pt},xshift=0pt,yshift=0pt,blue]
    (4, 1.2) -- (4, -2.4) {};
    
    \draw[blue] (3, 4.8) -- (4, 4.8) {};
    \draw[blue] (3, 1.2) -- (4, 1.2) {};
    \draw[blue] (3, -2.4) -- (4, -2.4) {};
    
   
   \draw [decorate,decoration={brace,amplitude=13pt,aspect=0.6666666, raise=-1pt},xshift=5pt,yshift=0pt,red]
   (3, 2.4) -- (3, -1.2) {};
  \draw [decorate,decoration={brace,amplitude=13pt,aspect=0.7777777, raise=-1pt},xshift=5pt,yshift=0pt,red]
   (3, -1.2) -- (3, -4.8) {};
   
    \end{tikzpicture}
    }
    \hspace{40pt}
    \scalebox{0.7}{
    \begin{tikzpicture}
    \coordinate (Origin)   at (0,0);
    \coordinate (XAxisMin) at (-2,-6.6);
    \coordinate (XAxisMax) at (-1,-6.6);
    \coordinate (YAxisMin) at (-2,-6.6);
    \coordinate (YAxisMax) at (-2,-5.6);
    \draw [thin,-latex] (XAxisMin) -- (XAxisMax) node [right] {$x$};
    \draw [thin,-latex] (YAxisMin) -- (YAxisMax) node [above] {$y$};
    
    \node[draw,circle,inner sep=2pt,fill] at (0,0) {};             
    \node[draw,circle,inner sep=2pt,fill] at (0,4) {};
    \node[draw,circle,inner sep=2pt,fill] at (0,-4) {};                    
    \node[draw,circle,inner sep=2pt,fill] at (3,0) {};                    
    \node[draw,circle,inner sep=2pt,fill] at (3,4) {};
    \node[draw,circle,inner sep=2pt,fill] at (3,-4) {};
    
    \node[above] at (1.5, 4.4) {$\alpha$};    
    \node[below] at (1.5, 3.2) {$\beta$};
    \node[above] at (1.5, 0.6) {$\alpha$};
    \node[below] at (1.5, -0.6) {$\beta$};
    \node[above] at (1.5, -3.2) {$\alpha$};
    \node[below] at (1.5, -4.4) {$\beta$};
    
    \draw[dashed, line width=0mm] (0, 5.5) -- (0,-5.5) node [below] {$i$};
    \draw[dashed, line width=0mm] (3, 5.5) -- (3,-5.5) node [below] {$i+1$};
    
    \draw[dashed, line width=0mm] (4, 4) -- (-0.5,4) node [left] {$j+1$};
    \draw[dashed, line width=0mm] (4, 0) -- (-0.5,0) node [left] {$j$};
    \draw[dashed, line width=0mm] (4, -4) -- (-0.5,-4) node [left] {$j-1$};
    
    \draw[->, blue] (0,4) -- (3, 4.8) node [draw,thick,minimum width=0.2cm,minimum height=0.2cm, fill] {};
    \draw[->, red] (0,4) -- (3, 2.4) node [draw,thick,minimum width=0.2cm,minimum height=0.2cm, fill] {};
    
    \draw[->, blue] (0,0) -- (3, 1.2) node [draw,thick,minimum width=0.2cm,minimum height=0.2cm, fill] {};
    \draw[->, red] (0,0) -- (3, -1.2) node [draw,thick,minimum width=0.2cm,minimum height=0.2cm, fill] {};
    
    \draw[->, blue] (0,-4) -- (3, -2.4) node [draw,thick,minimum width=0.2cm,minimum height=0.2cm, fill] {};
    \draw[->, red] (0,-4) -- (3, -4.8) node [draw,thick,minimum width=0.2cm,minimum height=0.2cm, fill] {};
    
    \draw [decorate, decoration={brace,amplitude=12pt,aspect=0.333333},xshift=4pt,yshift=0pt,black]
    (3, 4.8) -- (3, 2.4) {};
    
    \draw [decorate,decoration={brace,amplitude=12pt,aspect=0.5},xshift=4pt,yshift=0pt,black]
    (3, 1.2) -- (3, -1.2) {};
    
    
    \draw [decorate,decoration={brace,amplitude=12pt,aspect=0.6666666},xshift=4pt,yshift=0pt,black]
    (3, -2.4) -- (3, -4.8) {};

 
    \end{tikzpicture}
    }
    \caption{Schematic representation of linear interpolation where a curly bracket denotes the values used for interpolation and the gridpoint it influences. Interpolation for $u$, $p$ and $q$ can be done either using two distinct sets, formed by the two characteristic families (Approach 1, left), or using both sets combined (Approach 2, right).}
    \label{fig:ExplicitEulerStep}
    \end{figure}
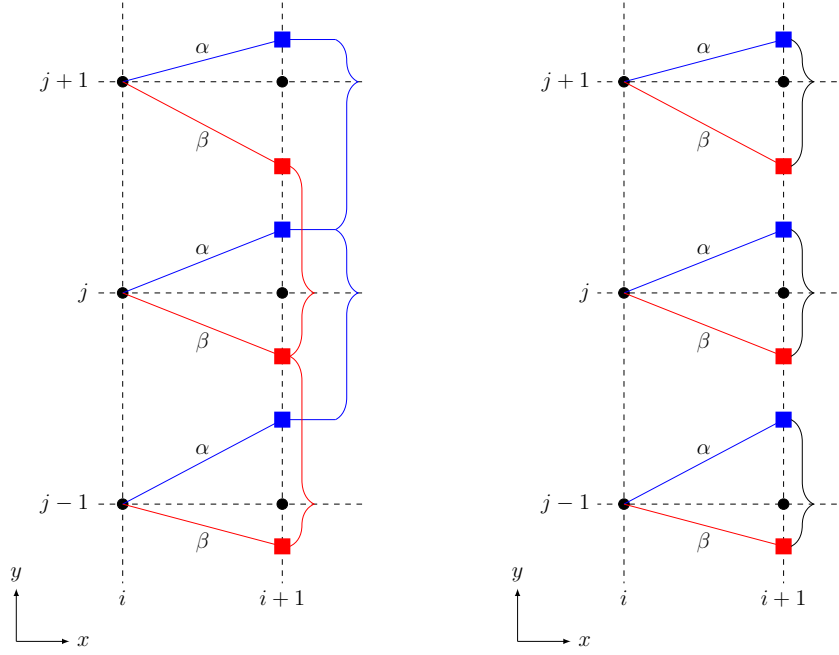

    The numerical stencil is schematically shown in Figure \ref{fig:ExplicitEuler}. The black dots represent the grid points. The solid blue and dashed red arrows correspond to the numerical approximation of the $\alpha$- and $\beta$ characteristic, respectively. At the grid points these characteristics are approximated using forward Euler, i.e., we approximate the characteristics as tangent lines passing through grid points $\mathbf{x}_{i,j}$ and having slope $a_{i,j}$ for the $\alpha$-characteristic and $b_{i,j}$ for the $\beta$-characteristic. 
    This implies that for a step size $h_x$, the two characteristics departing from $\mathbf{x}_{i,j}$ arrive at $(x_{i+1} ,\tilde{y}^\alpha_{i+1}(j))$ and $(x_{i+1}, \tilde{y}^\beta_{i+1}(j))$ for the $\alpha$- and $\beta$-characteristic respectively, where
    \begin{subequations}
        \label{eqn:characteristicPointsEE}
        \begin{align}
        x_{i+1} & = x_{i} + h_x, \\
        \tilde{y}^\alpha_{i+1}(j) & = y_{j} + h_x a_{i,j}, \\
        \tilde{y}^\beta_{i+1}(j) & = y_{j} + h_x b_{i,j}.
        \end{align}
    \end{subequations}
    The points $(x_{i+1} ,\tilde{y}^\alpha_{i+1}(j))$ and $(x_{i+1}, \tilde{y}^\beta_{i+1}(j))$ do generally not coincide with any gridpoint.
    Similarly we can calculate $u, p, q$ on both characteristics, while $a$ and $b$ can only be determined on one characteristic. More compactly written, we have the following forward Euler step:
    \begin{align}
    \widetilde{\mathbf{v}}^\alpha_{i+1}(j) = \mathbf{v}^\alpha_{i,j} + h_x \mathbf{g}^\alpha(\mathbf{v}^\alpha_{i,j}, a_{i,j}), \hspace{10pt} \widetilde{\mathbf{v}}^\beta_{i+1}(j) = \mathbf{v}^\beta_{i,j} + h_x \mathbf{g}^\beta(\mathbf{v}^\beta_{i,j}, b_{i,j}).
    \end{align}
    Here $\widetilde{\mathbf{v}}^\alpha_{i+1}(j)$ denotes the new values of $\mathbf{v}^\alpha$ at $(x_{i+1}, \tilde{y}^\alpha_{i+1}(j))$ for which the corresponding characteristic passes through the grid point $\mathbf{x}_{i,j}$ as shown in Figure~\ref{fig:ExplicitEuler}. Analogously we define $\widetilde{\mathbf{v}}^\beta_{i+1}(j)$.
    Because we are interested in obtaining $\mathbf{v}^\alpha(x_{i+1}, y_{j})$ and $\mathbf{v}^\beta(x_{i+1}, y_{j})$, we interpolate $\widetilde{\mathbf{v}}^\alpha_{i+1}(j)$ and $\widetilde{\mathbf{v}}^\beta_{i+1}(j)$. 
    Let the $y, u, p, q$-components of $\tilde{\mathbf{v}}^\alpha_{i+1}(j)$ be denoted by $\tilde{y}^\alpha_{i+1}(j), \tilde{u}^\alpha_{i+1}(j), \tilde{p}^\alpha_{i+1}(j)$ and $\tilde{q}^\alpha_{i+1}(j)$, respectively, and similarly for $\tilde{\mathbf{v}}^\beta_{i+1}(j)$. Two approaches to carry out the interpolation are shown in Figure~\ref{fig:ExplicitEulerStep}, using linear interpolation.
    
    \label{sec:interpolationOnFullGridLine}
    \textit{Approach 1:}
    We consider all values computed for each characteristic family, e.g. $\tilde{u}_{i+1}^\alpha(j)$, as one set and interpolate using local B-Splines (\hspace{-3pt}\cite[p. 90-97]{splinesBook}; see App.~\ref{sec:interpolationOnFullGridLine}) within that set, which yields the approximations $u_{i+1,j}^\alpha$ and $u_{i+1,j}^\beta$.
    Because we have no a priori preference whether $u_{i+1,j}^\alpha$ or $u_{i+1,j}^\beta$ approximates  $u_{i+1,j}$ better, we average the results and set
    \begin{align}
    u_{i+1,j} = \frac{u_{i+1,j}^\alpha + u_{i+1,j}^\beta}{2}.
    \end{align}
    Similarly we obtain $p_{i+1,j}$ and $q_{i+1,j}$.
    
    \textit{Approach 2:}
    We collect $\tilde{u}_{i+1}^\alpha(j)$ and $\tilde{u}_{i+1}^\beta(j)$ for all $j = 1, \dots, N_y$ in one set and interpolate using local B-Splines over that set to obtain $u_{i+1,j}$. Similarly we obtain $p_{i+1,j}$ and $q_{i+1,j}$.
    
    
    Numerical results show that neither of these methods outperforms the other significantly or consistently. Computationally both approaches are approximately equally expensive \cite{interpolationComplexity}.
    Approach 1 will be used throughout. 
    
    Various subtleties arise using the approach above. For one, $a$ and $b$ are known along one family of characteristics only, so we use one set of values for interpolation of these variables.
    In the case a grid point is not located between two numerically estimated characteristics, which may occur for a grid point on the boundary, then we supplement the missing boundary value as detailed in Section~\ref{sec:BoundaryConditionsHMA}. 
    
    
    To determine an appropriate spline interpolant we consider both the error associated with one forward Euler step and that of the interpolation. The local truncation error of the forward Euler method is $\mathcal{O}(h_x^2)$  \cite[p. 335]{NumericalAnalysisBook}. In Section~\ref{sec:dynamicStepSize} we introduce a method to control the step size $h_x$, which consequently implies that the local truncation error behaves as $\mathcal{O}(h_y^2)$. 
    The error of a spline interpolant of order $n$, also called a spline interpolant of degree $n-1$, is $\mathcal{O}(h_y^{n})$ \cite[p. 95]{splinesBook}. In general we choose the degree of the spline interpolant such that its order matches the order of the integration method. For the forward Euler method this implies a spline interpolant of order 2 which is standard local linear interpolation. 

    \subsection{Numerical method based on modified Euler}
    \label{sec:ModifiedEulerMethod}

    In this section we will discuss a numerical scheme based on the modified Euler scheme. The local truncation error of the modified Euler method is $\mathcal{O}(h_x^3)$.

    First we calculate a predictor by doing a forward Euler step of step size $\frac{h_x}{2}$ for~\eqref{eqn:odeSystem_vwg}, viz.
    \begin{align}
        \widetilde{\mathbf{v}}^\alpha_{i+\frac{1}{2}}(j) = \mathbf{v}^\alpha_{i,j} + \frac{h_x}{2} \mathbf{g}^\alpha(\mathbf{v}^\alpha_{i,j}, a_{i,j}), \hspace{30pt} \widetilde{\mathbf{v}}^\beta_{i+\frac{1}{2}}(j) = \mathbf{v}^\beta_{i,j} + \frac{h_x}{2}\mathbf{g}^\beta(\mathbf{v}^\beta_{i,j}, b_{i,j}).
    \end{align}
    We adhere to the same notation as in the previous section, where a tilde denotes the function value on the corresponding characteristic which is not necessarily located at a grid point.
    
    Because $a$ and $b$ are not known but are needed at $x_{i+1/2}$ along the $\alpha$- and $\beta$-characteristic, respectively, we approximate them using interpolation. 
    We interpolate the $a$-component of $\widetilde{\mathbf{v}}^\beta_{i+\frac{1}{2}}(j)$, known at the points $\widetilde{y}^\beta_{i+\frac{1}{2}}(j)$ to the points $\widetilde{y}^\alpha_{i+\frac{1}{2}}(j)$, which approximate $a(x_{i+\frac{1}{2}}, \tilde{y}^\alpha_{i+\frac{1}{2}}(j))$. We denote this approximation by $a_{i+\frac{1}{2},j}$. In the same way we approximate $b_{i+\frac{1}{2},j}$ from $\widetilde{\mathbf{v}}^\alpha_{i+\frac{1}{2}}(j)$.
    The modified Euler step is then given by
    \begin{align}
    \label{eqn:modifiedEulerFullStep}
    \widetilde{\mathbf{v}}^\alpha_{i+1}(j) = \mathbf{v}^\alpha_{i,j} + h_x \mathbf{g}^\alpha(\widetilde{\mathbf{v}}^\alpha_{i+\frac{1}{2}}(j), a_{i+\frac{1}{2},j}), \hspace{10pt}
    \widetilde{\mathbf{v}}^\beta_{i+1}(j) = \mathbf{v}^\beta_{i,j} + h_x \mathbf{g}^\beta(\widetilde{\mathbf{v}}^\beta_{i+\frac{1}{2}}(j), b_{i+\frac{1}{2},j}).
    \end{align}

    We conclude the modified Euler step by interpolating to the grid points using Approach 1, as discussed in Section~\ref{sec:interpolationOnFullGridLine}. We choose third order spline interpolation for the modified Euler method, as it corresponds to its local truncation error. 

    \subsection{Numerical method based on classic Runge-Kutta}
    \label{sec:RungeKuttaMethod}
    In this section we will introduce the classic Runge-Kutta method. Following a similar approach, we can generalize our integration methods to other higher order explicit Runge-Kutta methods. 
    
    First we make a forward Euler step of size $\frac{h_x}{2}$ for~\eqref{eqn:odeSystem_vwg} viz.
    \begin{align}
    \widetilde{\mathbf{v}}^\alpha_{i+\frac{1}{2}}(j) = \mathbf{v}^\alpha_{i,j} + \frac{h_x}{2} \mathbf{g}^\alpha(\mathbf{v}^\alpha_{i,j}, a_{i,j}), 
    \hspace{10pt}
    \widetilde{\mathbf{v}}^\beta_{i+\frac{1}{2}}(j) = \mathbf{v}^\beta_{i,j} + \frac{h_x}{2}\mathbf{g}^\beta(\mathbf{v}^\beta_{i,j}, b_{i,j}).
    \end{align}
    As in the case of the modified Euler based method, we interpolate the $a$- and $b$-components of $\widetilde{\mathbf{v}}^\beta_{i+\frac{1}{2}}(j)$ and $\widetilde{\mathbf{v}}^\alpha_{i+\frac{1}{2}}(j)$ to approximate $a_{i+\frac{1}{2},j}$ and $b_{i+\frac{1}{2},j}$, respectively.
    Second, we do a step of size $\frac{h_x}{2}$ with the slope based on the previously found values of $\widetilde{\mathbf{v}}^\alpha$ and $\widetilde{\mathbf{v}}^\beta$, viz.
    \begin{align}
    \widehat{\mathbf{v}}^\alpha_{i+\frac{1}{2}}(j) = \mathbf{v}^\alpha_{i,j} + \frac{h_x}{2} \mathbf{g}^\alpha(\widetilde{\mathbf{v}}^\alpha_{i+\frac{1}{2}}(j), a_{i+\frac{1}{2},j}), \hspace{10pt}
    \widehat{\mathbf{v}}^\beta_{i+\frac{1}{2}}(j) = \mathbf{v}^\beta_{i,j} + \frac{h_x}{2} \mathbf{g}^\beta(\widetilde{\mathbf{v}}^\beta_{i+\frac{1}{2}}(j), b_{i+\frac{1}{2},j}),
    \end{align}
    where we used a hat to distinguish between the different stages.
    Similar as before we interpolate $a$ and $b$ from the $\beta$- and $\alpha$-characteristics to the $\alpha$- and $\beta$-characteristics to obtain 
    $\hat{a}_{i+\frac{1}{2},j}$ and $\hat{b}_{i+\frac{1}{2},j}$, respectively. Using these slopes we do a step of size $h_x$ which yields
    \begin{align}
    \widehat{\mathbf{v}}^\alpha_{i+1}(j) = \mathbf{v}^\alpha_{i,j} + h_x \mathbf{g}^\alpha(\widehat{\mathbf{v}}^\alpha_{i+\frac{1}{2}}(j), \hat{a}_{i+\frac{1}{2},j}), \hspace{10pt}
    \widehat{\mathbf{v}}^\beta_{i+1}(j) = \mathbf{v}^\beta_{i,j} + h_x \mathbf{g}^\beta(\widehat{\mathbf{v}}^\beta_{i+\frac{1}{2}}(j), \hat{b}_{i+\frac{1}{2},j}).
    \end{align}
    Interpolating the $a$- and $b$-components of $\widehat{\mathbf{v}}^\beta_{i+1}(j)$ and $\widehat{\mathbf{v}}^\alpha_{i+1}(j)$ yields approximations for  $\hat{a}^\alpha_{i+1}(j)$ and $\hat{b}^\beta_{i+1}(j)$, respectively.
    Finally the full Runge-Kutta step is given by
    \begin{subequations}
    \begin{align}
    \begin{split}
    \widetilde{\mathbf{v}}^\alpha_{i+1}(j) = \mathbf{v}^\alpha_{i,j} + \frac{h_x}{6} \Big(&
    \mathbf{g}^\alpha(\mathbf{v}^\alpha_{i,j}, a_{i,j}) + 
    2 \mathbf{g}^\alpha(\widetilde{\mathbf{v}}^\alpha_{i+\frac{1}{2}}(j), a_{i+\frac{1}{2},j}) + 
    \\
    &
    2 \mathbf{g}^\alpha(\widehat{\mathbf{v}}^\alpha_{i+\frac{1}{2}}(j), \hat{a}_{i+\frac{1}{2},j}) + 
    \mathbf{g}^\alpha(\widehat{\mathbf{v}}^\alpha_{i+1}(j), \hat{a}^\alpha_{i+1}(j)) 
    \Big),
    \end{split}
    \\
    \begin{split}
    \widetilde{\mathbf{v}}^\beta_{i+1}(j) = \mathbf{v}^\beta_{i,j} + \frac{h_x}{6} \Big(&
    \mathbf{g}^\beta(\mathbf{v}^\beta_{i,j}, b_{i,j}) + 
    2 \mathbf{g}^\beta(\widetilde{\mathbf{v}}^\beta_{i+\frac{1}{2}}(j), b_{i+\frac{1}{2},j}) + 
    \\
    &
    2 \mathbf{g}^\beta(\widehat{\mathbf{v}}^\beta_{i+\frac{1}{2}}(j), \hat{b}_{i+\frac{1}{2},j}) + 
    \mathbf{g}^\beta(\widehat{\mathbf{v}}^\beta_{i+1}(j), \hat{b}^\beta_{i+1}(j)) 
    \Big).
    \end{split}
    \end{align}
    \end{subequations}
    
    We conclude the Runge-Kutta step by interpolating $\widetilde{\mathbf{v}}^\alpha_{i+1}(j)$ and $\widetilde{\mathbf{v}}^\beta_{i+1}(j)$ using fifth order spline interpolation to the grid points using Approach 1, as discussed in Section~\ref{sec:interpolationOnFullGridLine}.
        
%

    \subsection{Dynamic step size control}
        \label{sec:dynamicStepSize}
        In this section we introduce a procedure to choose the step size $h_x$ adaptively. We aim to reduce the computational error while maintaining convergence, which depends both on the integration method and the interpolation methods. 
        
\begin{figure}[H]
            \centering
            \captionsetup{width=0.85\linewidth}
            \scalebox{0.8}{
                \begin{tikzpicture}
                \coordinate (XAxisMin) at (-2,-2);
                \coordinate (XAxisMax) at (-1,-2);
                \coordinate (YAxisMin) at (-2,-2);
                \coordinate (YAxisMax) at (-2,-1);
                
                \draw [thin,-latex] (XAxisMin) -- (XAxisMax) node [right] {$x$};
                \draw [thin,-latex] (YAxisMin) -- (YAxisMax) node [above] {$y$};
                
                \foreach \x in {0, 4}{
                    \foreach \y in {0, 1,...,5}{
                        \node[draw,circle,inner sep=2pt,fill] at (\x, \y) {};      
                    }
                }
                
                \draw[line width=0.2mm] (4,0) -- (0,0) {};                
                \draw[line width=0.2mm] (0,0) -- (0,5) {};
                
                \draw[dashed, line width=0mm] (4,0) -- (-0.5,0) node [left] {$j=1$};
                \draw[dashed, line width=0mm] (4,1) -- (-0.5,1) node [left] {$j=2$};
                
                \draw[dashed, line width=0mm] (0,5) -- (0,-0.5) node (i1) [below] {$i$};
                \draw[dashed, line width=0mm] (4,5) -- (4,-0.5) node (i2) [below] {$i+1$};
                
                \draw[dashed, blue] (0,0) .. controls (1.2, 0.5) and (3.2, 2.1) .. (4,3.3) node [midway, below,text=black, xshift=5pt,yshift=5pt] {$\beta$};
                \draw[->, blue] (0,0) -- (4,3.6) node [midway,above,align=left,text=black,xshift=-7pt,yshift=12pt] {$\beta$ estimate};

                \draw [decorate,decoration={brace,amplitude=16pt},xshift=4pt,yshift=0pt]
                (4, 3.6) -- (4,0) node [midway,right,xshift=12pt] {$\Delta y(1)$};
                
                \end{tikzpicture}
                \hspace{40pt}
                \begin{tikzpicture}
                \coordinate (XAxisMin) at (-2,-2);
                \coordinate (XAxisMax) at (-1,-2);
                \coordinate (YAxisMin) at (-2,-2);
                \coordinate (YAxisMax) at (-2,-1);
                
                \draw [thin,-latex] (XAxisMin) -- (XAxisMax) node [right] {$x$};
                \draw [thin,-latex] (YAxisMin) -- (YAxisMax) node [above] {$y$};
                
                \foreach \x in {0, 1}{
                    \foreach \y in {0, 1,...,5}{
                        \node[draw,circle,inner sep=2pt,fill] at (\x, \y) {};      
                    }
                }
                
                \draw[line width=0.2mm] (1,0) -- (0,0) {};                
                \draw[line width=0.2mm] (0,0) -- (0,5) {};
                
                \draw[dashed, line width=0mm] (1,0) -- (-0.5,0) node [left] {$j=1$};
                \draw[dashed, line width=0mm] (1,1) -- (-0.5,1) node [left] {$j=2$};
                
                \draw[dashed, line width=0mm] (0,5) -- (0,-0.5) node (i1) [below] {$i$};
                \draw[dashed, line width=0mm] (1,5) -- (1,-0.5) node (i2) [below] {$i+1$};
                
                \draw[dashed, blue] (0,0) .. controls (1.2, 0.5) and (3.2, 2.1) .. (4,3.3) node [midway, below,text=black, xshift=5pt,yshift=5pt] {$\beta$};
                \draw[->, blue] (0,0) -- (4,3.6) node [midway,above,align=left,text=black,xshift=5pt,yshift=23pt] {$\beta$ estimate};
                
                \draw [decorate,decoration={brace,amplitude=6pt},xshift=4pt,yshift=0pt]
                (1, 0.9) -- (1,0) node [midway,right,xshift=12pt] {$\Delta y(1)$};
                
                \end{tikzpicture}                
            }
            \caption{Two schematic situations of $\Delta y$, as function of $(h_x)_i$.}	
            \label{fig:dynamicStepSize}
        \end{figure}
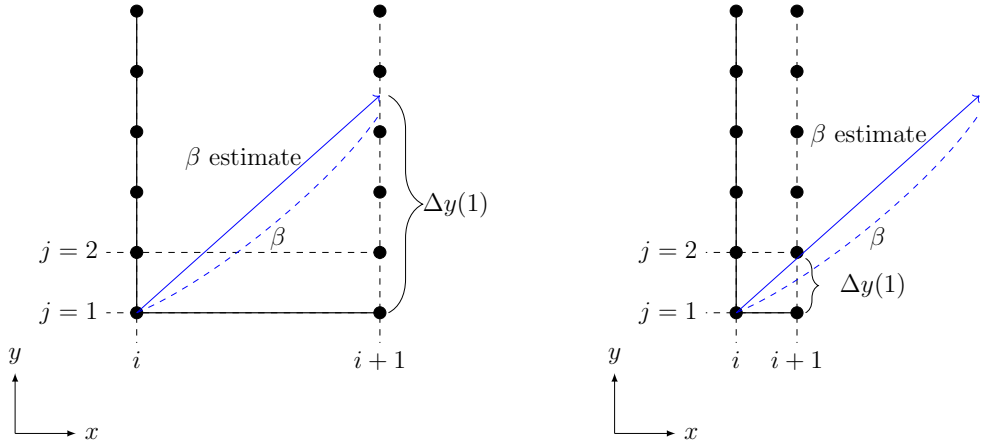

        Because integration is done in the positive $x$-direction, the corresponding discretization error is a function of $h_x$. On the other hand, we interpolate in the $y$-direction, but the interpolation error is not solely a function of $h_y$ as we will see. Ideally we want both the integration and interpolation error to be of the same order, such that neither of them dominates. Asymptotically, this is obtained most easily by using both an integration and interpolation method of the same order, and choosing the discretization steps in the $x$- and $y$-direction of the same order of magnitude as well. In the $x$-direction the discretization step size is $h_x$, though the distance between the interpolation nodes, those points that are the intersection between a grid line and the numerical approximation of the characteristics, is not equidistant, but rather follows from both the evolution along the characteristics~\eqref{eqn:evolutionFullyReduced} and the integration method used.
        Without loss of generality, we solely consider the $\beta$-characteristic. Let $\Delta y(j) = |\tilde{y}^\beta_{i+1}(j) - y_j|$ denote the distance between, on the one hand, the intersection point of (the approximation of) the $\beta$-characteristic at $x = x_{i+1}$, and on the other hand, the point $(x_{i+1}, y_j)$, as shown in Figure~\ref{fig:dynamicStepSize} for $j=1$.
        Approximating $\Delta y(j)$ by a forward Euler step yields
        \begin{align}
        \Delta y(j) = |y_j + (h_x)_i \, b_{i,j} - y_j| = |b_{i,j}| (h_x)_i.
        \end{align}
        Hence $\Delta y(j) \leq h_y$ is obtained if we choose 
        \begin{align}
        (h_x)_i < h_y \min_{j \in \{1, \dots, N_y\}}(1, 1/|b_{i,j}|),
        \end{align}
        where the constant ``$1$'' is chosen such that $(h_x)_i < h_y$, regardless of the slope of the characteristics. This allows us to control the error of the numerical methods by solely controlling $h_y$.
        Similarly we would like to have $(h_x)_i < h_y \min(1, 1/|a_{i,j}|)$ for all $j = 1,\dots, N_y$. Therefore we choose
        \begin{align}
        \label{eqn:h_xAdaptive}
        (h_x)_i = \gamma \, h_y \cdot \min_{j \in \{1, \dots, N_y\}} \left\{1, \left|\frac{1}{a_{i,j}}\right|, \left|\frac{1}{b_{i,j}}\right| \right\},
        \end{align}
        where $0 \leq \gamma \leq 1$ is a tuning parameter. Generally we choose $\gamma = 0.95$, as this implies strict inequality.
        
    \subsection{Residual of the \MAe{}}   
    \label{sec:numericalResidual} 
    We solve the mutually coupled ODE systems~\eqref{eqn:evolutionFullyReduced} without calculating $r, s$ or $t$ nor by calculating any numerical derivatives. Ideally we should be able to determine the numerical residual without numerically taking derivatives either. We do so by formulating the \MAe{} as an integral equation, which we evaluate numerically.

\begin{figure}[H]
        \centering
        \captionsetup{width=0.8\linewidth}
        \begin{tikzpicture}[scale=0.8]
        \node (Aij)  at (2,2) {$A_{i,j}$};
        \draw [-] (0,0) |- (0,4) |- (4,4) |- (4, 0) |- (0,0);
        \draw [dashed] (0,0) |- (-1/2,0) node (ymin) [left] {$y_{j-\frac{1}{2}}$};
        \draw [dashed] (0,4) |- (-1/2,4) node (ymax) [left] {$y_{j+\frac{1}{2}}$};
        \draw [dashed] (0,0) -- (0,-1/2) node (xmin) [below] {$x_{i-\frac{1}{2}}$};
        \draw [dashed] (4,0) -- (4,-1/2) node (xmax) [below] {$x_{i+\frac{1}{2}}$};
        
        \draw[->, blue, ultra thick] (0,2.5) -- (0,1.5) node [midway, right, text=black] {$\hat{\pmb{\tau}}$} node [midway, left, text=black, xshift=-10pt] {W};
        \draw[->, blue, ultra thick] (1.5, 0) -- (2.5,0) node [midway, above, text=black] {$\hat{\pmb{\tau}}$} node [midway, below, text=black, yshift=-10pt] {S};
        \draw[->, blue, ultra thick] (4,1.5) -- (4,2.5) node [midway, left, text=black] {$\hat{\pmb{\tau}}$} node [midway, right, text=black, xshift=10pt] {E};
        \draw[->, blue, ultra thick] (2.5,4) -- (1.5,4) node [midway, below, text=black] {$\hat{\pmb{\tau}}$} node [midway, above, text=black, yshift=10pt] {N};
        \end{tikzpicture}
        \caption{ Schematic overview of domain $A_{i,j}$ with corresponding tangent vectors.}
        \label{fig:windDirectionsIntegrals}
    \end{figure}
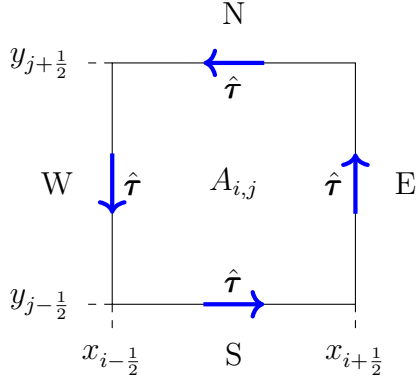

    The \MAe{} can be written as $p_x q_y - p_y q_x + f^2 = 0$, which only depends on $f^2$ and the derivatives of $p$ and $q$. We can rewrite it in terms of $\nabla p$ and $\nabla q$. To this end let $\mathbf{J}$ be the symplectic matrix $\begin{pmatrix} 0 & 1 \\ -1 & 0\end{pmatrix}$, from wich the two equivalent formulations
    \begin{subequations}
    \label{eqn:MA_pq_reformulated}
        \begin{align}
        \label{eqn:MA_pq_reformulated_a}
        -f^2 & = \nabla \boldsymbol{\cdot} (p \mathbf{J} \nabla q), \\
        \label{eqn:MA_pq_reformulated_b}
         f^2 & = \nabla \boldsymbol{\cdot} (q \mathbf{J} \nabla p),
        \end{align}
    \end{subequations}    
    follow, where we used $\nabla \boldsymbol{\cdot} \big(\mathbf{J} \nabla \phi \big) = 0$ for a scalar function $\phi$, and $\nabla \boldsymbol{\cdot} (\phi \mathbf{P}) = \phi \nabla \boldsymbol{\cdot} \mathbf{P} + \mathbf{P} \boldsymbol{\cdot} \nabla \phi$ with $\mathbf{P}$ a vector-valued function.
    Let $ A \subseteq \Omega$ be an orientable domain and let $\hat{\mathbf{n}}$ be the outward unit normal on $\partial A$. By subsequently integrating the right-hand side of~\eqref{eqn:MA_pq_reformulated_a} over $A$ and applying Gauss's theorem we obtain
    \begin{align}
    \begin{split}
    \iint_A \nabla \boldsymbol{\cdot} \left(p \mathbf{J} \nabla q\right) \diff A = \oint_{\partial A} p \mathbf{J} \nabla q  \cdot \hat{\mathbf{n}} \diff s = \oint_{\partial A} p \nabla q  \cdot \hat{\pmb{\tau}} \diff s, 
    \end{split}
    \end{align}
    where we defined $\hat{\pmb{\tau}} = \mathbf{J}^\text{T} \hat{\mathbf{n}}$. Note that  $\hat{\pmb{\tau}}$ is the unit tangent vector to the domain taken in the counter clockwise direction.
    It follows from~\eqref{eqn:MA_pq_reformulated} that
    \begin{subequations}
    \label{eqn:MA_Integrals}
    \begin{align}
    \iint_A f^2 \diff A & = - \oint_{\partial A} p \nabla q  \cdot \hat{\pmb{\tau}} \diff s, \label{eqn:MA_IntegralsA} 
    \end{align}
    and in a similar way we find:
    \begin{align}
    \iint_A f^2 \diff A & = \oint_{\partial A} q \nabla p \cdot \hat{\pmb{\tau}} \diff s. \label{eqn:MA_IntegralsB} 
        \end{align}
    \end{subequations}
    Note that adding these equations yields
    \begin{align}
    \label{eqn:MA_Integral}
    \oint_{\partial A} (p \nabla q  + q \nabla p ) \cdot \hat{\pmb{\tau}} \diff s = 0,
    \end{align}
    which trivially holds by Stokes' theorem since $p \nabla q + q \nabla p = \nabla(pq)$ is conservative.
    Because the integral formulations~\eqref{eqn:MA_Integrals} are both equivalent to the \MAe{}, we use numerical approximations of~\eqref{eqn:MA_Integrals} as a measure for the residual of the numerical solution. To this purpose, let 
    \begin{subequations}
    \begin{gather}
    \begin{alignat}{5}
    \mathbf{H}_1 & = - & p \nabla q & =  - & p \begin{pmatrix}
    s \\ t
    \end{pmatrix} & = \frac{p f }{ a - b} \begin{pmatrix} a + b \\  -2  \end{pmatrix}, \\
    \mathbf{H}_2 & =  & q \nabla p & = & q \begin{pmatrix}
    r \\ s
    \end{pmatrix} & = \frac{q f }{ a - b} \begin{pmatrix} 2 a b \\ - a - b\end{pmatrix}.     
    \end{alignat}
    \end{gather}
    \end{subequations}
    Equations~\eqref{eqn:MA_Integrals} are equivalent to 
    \begin{align}
    \label{eqn:MA_Integral_HxHy}
    \iint_A f^2 \diff A & = \oint_{\partial A} \mathbf{H}_k \boldsymbol{\cdot} \hat{\pmb{\tau}} \diff s,
    \end{align}
    for $k = 1,2$.
    Choosing an appropriate domain $A$, this can therefore be used to determine the numerical residual.
    We choose the control volume $A = A_{i,j} = [x_{i-\frac{1}{2}}, x_{i+\frac{1}{2}}]\times[y_{j-\frac{1}{2}},y_{j+\frac{1}{2}}]$, and write~\eqref{eqn:MA_Integral_HxHy} as
    \begin{align}
    \label{eqn:MA_Integral_HxHy_split}
    I_k^N + I_k^S + I_k^W + I_k^E - \iint_{A_{i,j}} f^2 \diff A = 0,
    \end{align}
    where
    \begin{align}
    \label{eqn:MA_Integral_HxHy_splitIntegrals}
    \begin{split}
     I_k^N = - \int_{x_{i-\frac{1}{2}}}^{x_{i+\frac{1}{2}}} H_{k,x}(x, y_{j+\frac{1}{2}}) \diff x, & \hspace{20pt}
    I_k^S = \int_{x_{i-\frac{1}{2}}}^{x_{i+\frac{1}{2}}} H_{k,x}(x, y_{j-\frac{1}{2}}) \diff x, 
    \\
     I_k^W = - \int_{y_{j-\frac{1}{2}}}^{y_{j+\frac{1}{2}}} H_{k,y}(x_{i-\frac{1}{2}}, y) \diff y, & \hspace{20pt}
    I_k^E = \int_{y_{j-\frac{1}{2}}}^{y_{j+\frac{1}{2}}} H_{k,y}(x_{i+\frac{1}{2}}, y) \diff y,
    \end{split}
    \end{align}
    and $H_{k,x}$ denotes the $x$-component of $\mathbf{H}_k$, and the line integrals are carried out over the North, South, West and East part of the control volume $A_{i,j}$, as shown in Figure~\ref{fig:windDirectionsIntegrals}.
    
    We approximate the integrals $I_k^N$, $I_k^S$, $I_k^W$ and $I_k^E$ using Gauss-Legendre quadrature rules.
    To this end, let $z_1, z_2 \in \mathbb{R}$, $z_1 < z_2$ and $g\in C([z_1, z_2])$. Given the $n$ points $\xi_i$ and $n$ weights $w_i$ (see \cite{stroud1966gaussian} for example), we approximate
    \begin{align}
    \label{eqn:GaussquadRule}
    \int_{z_1}^{z_2} g(x) \diff x \approx \frac{z_2 - z_1}{2} \sum_{i=1}^{n} w_i g\Big(\frac{z_2 - z_1}{2} \xi_i + \frac{z_1 + z_2}{2}\Big).
    \end{align}
    In case $g$ is $2n$ times continuously differentiable, the error associated with~\eqref{eqn:GaussquadRule} is $\mathcal{O}((z_2 - z_1)^{2n+1})$ \cite[\S 5.2]{kahaner1988numerical}. Therefore, we choose $n = 3$ such that the calculation of the residual is asymptotically more accurate than the integration methods used for stepping from one grid line to the next. Note that we generally do not know $\mathbf{H}_k$ in the \mbox{points $(\tfrac{z_2 - z_1}{2} \xi_i + \tfrac{z_1 + z_2}{2})$}, where either $(z_1, z_2) = (x_{i-1/2}, x_{i+1/2})$ or $(z_1, z_2) = (y_{i-1/2}, y_{i+1/2})$. Therefore, we interpolate $\mathbf{H}_k$ using splines of order 5, such that the interpolation is at least as accurate as the step method. Similarly, we approximate $\iint_{A_{i,j}} f^2 \diff A$ by subsequently integrating over $x$ and $y$ using the Gauss-Legendre quadrature rule.
    We normalize the absolute residual of~\eqref{eqn:MA_Integral_HxHy_split} over $A_{i,j}$ by dividing it by the area $|A_{i,j}|$, and denote the result by $\epsilon_k(i,j)$. Lastly, we measure the residual over the whole grid by
    \begin{align}
    \label{eqn:epsilon_k}
    \epsilon_k = \max_{ \substack{ i \in \{2, \dots, N_x -1\} \\ j \in \{2, \dots, N_y -1 \}}} \epsilon_k(i,j).
    \end{align}

    \section{Numerical results}
    \label{sec:numericalResults}
    In this section we present numerical results for the \MAe{}. We will present results for the forward Euler, modified Euler and classic Runge-Kutta based methods for a default test case (Section~\ref{sec:defaultTestCase}), an example of which the analytical solution is known.
    For the modified Euler method we furthermore consider $2^\text{nd}$, $3^\text{rd}$ and $5^\text{th}$ order splines. 
    In Section~\ref{sec:aggregatedExample}, we compare the methods for a different neat example .
    Additionally, we show numerical results for the Runge-Kutta based method, for which we prescribe either two or zero boundary conditions per boundary segment (Section~\ref{sec:2InitStrips}), a case where the number of boundary conditions varies along the boundary (Section~\ref{sec:varyingBoundaryConditions}) and a case where one boundary value is nonsmooth (Section~\ref{sec:nonsmoothBoundary}). 
    
    \subsection{Default test case}
    \label{sec:defaultTestCase}
    To validate the numerical methods we design a default test case. To this end we let $\Omega = [0,1]\times[-0.5,0.5]$ and we calculate a right-hand side-solution pair $(f,u)$ using the method outlined in App.~\ref{app:generateSolutionPair} with $w(z) = \cos(i z)$, giving
    \begin{equation}
    \label{eqn:defaultTestCase}
    \begin{aligned}
        u(x,y) &= \cos(y) \cosh(x), \\
        f(x,y) &= \sqrt{\frac{\cos(2 y) + \cosh(2 x)}{2}}.
    \end{aligned}
    \end{equation}
    From the exact solution~\eqref{eqn:defaultTestCase} and~\eqref{eqn:introductionVariables_ab} we find $a$ and $b$ on the whole domain, viz.
    \begin{subequations}
    \label{eqn:defaultExampleCaseUpperLower}
    \begin{align}
    \label{eqn:defaultExampleCaseUpperLower_a}
    a(x,y) = & -\frac{\sin(y) \sinh(x) + f(x,y)}{\cos(y) \cosh(x)}, \\
    \label{eqn:defaultExampleCaseUpperLower_b}
    b(x,y) = & \frac{-\sin(y) \sinh(x) + f(x,y)}{\cos(y) \cosh(x)}.
    \end{align}
    \end{subequations}
    We impose the corresponding initial conditions
    \begin{align}
    \label{eqn:defaultExampleBasis}
    u(0,y) = \cos(y), \hspace{20pt} p(0,y) = 0.
    \end{align}
    By~\eqref{eqn:initialStripCalculations} we find
    \begin{equation}
    \label{eqn:defaultExampleCaseInitialStrip}
    \begin{aligned}
    q(0,y) & = -\sin(y), \hspace{5pt} & 
    t(0,y) & = -\cos(y), \hspace{5pt} &
    s(0,y) & = 0, \\
    a(0,y) & = -1, \hspace{5pt} &
    b(0,y) & = 1,  \hspace{5pt} &
    r(0,y) & = \cos(y).
    \end{aligned}
    \end{equation}
    To justify~\eqref{eqn:defaultExampleBasis}, note that the outward unit normal vector $\hat{\mathbf{n}}$ on the initial strip is $\hat{\mathbf{n}}(0,y) = (-1,0)^\text{T}$, $\mathbf{x}_\alpha(0,y) = (1,-1)^\text{T}$ and $\mathbf{x}_\beta(0,y) = (1,1)^\text{T}$. Therefore $\mathbf{x}_\alpha(0,y) \cdot \hat{\mathbf{n}}(0,y) = \mathbf{x}_\beta(0,y) \cdot \hat{\mathbf{n}}(0,y) = -1$, and hence we should prescribe $u$ and $p$ on the initial strip according to Case 2 in Section~\ref{sec:BoundaryConditionsHMA}.
    A similar calculations show that at $x = 1$ no boundary conditions need to be prescribed.
%
    On the upper boundary, we have $\hat{\mathbf{n}} = (0,1)^\text{T}$, $\mathbf{x}_\alpha = (1,a)^\text{T}$, $\mathbf{x}_\beta = (1,b)^\text{T}$, $a < 0$ and $b>0$, so that the $\alpha$-characteristic is entering the domain. Hence we need to prescribe $b$ and $a$ is known. In the same way, we need to prescribe $a$ at the lower boundary and $b$ is known.
    
    \subsubsection{Forward Euler based method}
    \label{sec:forwardEuler}
\begin{figure}[H]
        \centering
        \includegraphics[width = 0.45\linewidth]{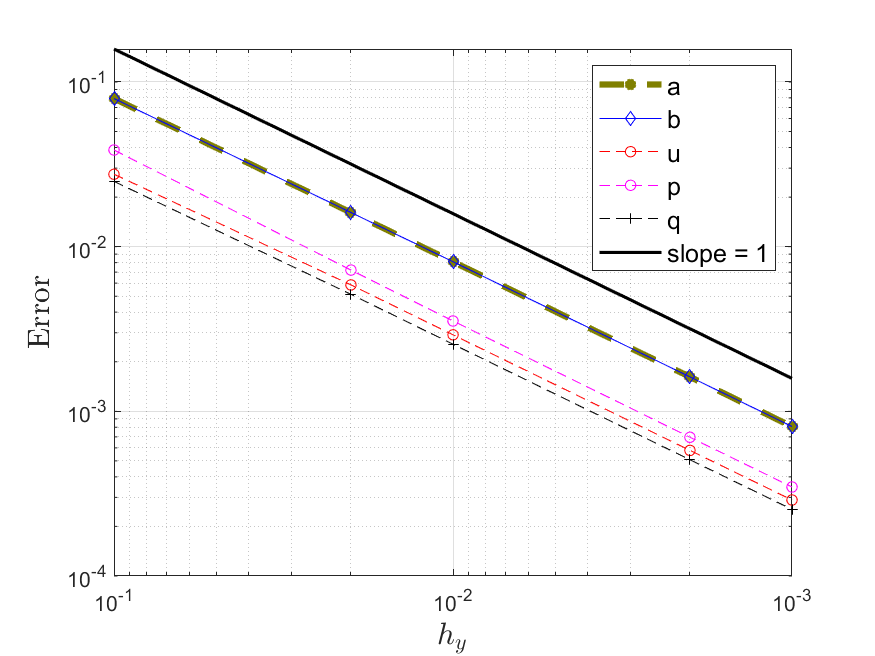}
        \hspace{20pt}
        \includegraphics[width = 0.45\linewidth]{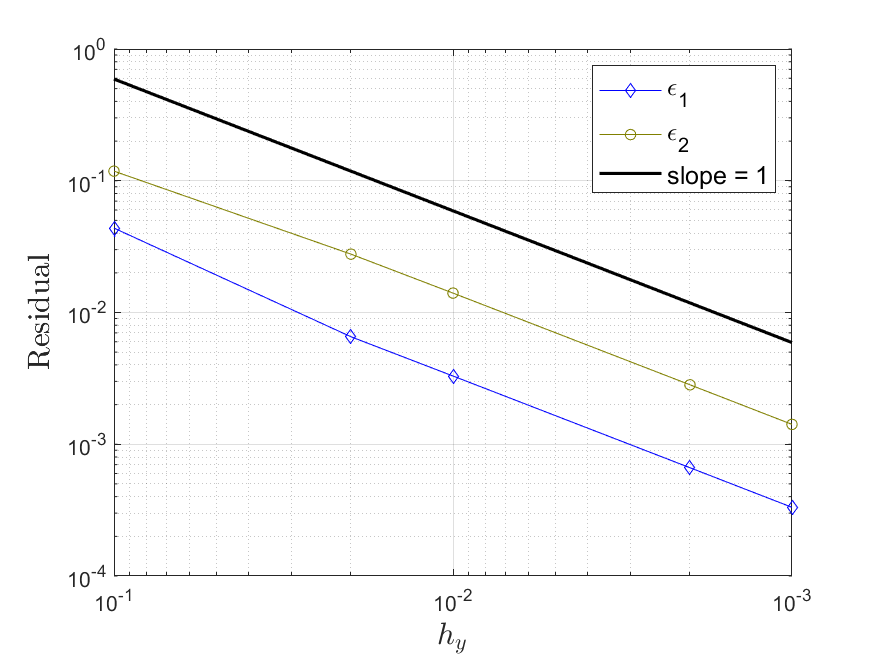}
        \caption{ Convergence of the global error (left) and the residual (right) for the forward Euler based method. }
        \label{fig:EE_global_res_convergence}
    \end{figure}

    We present the results for the forward Euler based method for which we use splines of first degree, i.e., linear interpolants.    
    The convergence of the forward Euler scheme is shown in Figure~\ref{fig:EE_global_res_convergence} as function of $h_y$, which controls $(h_x)_i$; see Section~\ref{sec:dynamicStepSize}.
    In the left figure the maximum absolute differences at $x=1$ between the function value and its numerical approximation for several variables are shown. More precisely, the figure shows 
    \begin{equation}
    E[a] := \max_{j = 1, \dots, N_y} \tfrac{1}{N_y} \left|a(x_{N_x}, y_j) - a_{N_x, j}\right|,
    \end{equation}
    and similar errors for $b, u, p$ and $q$ for varying $h_y$. 
    The dynamic step size control implies $(h_x)_i \approx \mathcal{O}(h_y)$, which allows us to quantify the error solely in terms of $h_y$.
    
    It is well known that the forward Euler method is locally second order accurate and globally first order, if the solutions are sufficiently smooth. Because the interpolation error and the local discretization error are both second order accurate, we expect the global convergence to be that of the forward Euler method, i.e., first order, which is also seen in the figure.

     The residuals are also shown and show first order convergence for $\epsilon_1$, $\epsilon_2$ defined by~\eqref{eqn:epsilon_k}. 
    To understand this, note that as we divide by the area of the control volume, that is we divide by $|A_{i,j}|$, by which we effectively normalize $\epsilon_1$, $\epsilon_2$ such that they convergence as the integrand does, which in this case is first order (for $\mathbf{H}_1$ and $\mathbf{H}_2$).
    
    
\begin{figure}[H]
        \centering
        \includegraphics[width = 0.45\linewidth]{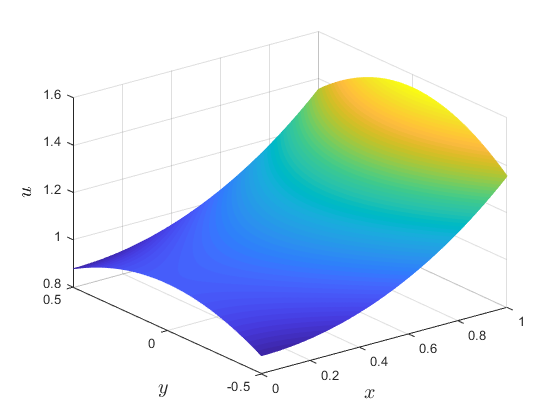}
        \hspace{20pt}
        \includegraphics[width = 0.45\linewidth]{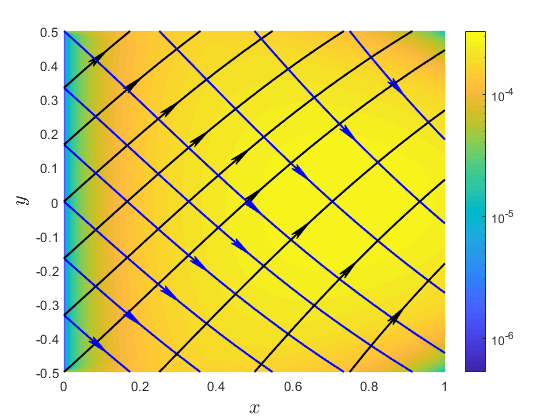}
        \caption{ Surface $u$ (left) and a color map of the residual $\epsilon_1$ with characteristics (right).}
        \label{fig:EE_surfaceAndIntRes}
    \end{figure}

    Figure~\ref{fig:EE_surfaceAndIntRes} shows both the solution surface $u$ (left), and a color map of the residual $\epsilon_1$ for the case $N_y = 1000$. The surface $u$ clearly is both smooth and a saddle surface. The right image shows the residual along with some characteristics. The shown characteristics are calculated after the simulation is done, and chosen such that they enter at 7 equidistant points on the initial strip and 5 on the upper and lower boundary. The direction of the characteristics clearly shows that both the blue characteristic, i.e., the $\alpha$-characteristic, and the black characteristic, i.e., the $\beta$-characteristic, enter at the initial strip. Hence, both $u$ and $p$ need to be prescribed at the initial strip. This is in agreement with~\eqref{eqn:defaultExampleBasis}. Furthermore, it shows that the $\alpha$-characteristics and $\beta$-characteristics leave the domain at the lower and upper boundary, respectively. Therefore $a$ and $b$ should be prescribed at the lower and upper boundary, respectively, agrees with the discussion on boundary conditions above.
    
    \subsubsection{Modified Euler based method}
    \label{sec:ME}
\begin{figure}[H]
        \centering
        \includegraphics[width = 0.45\linewidth]{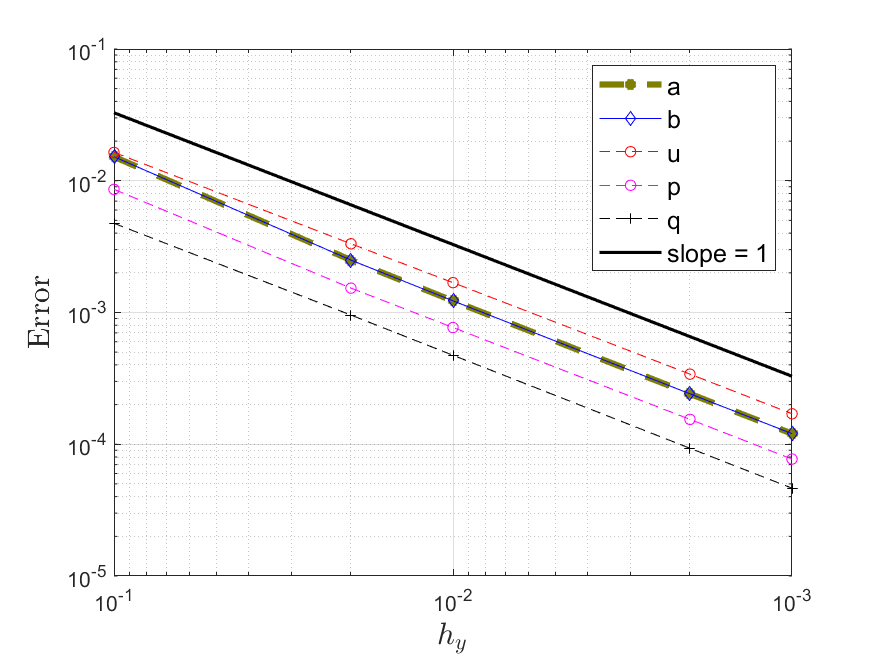}
        \hspace{20pt}
        \includegraphics[width = 0.45\linewidth]{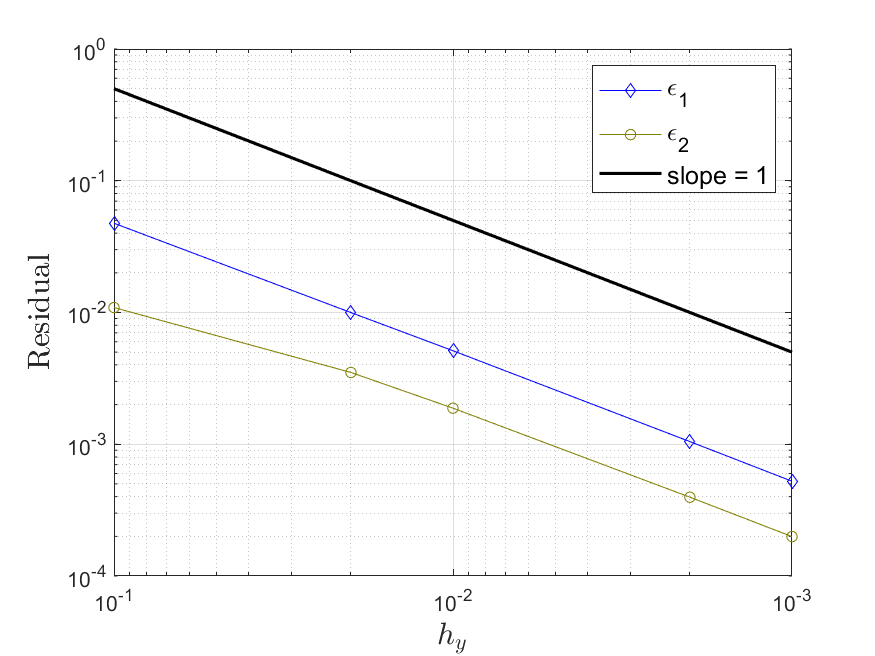}
        \caption{ Convergence for the global error (left) and the residual (right) for the modified Euler based method with a second order accurate interpolant. }
        \label{fig:ME_global_res_convergence_order2}
    \end{figure}

\begin{figure}[H]
    \centering
    \includegraphics[width = 0.45\linewidth]{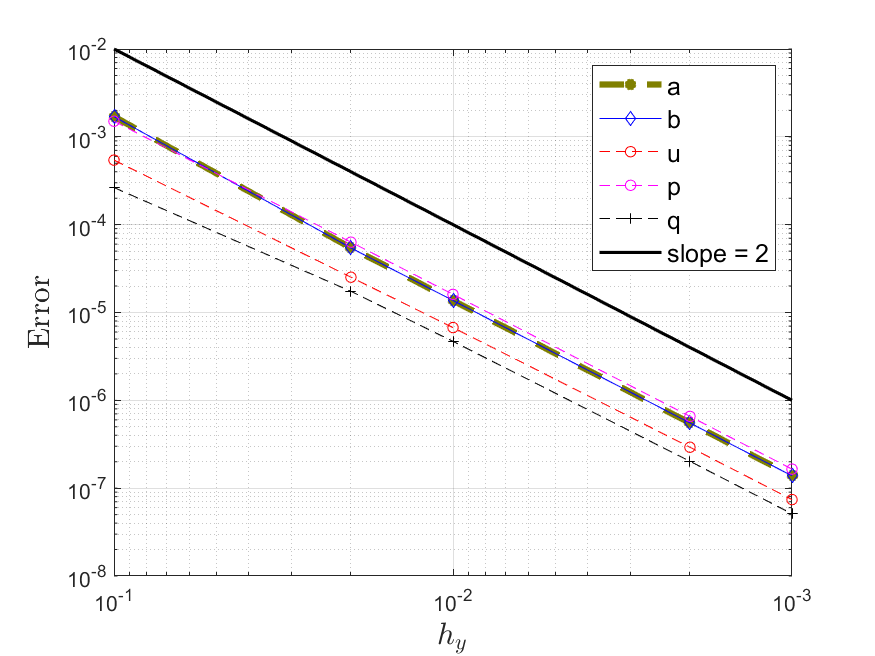}
    \hspace{20pt}
    \includegraphics[width = 0.45\linewidth]{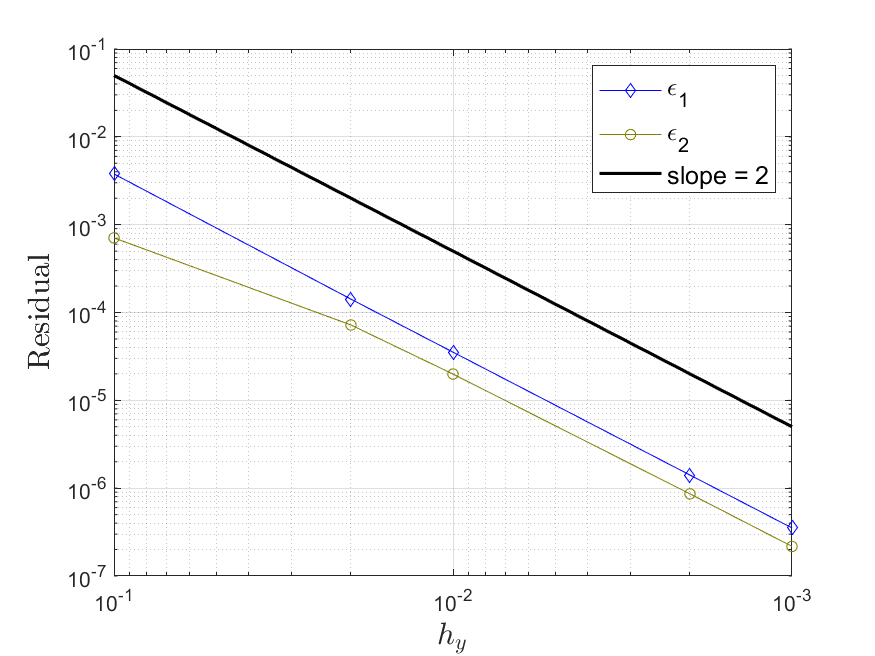}
    \caption{ Convergence for the global error (left) and the residual (right) for the modified Euler based method with a third order accurate interpolant. }
    \label{fig:ME_global_res_convergence_order3}
    \end{figure}

\begin{figure}[H]
    \centering
    \includegraphics[width = 0.45\linewidth]{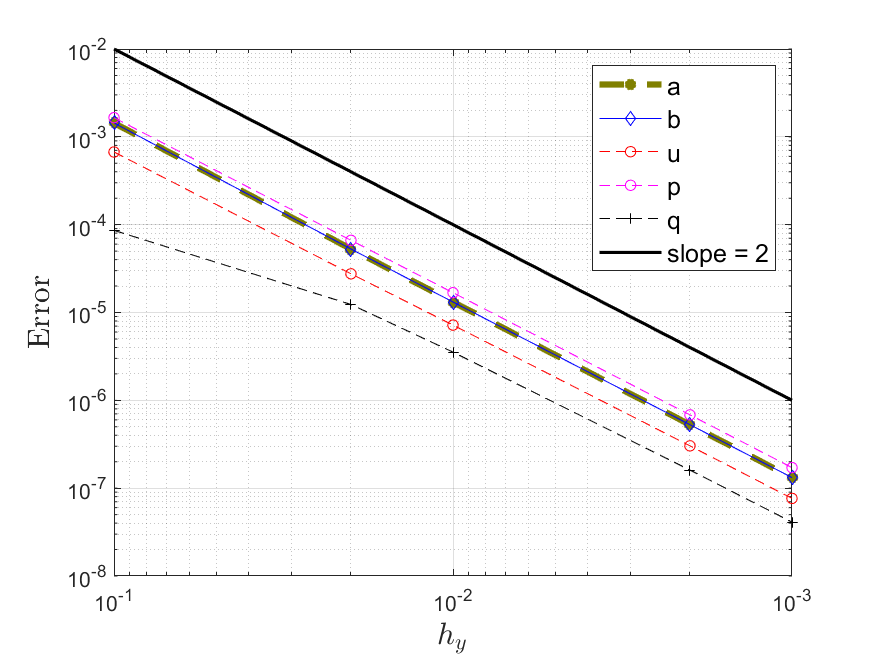}
    \hspace{20pt}
    \includegraphics[width = 0.45\linewidth]{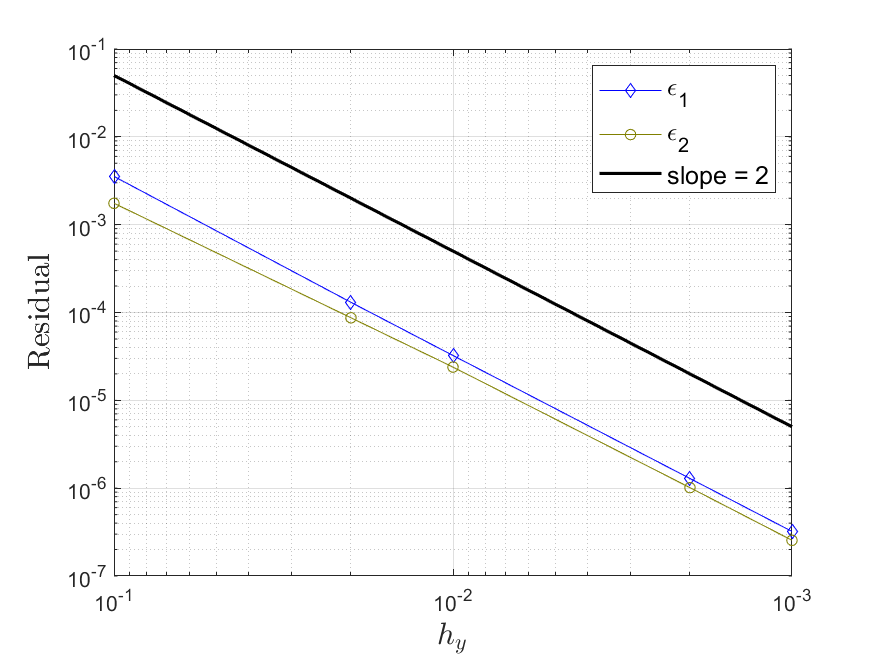}
    \caption{ Convergence for the global error (left) and the residual (right) for the modified Euler based method with a fifth order accurate interpolant. }
    \label{fig:ME_global_res_convergence_order5}
    \end{figure}

    In this section we will discuss a few results for the modified Euler based method. We demonstrate the importance of choosing an appropriate interpolation routine and show accompanying convergence results.
    Generally we use an interpolant which is as accurate as the integration routine because a more accurate interpolant will not increase the overall accuracy while being computationally more expensive, and a lower order interpolant will lower the convergence. In Figures~\ref{fig:ME_global_res_convergence_order2},~\ref{fig:ME_global_res_convergence_order3} and~\ref{fig:ME_global_res_convergence_order5} the convergence is shown for splines of order 2, 3, and 5, respectively. 
    Using splines of second order yields first order convergence for both the global error and the residual. This is in agreement with the expectation due to the local discretization error after interpolation being second order. Henceforth, second order splines, i.e., linear B-splines, reduce the rate of convergence, and higher order splines are preferred.
    
    Figures~\ref{fig:ME_global_res_convergence_order3} and ~\ref{fig:ME_global_res_convergence_order5} show the same order of convergence because the interpolants are at least as accurate as the local integration error. For a spline of order 3, i.e., quadratic B-splines, the interpolation error is as accurate as the local error of the modified Euler method. The accuracy of the global error is one order lower and equal to that of the residual. For the spline of order 5, i.e., for polynomials of fourth degree as B-splines, the local error of the modified Euler method is the limiting factor and the global error and the residual are second order accurate.

    \subsubsection{Classic Runge-Kutta based method}
    \label{sec:RK4}
\begin{figure}[H]
        \centering
        \includegraphics[width = 0.45\linewidth]{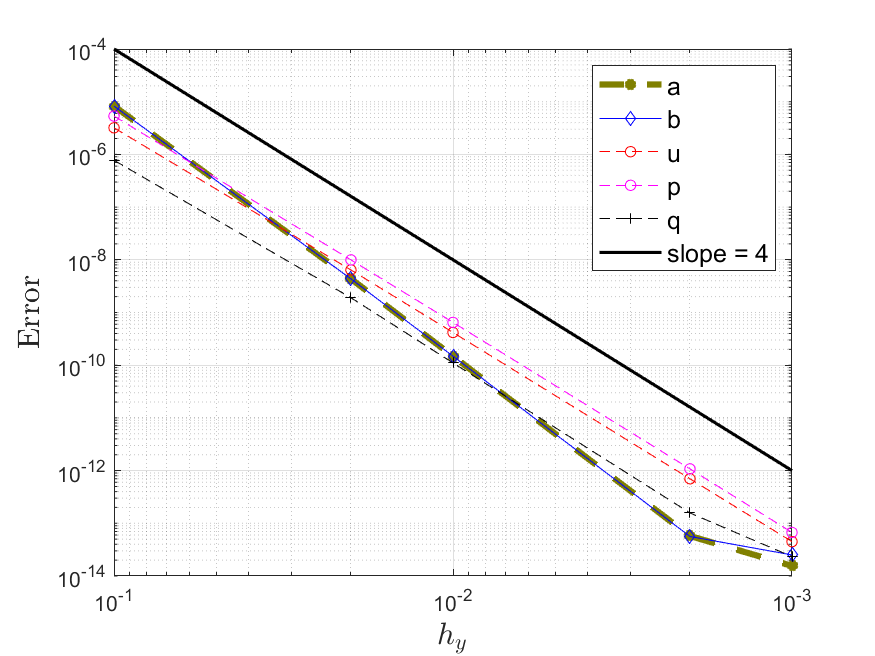}
        \hspace{20pt}
        \includegraphics[width = 0.45\linewidth]{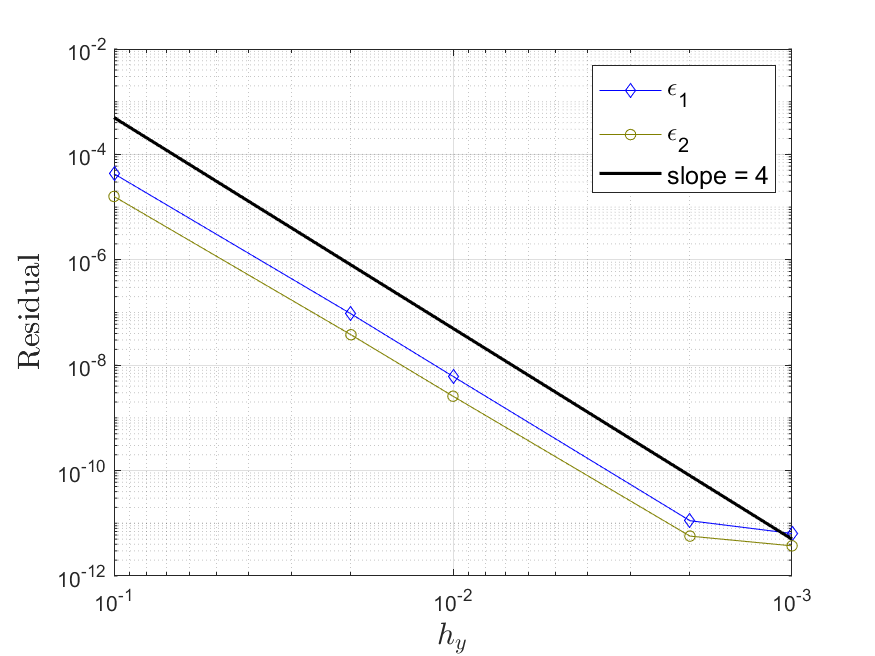}
        \caption{ Convergence for the global error (left) and the residual (right) for the classic Runge-Kutta based method with a fifth order accurate interpolant. }
        \label{fig:RK_global_res_convergence_order5}
    \end{figure}

    Analogously to the previous sections, we consider the default test case. Figure~\ref{fig:RK_global_res_convergence_order5} shows the results for using the Runge-Kutta method with spline interpolants of order 5. Because the Runge-Kutta method is locally fifth order accurate, which coincides with the accuracy of the splines, we expect a fourth order global convergence. This is indeed shown in the figure. The convergence of the residuals is also expected to be of order 4, as also seen in the figure. The convergence seems to slow down for $h_y \approx 1/1000$, which is due to round-off errors as the solutions reaches the used computer precision.
    
    \subsection{An aggregated example}
    \label{sec:aggregatedExample}
\begin{figure}[H]
        \centering
        \includegraphics[width = 0.45\linewidth]{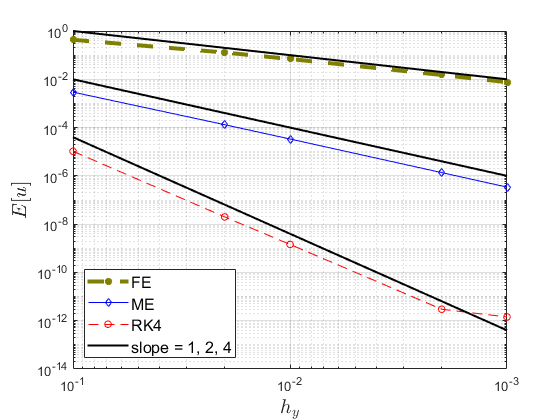}
        \hspace{20pt}
        \includegraphics[width = 0.45\linewidth]{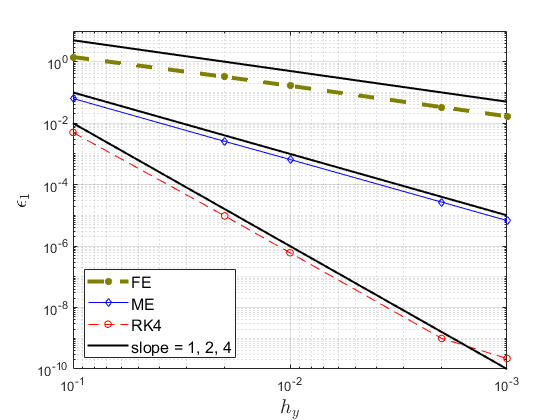}
        \caption{ Convergence for the global error (left) and the residual (right). }
        \label{fig:compareMethods}
    \end{figure}

\begin{figure}[H]
    \centering
    \includegraphics[width = 0.7\linewidth]{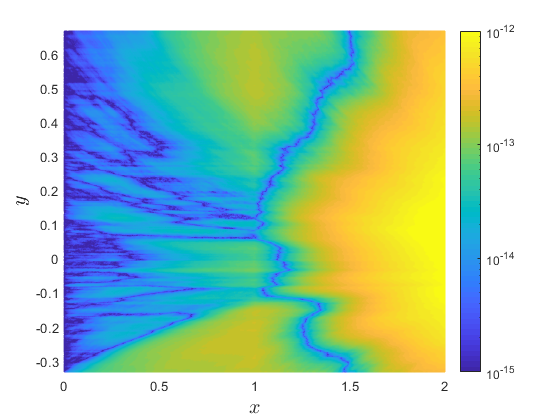}
    \caption{ The error $|u_{i,j} - u(\mathbf{x}_{i,j})|$ over the domain for $h_y = 10^{-3}$.}
    \label{fig:uError}
    \end{figure}

    Next we compare the developed numerical methods for the example
    \begin{align}
    \label{eqn:ufCase5}
    u(x,y) = e^x \cos(y), \hspace{20pt} f(x,y) = e^x,
    \end{align}
    which is constructed using $w(z) = e^z$; see App.~\ref{app:generateSolutionPair}.
    Let $\Omega = [0, 2] \times [-1/3, 2/3]$ be the computational domain. A straightforward calculation using~\eqref{eqn:ufCase5} shows
    \begin{equation}
    \label{eqn:combinedInitialStrip}
    \begin{aligned}
        p(x,y) & = e^x \cos(y), \hspace{5pt} & 
        q(x,y) & = -e^x \sin(y)\\
        a(x,y) & = -\frac{\sin(y)+1}{\cos(y)}, \hspace{5pt} &
        b(x,y) & = \frac{-\sin(y)+1}{\cos(y)}, 
    \end{aligned}
    \end{equation}
    which we use, along with~\eqref{eqn:ufCase5}, to prescribe $u, p, q, a, b$ on the initial strip $x=0$, $-\tfrac{1}{3} \leq y \leq \tfrac{2}{3}$, accordingly. Equations~\eqref{eqn:combinedInitialStrip} show $a < 0$, $b >0$ on $\Omega$. Analogously to the default test case, we prescribe $a$ on the lower and $b$ on the upper boundary as dictated by~\eqref{eqn:combinedInitialStrip}.
    We compare the forward Euler, modified Euler and Runge-Kutta methods using second, third and fifth order splines, respectively. The results are shown in Figure~\ref{fig:compareMethods}. The left figure shows first, second and fourth order convergence of the global error of $u_{i,j}$ for the forward Euler, modified Euler and Runge-Kutta-method, respectively. These rates of convergence are in agreement with the previous sections.  The convergence of $\epsilon_1$, shown in the right figure, also shows first, second and fourth order convergence. Furthermore, the convergence of the global error and residual stagnates for the Runge-Kutta-method at $h_y\approx 10^{-3}$. Figure~\ref{fig:uError} shows that the error accumulates for increasing $x$, i.e., the further in the domain, as measured from the initial strip, the higher the error. Due to the error being of computer precision near $x = 0$ for $h_y = 10^{-3}$, this accumulation of errors bounds the global error $E[u]$ from below. Closed form bounds for $E[u]$ are not known, though it is evident that it depends on the computational domain and the boundary conditions prescribed.

    \subsection{An initial strip extended over two edges}
    \label{sec:2InitStrips}
    We will demonstrate an example for the forward Euler based method for which we have two entering characteristics at both the Western and Northern boundary, and two leaving characteristics on both the Southern and Eastern boundary. In this case the Northern boundary is also an initial strip as discussed in Section~\ref{sec:BoundaryConditionsHMA}.
    To this end let
    \begin{align}
    \label{eqn:ufCase23}
    u(x,y) = x^3 y^2 + 1, \hspace{20pt} f(x,y) = 2 \sqrt{6} x^2 y,
    \end{align}
    on the domain $\Omega = [1,2]^2$. 
    A straightforward calculation shows that~\eqref{eqn:ufCase23} satisfies the \MAe{}~\eqref{eqn:MASimple}, and that
    \begin{align}
    a(x,y) = \frac{(-3 + \sqrt{6}) y}{x}, \hspace{20pt} 
    b(x,y) = -\frac{(3 + \sqrt{6}) y}{x},
    \end{align}
    which implies $a,b < 0$ on $\Omega$. Let $\hat{\mathbf{n}}_\text{W}$ denote the normal at the Western boundary segment, and likewise for the other subscripts. It follows that
    \begin{equation}
    \begin{aligned}
    \mathbf{x}_\alpha(1,y) \boldsymbol{\cdot} \hat{\mathbf{n}}_\text{W} < 0, \hspace{10pt} & \hspace{10pt} \mathbf{x}_\beta(1,y) \boldsymbol{\cdot} \hat{\mathbf{n}}_\text{W} < 0, \\
    \mathbf{x}_\alpha(x,1) \boldsymbol{\cdot} \hat{\mathbf{n}}_\text{S} > 0, \hspace{10pt} & \hspace{10pt} \mathbf{x}_\beta(x,1) \boldsymbol{\cdot} \hat{\mathbf{n}}_\text{S} > 0, \\
    \mathbf{x}_\alpha(2,y) \boldsymbol{\cdot} \hat{\mathbf{n}}_\text{N} > 0, \hspace{10pt} & \hspace{10pt} \mathbf{x}_\beta(2,y) \boldsymbol{\cdot} \hat{\mathbf{n}}_\text{N} > 0, \\ 
    \mathbf{x}_\alpha(x,2) \boldsymbol{\cdot} \hat{\mathbf{n}}_\text{E}< 0, \hspace{10pt} & \hspace{10pt} \mathbf{x}_\beta(x,2) \boldsymbol{\cdot} \hat{\mathbf{n}}_\text{E} < 0. 
    \end{aligned}
    \end{equation} 
    The classification of the boundary conditions in Section~\ref{sec:BoundaryConditionsHMA}, imply we should prescribe two boundary conditions at the Western and Northern boundaries and zero boundary conditions at the Eastern and Southern segments. By prescribing two boundary conditions to the Northern boundary, it is an initial strip. 
        The total set of prescribed boundary conditions thus read
    \begin{align}
    &     u(1,y) = y^2 + 1, 
    &&    p(1,y) = 3 y^2,
    &&    u(x,2) = 4x^3 + 1,
    &&    q(x,2) = 4 x^3.
    &        
    \end{align}
    From $u(1,y)$ and $p(1,y)$ we obtain $q,r,s,t,a,b$ at the Western boundary by~\eqref{eqn:initialStripCalculations}. Analogously, $u(x,2)$ and $q(x,2)$ determine $p,a,b,r,s,t$ at the northern boundary.

\begin{figure}[H]
        \centering
        \includegraphics[width = 0.45\linewidth]{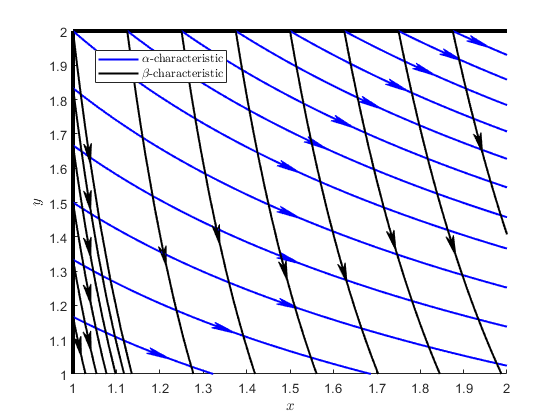}
        \hspace{20pt}
        \includegraphics[width = 0.45\linewidth]{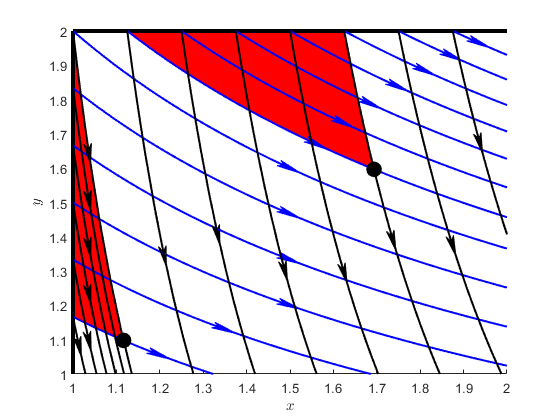}
        \caption{ Domain $\Omega$ and some characteristics, which enter the domain at the Western and Northern boundaries, with the domain of dependence for the points $(1.7, 1.6)$ and $(1.1, 1.1)$. }
        \label{fig:boundaryProblemCharacteristics}
    \end{figure}

    Figure~\ref{fig:boundaryProblemCharacteristics} shows some characteristics for this example, which nicely demonstrates where characteristics are entering or leaving the domain. Furthermore, the figure shows that $\Omega$ is fully covered by the domain of dependence of the two initial strips. The figure on the right also shows the points $(1.7, 1.6)$ and $(1.1, 1.1)$ with their corresponding domain of dependence colored red.

\begin{figure}[H]
        \centering
        \includegraphics[width = 0.45\linewidth]{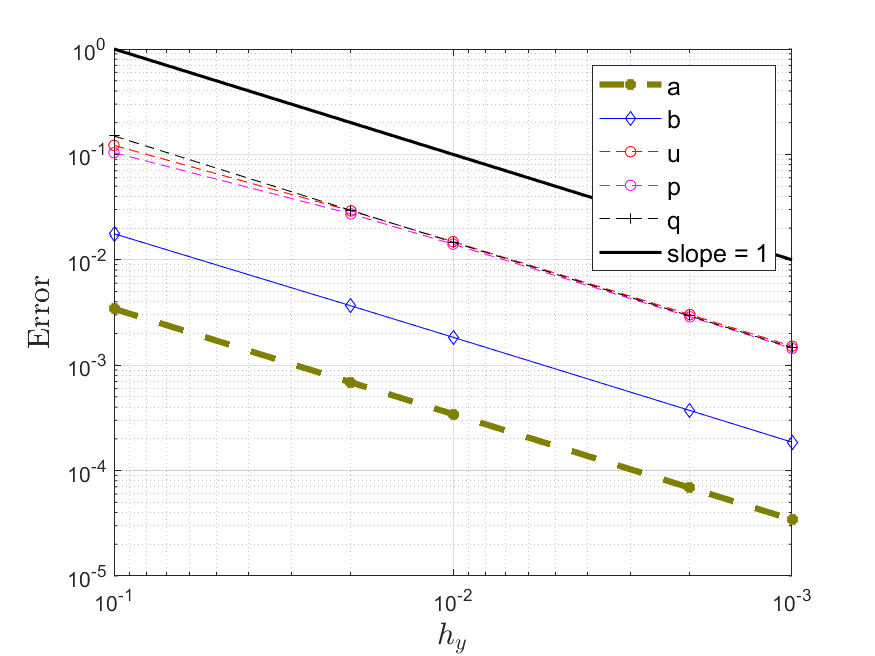}
        \hspace{20pt}
        \includegraphics[width = 0.45\linewidth]{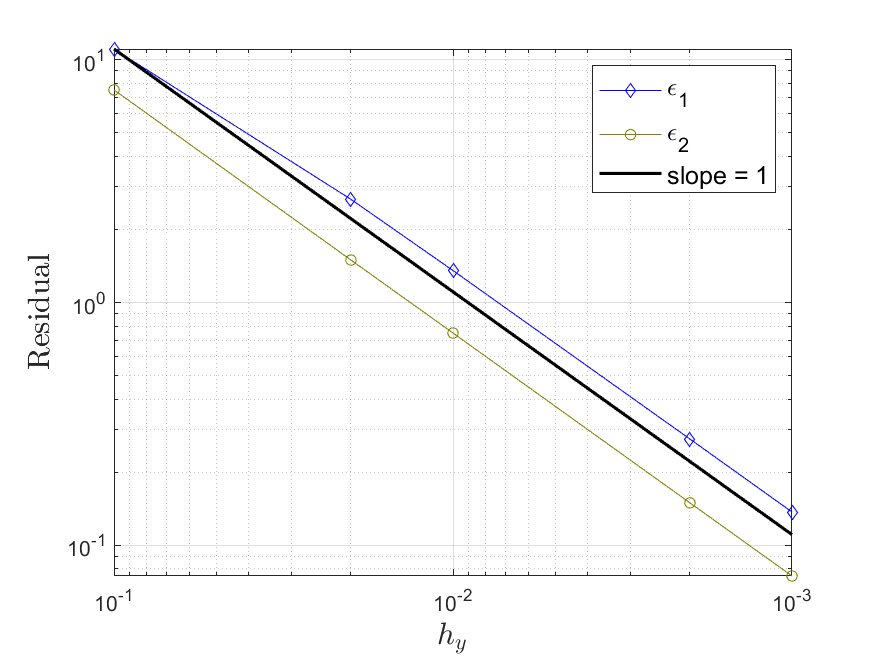}
        \caption{ Convergence for the global error (left) and the residual (right) for the forward Euler based method with a second order accurate interpolant. }
        \label{fig:Case23}
    \end{figure}

\begin{figure}[H]
    \centering
    \includegraphics[width = 0.45\linewidth]{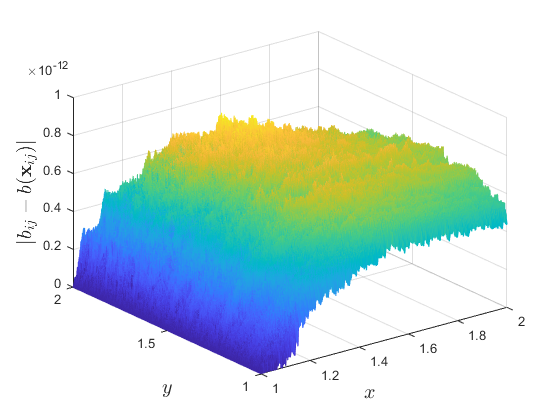}
    \hspace{20pt}
    \includegraphics[width = 0.45\linewidth]{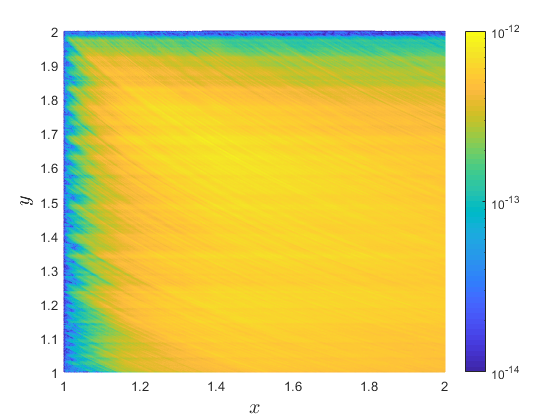}
    \caption{ The numerical error $|b_{i,j} - b(\mathbf{x}_{i,j})|$. }
    \label{fig:Case23ErrorInB}
    \end{figure}

    Figure~\ref{fig:Case23} shows the convergence for the Euler based method for this example, for both the global error and the residual. As expected, both show first order convergence.
    Figure~\ref{fig:Case23ErrorInB} shows the error in $b$ for this example, calculated using the he Runge-Kutta method with fifth order splines with $h_y = 1/1000$.
    The figure shows the accumulation of numerical errors over the domain and shows $b$ is most accurate near the boundaries where both $a$ and $b$ are prescribed.
    
    \subsection{Varying number of boundary conditions}
    \label{sec:varyingBoundaryConditions}
\begin{figure}[H]
    \centering
    \captionsetup{width=0.8\linewidth}
    \begin{tikzpicture}[scale=1.5]
    \draw [thin,-latex, ->] (-1.2, -1.2) -- (-2/3, -1.2) node [right] {$x$};
    \draw [thin,-latex, ->] (-1.2, -1.2) -- (-1.2, -2/3) node [above] {$y$};

    \draw [-] (0,2) |- (0,0) |- (6,0) |- (6,2) |- (0,2);

    \draw [dashed] (1,0) |- (-1/2,0) node (ymin) [left] {-2};
    \draw [dashed] (0,2) |- (-1/2,2) node (ymax) [left] {-1.5};
    \draw [dashed] (0,0) -- (0,-1/2) node (xmin) [below] {1};
    \draw [dashed] (6,0) -- (6,-1/2) node (xmax) [below] {2.5};
    
    \draw [dashed] (4,0) -- (4,-1/2) node (xmax) [below] {2};
    \draw [dashed] (2,2) -- (2,-1/2) node (xmax) [below] {1.5};

    \draw [line width=0.8mm] (0,0) |- (0,2) {};
    \draw [line width=1.5mm, dashed] (0,0) |- (0,2) {};
    
    \draw [line width=0.8mm] (0,2) |- (6,2) {};
    \draw [line width=1.5mm, dashed] (0,2) |- (6/3,2) {};

    \draw [line width=1.5mm, dashed] (12/3,0) |- (6,0) {};
    
    \draw [-] (6.5, 2) |- (7,2) node [right] {No b.c.};
    \draw [line width=0.8mm]  (6.5, 1.7) |- (7,1.7) node [right] {$a$};
    \draw [line width=1.5mm, dashed] (6.5, 1.4) |- (7,1.4) node [right] {$b$};
    
    \end{tikzpicture}
    \caption{ Schematic overview of the prescribed boundary conditions and their locations. }
    \label{fig:boundaryConditionsVarying}
    \end{figure}
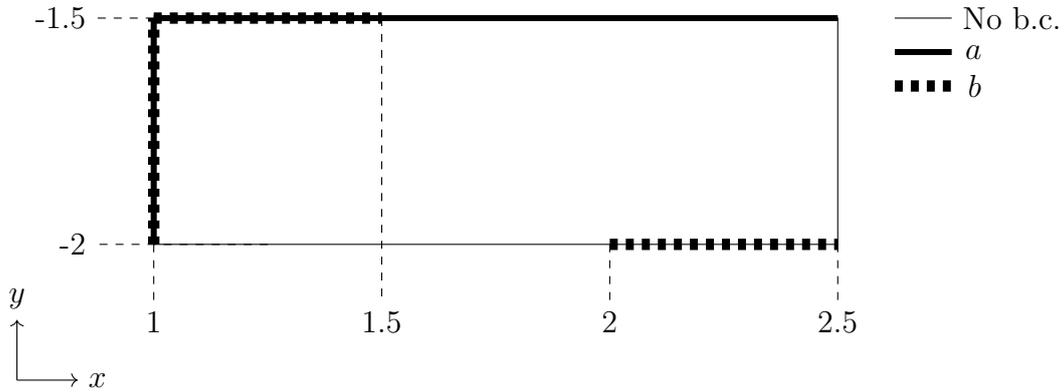

    Next we will show an example for which the number of boundary conditions we prescribe changes along the boundary. To this end let $\Omega = [1, 2.5] \times [-2, -1.5]$ and
    \begin{align}
    u(x,y) = 1 + e^{2 y / x}, \hspace{20pt} f(x,y) = \frac{2}{x^2} e^{2 y / x}.
    \end{align}
    It follows that
    \begin{align}
    a(x,y) = 1 + \frac{y}{x}, \hspace{20pt} b(x,y) = \frac{y}{x},
    \end{align}
    which implies $a \leq 0$ for $x \geq -y$, $a > 0$ for $x < - y$ and $b < 0$ on $\Omega$. Let $\hat{\mathbf{n}}$ be the outward unit normal, then the boundary conditions to be prescribed are
    \begin{equation*}
    \setlength{\arraycolsep}{0pt}%
    \renewcommand{\arraystretch}{1.2}%
    \begin{array}{ *{4}llll }
    \text{Boundary segment} & \multicolumn{2}{l}{\hspace{40pt} \text{Classification}} &\hspace{40pt} \text{Boundary condition(s)} \\
    \hline
       \text{Western}, & \hspace{40pt} (\mathbf{x}_\alpha, \hat{\mathbf{n}}) < 0, & \hspace{7pt} (\mathbf{x}_\beta, \hat{\mathbf{n}}) < 0, & \hspace{40pt} a \text{ and } b, \\
        \text{Southern, left},       & \hspace{40pt} (\mathbf{x}_\alpha, \hat{\mathbf{n}}) > 0, & \hspace{7pt} (\mathbf{x}_\beta, \hat{\mathbf{n}}) > 0, & \hspace{40pt} \text{None}, \\
        \text{Southern, right},  & \hspace{40pt} (\mathbf{x}_\alpha, \hat{\mathbf{n}}) < 0, & \hspace{7pt} (\mathbf{x}_\beta, \hat{\mathbf{n}}) > 0, & \hspace{40pt} b, \\
        \text{Northern, left},     & \hspace{40pt} (\mathbf{x}_\alpha, \hat{\mathbf{n}}) < 0, & \hspace{7pt} (\mathbf{x}_\beta, \hat{\mathbf{n}}) < 0, & \hspace{40pt} a \text{ and } b,\\
        \text{Northern, right},  & \hspace{40pt} (\mathbf{x}_\alpha, \hat{\mathbf{n}}) > 0, & \hspace{7pt} (\mathbf{x}_\beta, \hat{\mathbf{n}}) < 0, & \hspace{40pt} a, \\
        \text{Eastern}, & \hspace{40pt} (\mathbf{x}_\alpha, \hat{\mathbf{n}}) > 0, & \hspace{7pt} (\mathbf{x}_\beta, \hat{\mathbf{n}}) > 0, & \hspace{40pt} \text{None,}
    \end{array}
    \end{equation*}
    as illustrated in Figure~\ref{fig:boundaryConditionsVarying}.
    
\begin{figure}[H]
        \centering
        \includegraphics[width = 0.45\linewidth]{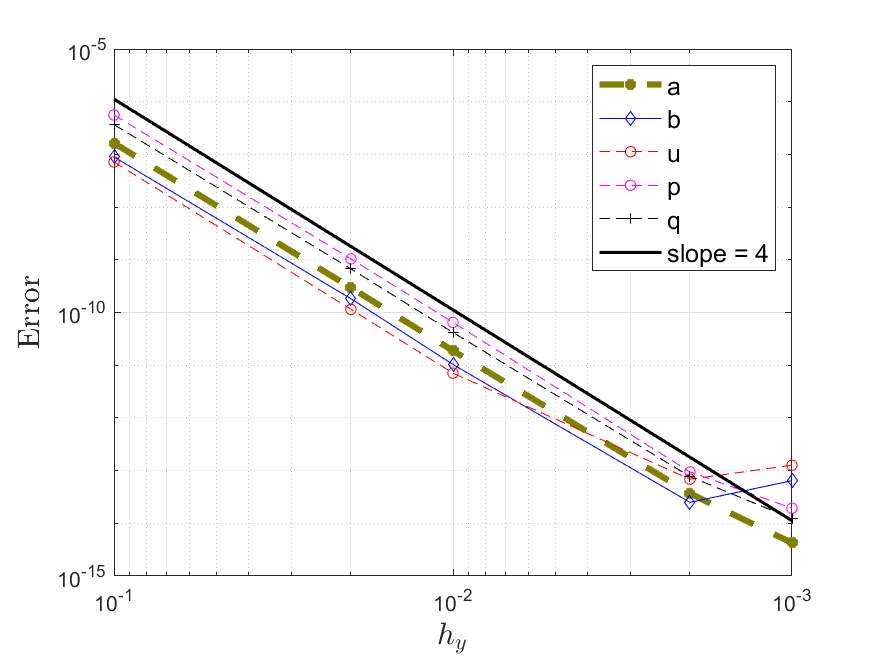}
        \hspace{20pt}
        \includegraphics[width = 0.45\linewidth]{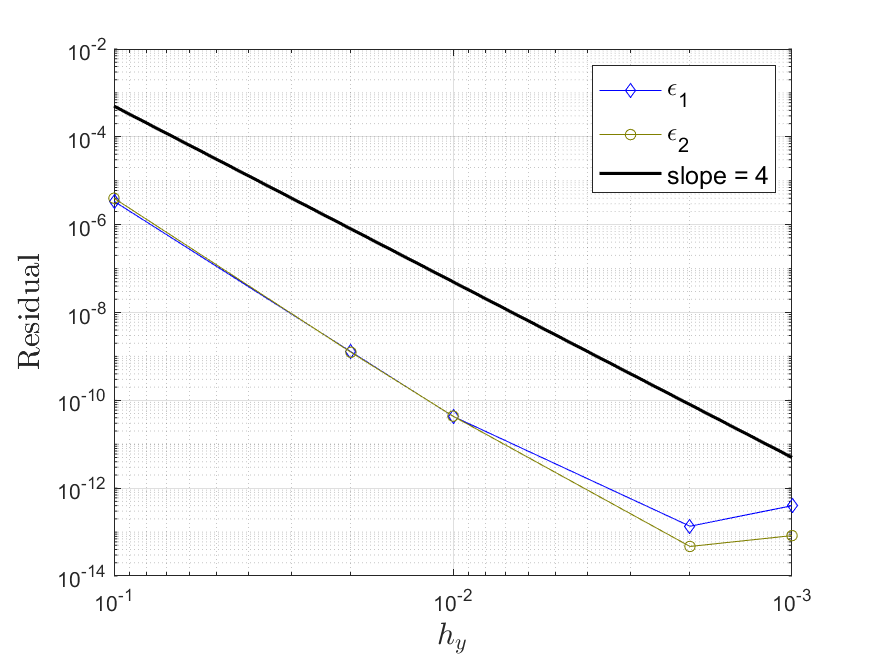}
        \caption{ Convergence of the global error (left) and the residual (right).}
        \label{fig:ChangingBCs}
    \end{figure}

\begin{figure}[H]
        \centering
        \includegraphics[width = 0.45\linewidth]{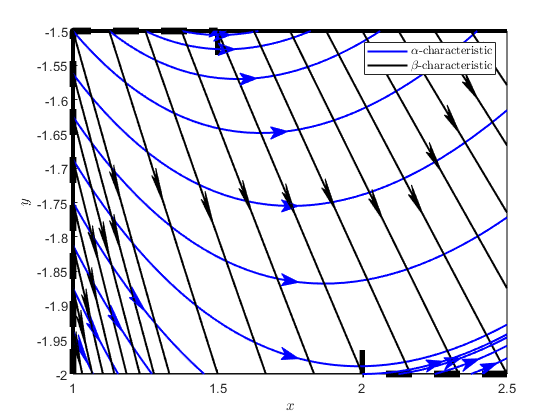}
        \hspace{20pt}
        \includegraphics[width = 0.45\linewidth]{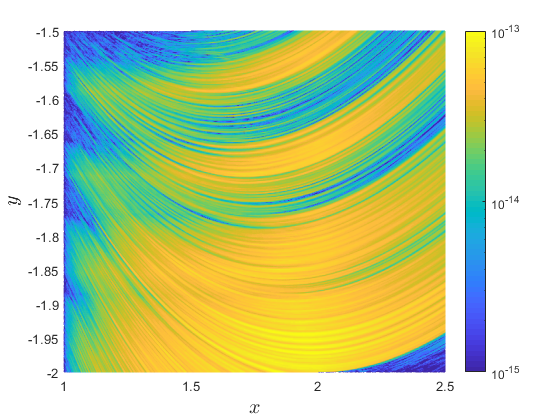}
        \caption{ Characteristics in the domain (left) and a heatmap of the numerical error in $b$ (right).}
        \label{fig:ChangingBCs2}
    \end{figure}

    Figure~\ref{fig:ChangingBCs} shows the convergence of the global error and the residual for the Runge-Kutta based method with fifth order splines. The convergence is fourth order as expected, and slowly comes to a halt for a fine grid, as also discussed in the previous section. Figure~\ref{fig:ChangingBCs2} shows the characteristics in the domain (left), and a heatmap of the error $|b_{i,j} - b(\mathbf{x}_{i,j})|$ (right). The heatmap clearly shows the swirling influence of the $\alpha$-characteristics, and a lower error near the segments of the boundary where $b$ is prescribed.
    
    \subsection{A problem with nonsmooth boundary conditions}
    \label{sec:nonsmoothBoundary}
    The last example is similar to the default test case of Section~\ref{sec:defaultTestCase}, but now with nonsmooth boundary conditions. As before, we set $\Omega = [0, 1] \times [-0.5, 0.5]$. The initial conditions are the same as before, given by~\eqref{eqn:defaultTestCase} and~\eqref{eqn:defaultExampleBasis}. At the lower and upper boundary, we now use
    \begin{subequations}
    \label{eqn:a_nonsmooth_prescribed}
    \begin{align}
    \label{eqn:a_nonsmooth_prescribed_a}
    a(x,-0.5) & = - e^{-3 x / 2} (x^2 + 1), \\
    \label{eqn:a_nonsmooth_prescribed_b}
    b(x,0.5) & = \frac{-\sin(0.5) \sinh(x) + f(x,0.5)}{\cos(0.5)\cosh(x)},
    \end{align}
    \end{subequations}
    respectively, where purposely $a(x, -0.5)$ as given by \ref{eqn:a_nonsmooth_prescribed}, does not correspond to $a$ of the default test case,~\eqref{eqn:defaultExampleCaseUpperLower_a}. The prescribed $b$ is continuous and smooth while the prescribed $a$ is continuous but nonsmooth at the point $(0,-0.5)$, i.e.,
    \begin{align}
    \odiff{b(x, 0.5)}{x} \Big|_{x = 0} & = \odiff{b(0,y)}{y} \Big|_{y = -0.5}, \\
    \odiff{a(x, -0.5)}{x} \Big|_{x = 0} & \neq \odiff{a(0,y)}{y} \Big|_{y = -0.5}.
    \end{align}
    The latter statement can be established by applying the identities
    \begin{align}
    a_x = \frac{-s_x + f_x}{t} - \frac{-s+f}{t^2} t_x = \frac{-s_x + f_x - a t_x}{t}, \hspace{20pt} s_x = r_y, \hspace{20pt} t_x = s_y,
    \end{align}
    to the initial strip~\eqref{eqn:defaultExampleCaseInitialStrip} to obtain $a_x(0, -0.5)  \approx 0.546$,
    while the derivative of $a$ as given by~\eqref{eqn:a_nonsmooth_prescribed} is $a_x(0, -0.5) \approx - 0.878$.
    Hence $a$ is nonsmooth. Generally the error terms of our numerical methods depend on derivatives of the functions to be estimated. As $a$ is nonsmooth, we do not necessarily expect convergence as we did before. Furthermore, no analytical solution is known for this particular example, therefore we base convergence on the residual values. 
    
\begin{figure}[H]
        \centering
        \includegraphics[width = 0.45\linewidth]{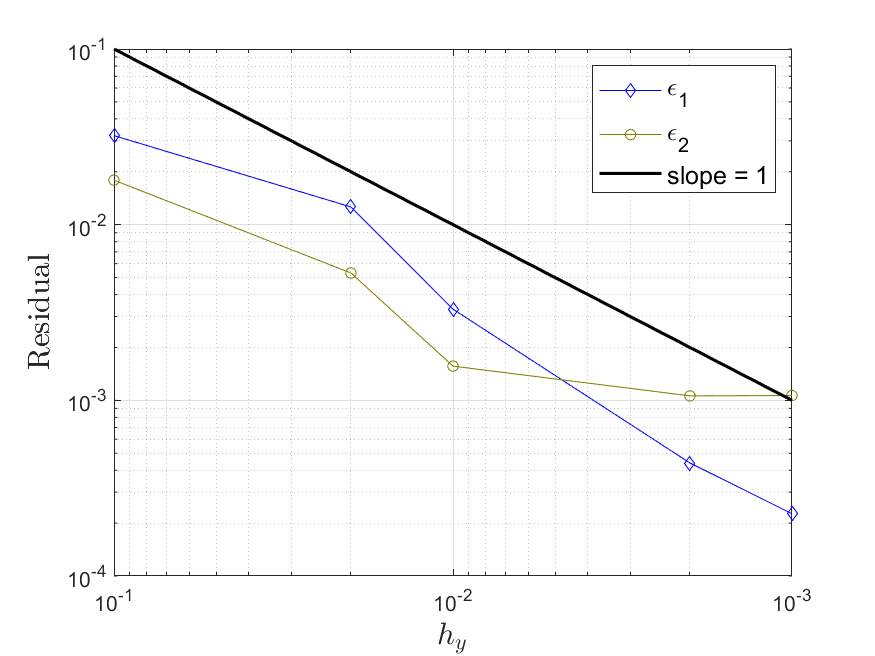}
        \hspace{20pt}
        \includegraphics[width = 0.45\linewidth]{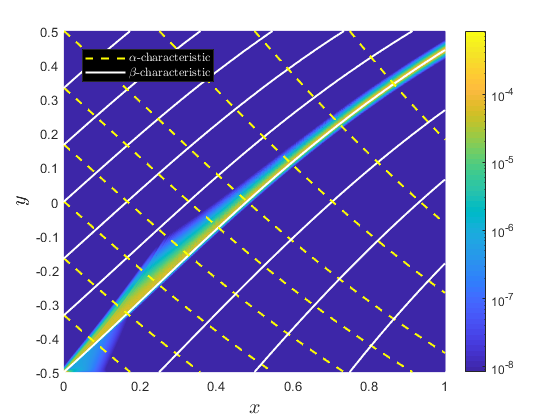}
        \caption{ Convergence of the residual (left) and a heat map of the residual $\epsilon_1$ (right) for nonsmooth boundary value $a$.}
        \label{fig:discontinuousBoundary}
    \end{figure}

\begin{figure}[H]
    \centering
    \includegraphics[width = 0.45\linewidth]{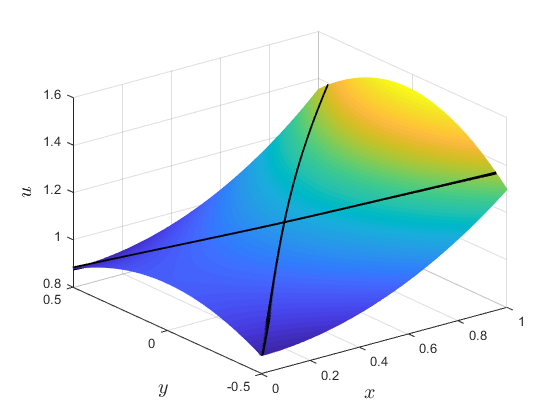}
    \hspace{20pt}
    \includegraphics[width = 0.45\linewidth]{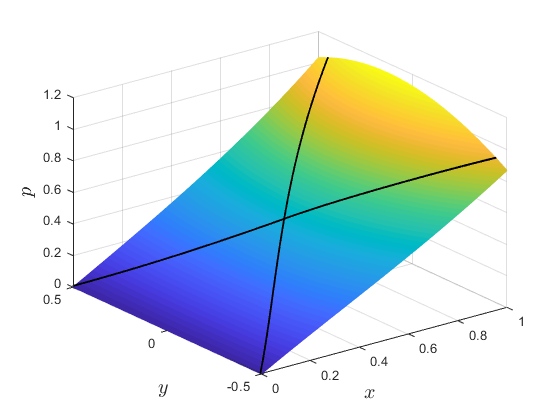}
    \caption{ Smooth solutions $u$ (left) and $p$ (right) for nonsmooth boundary data.}
    \label{fig:discontinuousBoundarySurface}
    \end{figure}

\begin{figure}[H]
    \centering
    \includegraphics[width = 0.45\linewidth]{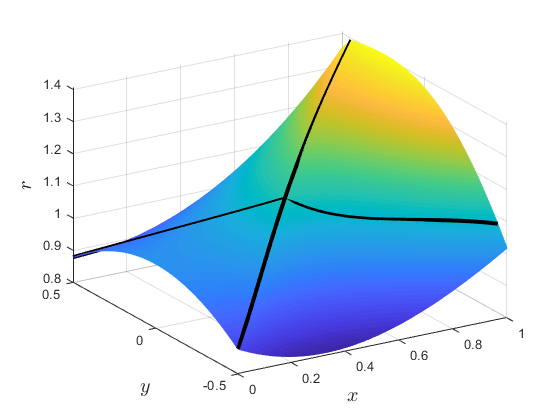}
    \hspace{20pt}
    \includegraphics[width = 0.45\linewidth]{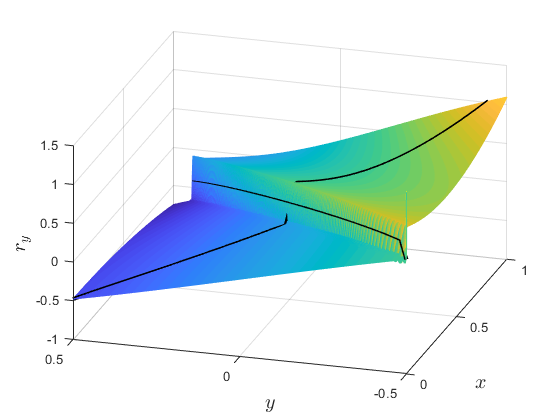}
    \caption{ Nonsmooth solutions $r$ (left) and its derivative $r_y$ (right) for nonsmooth boundary data.}
    \label{fig:discontinuousBoundarySurface_s_tx}
    \end{figure}

    We use the Runge-Kutta method with fifth order splines for this example.
    Figure~\ref{fig:discontinuousBoundary} shows the convergence of the residual (left), and a heat map of $\epsilon_1$, i.e., the magnitude of the residual on a color scale, in the right figure. The heat map, and surface plots in this section, are constructed for $N_y = 1001$. The figure shows convergence of the solution, although at a slower rate than for continuous boundary values. The heat map also shows a few characteristics given by the solid white and dashed yellow curves. Furthermore, it shows that the discontinuity of the derivatives of $a$ in $(0, -0.5)$ yields a locally distinct residual. This difference in residual is propagated along the characteristic starting in $(0, -0.5)$. This coincides with an alternative equivalent definition of a characteristic, from~\cite[p. 408]{HilbertCourant}: ``Discontinuities (of a nature to be specified later) of a solution
    cannot occur except along characteristics.''. We add to this that the discontinuities mentioned, only arise in the second order derivatives, and $u$, $p$ and $q$ are smooth as seen in Figure~\ref{fig:discontinuousBoundarySurface}. The figures also show two characteristics, departing from the end points of the initial strip. Furthermore, Figure~\ref{fig:discontinuousBoundarySurface_s_tx} shows $r$ to be continuous but nonsmooth (left), and its derivative $r_y$ to be discontinuous (right). This is to be expected as the second derivatives $r,s,t$ correspond to $a$ and $b$ via~\eqref{eqn:equivalenceRelations}, while roughly speaking the derivatives of $r,s,t$ correspond to the derivatives of $a$, $b$ and $f$.
    
    \section{Conclusions}
    \label{sec:Conclusion}
    In this paper, we presented the general formulation of the method of characteristics for a nonlinear second order hyperbolic partial differential equation (PDE) in two variables. We derived conditions which determine whether a curve is characteristic. Along these characteristics the PDE reduces to two systems of ordinary differential equations which are mutually coupled. These ODE systems can be solved using explicit numerical integrators. We presented three such integrators, which are based on one-step methods. The computed characteristics will not pass through grid points. Therefore interpolation is necessary which should be handled carefully as not to spoil numerical convergence. We discussed how the direction of the characteristics at the boundary determines the number of boundary conditions which should be prescribed.
    
    For test cases with known analytical solutions the developed methods are shown to converge to the analytical solution up and till computer precision. Furthermore, two measures for the residual are formulated which seemingly converge to computer precision. The method is shown to work for an example where the initial strip is extended over two boundary segments, and for another example for which the number of boundary conditions necessary vary along a boundary segment. A nonsmooth boundary condition was imposed to show the propagation of the nonsmooth boundary data along the characteristics while the algorithm remained stable.
    
    In future work we would like to extend the algorithm to more general variants of the \MAe{}. We intend to use the numerical methods to design hyperbolic optical freeform surfaces as has been done in \cite{CorienThesis} for the elliptic \MAe{}.

    \bibliographystyle{unsrt}
    \bibliography{refs}
    
    \begin{appendices}
    \section{Interpolation}
    \label{sec:numericalInterpolation}
    In this section we will briefly introduce splines, which we use for numerical interpolation. To understand spline interpolation, we first introduce knot sequences, which generate B-splines, and in term determine the spline interpolant. To this end consider a set of $N+1$ numbers $\pmb{\xi} = \{\xi_0, \xi_1, \dots, \xi_N\}$ with $\xi_0 \leq \xi_1 \leq \dots \leq \xi_N$. Such a sequence is called a knot sequence and each member of the sequence is called a knot. 
    B-splines $b^n_k$ of degree $n$ for the knot sequence $\pmb{\xi}$ are recursively defined on the interval $[\xi_k, \xi_{k+n+1})$ by \cite[p. 52]{splinesBook}
    \begin{align}
    b_k^n(t) = \gamma^n_k(t) b^{n-1}_k(t) + \Big(1-\gamma^n_{k+1}(t)\Big) b^{n-1}_{k+1}(t), \hspace{10pt} \gamma^n_k(t) = \frac{t-\xi_k}{\xi_{k+n} - \xi_k},
    \end{align}
    for $ 0 \leq k$, $k + n + 1 \leq N$, $0 \leq n \leq k$ with initial values
    \begin{align}
    b^0_k(t) = \begin{cases}
    1, & \text{if } \xi_k \leq t < \xi_{k+1}, \\
    0, & \text{otherwise.}
    \end{cases}
    \end{align}
    Each B-spline $b^n_k$ is a polynomial, of degree $\leq n$ on its knot interval $[\xi_k, \xi_{k+n+1})$, and vanishes outside this interval. \\
    
    Let $m \in \mathbb{N}_{+}$, let $g$ be a sufficiently smooth function and let $t_0 \leq \dots \leq t_{m-1}$ be a set of points such that $g(t_0), \dots, g(t_{m-1})$ are known. Furthermore, let $N = m + n$, so $\pmb{\xi} = \{\xi_0, \xi_1, \dots, \xi_{m+n}\}$ and let $\pmb{\xi}$ be such that 
    \begin{align}
    \label{eqn:knotSequenceCondition}
    \xi_k < t_k < \xi_{k+n+1} \text{ for  }k = 0, \dots, m-1.
    \end{align}
    These conditions, also known as the Schoenberg-Whitney conditions \cite[p. 91]{splinesBook}, imply that there exists a unique interpolating spline 
    \begin{align}
    P(t) = \sum_{k=0}^{m-1} c_k b_k^n(t),
    \end{align}
    of degree $n$, which interpolates $g$ in the interval $t\in[t_0, t_{m-1})$. The coefficients $c_k$ are calculated via the implicit relation
    \begin{align}
    \mathbf{A} \mathbf{c} = \mathbf{g}, \hspace{10pt}
    (a_{j, k}) = b_k^n(t_j), \hspace{10pt} \mathbf{c} = (c_k),\hspace{10pt} \mathbf{g} = (g(t_k)),
    \end{align}
    for $j,k = 0, 1, \dots, m-1$.
    Furthermore, the associated error can be estimated by
    \begin{align}
    \label{eqn:convergeSpline}
    |g(t) - P(t)| \leq C\Big(n, \|\mathbf{A}^{-1}\|_{\infty}\Big) \, \|g^{(n+1)}\|_{\infty, R} \, h^{n+1}, \hspace{10pt} t \in R,
    \end{align}
    where $\|g^{(n+1)}\|_{\infty, R}$ is the maximum norm of the derivative of $g$ of order $n+1$ on the interval $R = [\xi_n, \xi_{n+m}]$, $h = (\xi_{j+1} - \xi_j)$ and the constant $C$ depends on the degree $n$ of the B-splines and the infinity-norm of the inverse of $\mathbf{A}$. 
    A convergence order for odd and even $n$ has been established, where generally convergence for odd $n$ is of order $\mathcal{O}(h^{n+1})$ as given by~\eqref{eqn:convergeSpline}. Convergence for even $n$ has been observed to be of order $\mathcal{O}(h^{n+2})$ instead of the theoretical established upper bound \eqref{eqn:convergeSpline}. This observed superior convergence for even $n$  is not fully understood at the time of writing \cite{Volkov2019}.
    
    What remains is to construct a suitable knot sequence such that the \hfill\break\mbox{Schoenberg-Whitney} condition holds.  Let data points $t_0 \leq \dots \leq t_{m-1}$ be such that $m > n$, i.e., let the number of data points exceeds the degree of the B-splines used for interpolation. Then we choose the knot sequence $\pmb{\xi}$ according to
    \begin{align}
    \label{eqn:aptknt}
    \begin{cases}
    \xi_{i} = t_0 & i = 0, \dots, n-1, \\
    \xi_{i} = \frac{1}{n-1} \sum_{j = 1}^{n-1} t_{i - n + j}, & i = n, \dots, m, \\
    \xi_{i} = t_{m} & i = m+1, \dots, m + n - 1.
    \end{cases}
    \end{align}
    The first and last terms of $\pmb{\xi}$ are the original, possibly duplicated, starting and end values $t_0$ and $t_{m}$, respectively. The remaining components are running averages of size $n-1$ of $\{t_j\}$, which ensures~\eqref{eqn:knotSequenceCondition} holds.
    
    \textit{Example: } Consider the ordered sequence of data points $t: \{1, 2, 4, 5, 7, 10\}$ and let the desired spline order be $n = 3$. Hence $m=6$ and the knot sequence $\pmb{\xi}$ according to~\eqref{eqn:aptknt} is given by $\xi: \{1, 1, 1, 3, 4.5, 6, 10, 10, 10\}$. \\
    
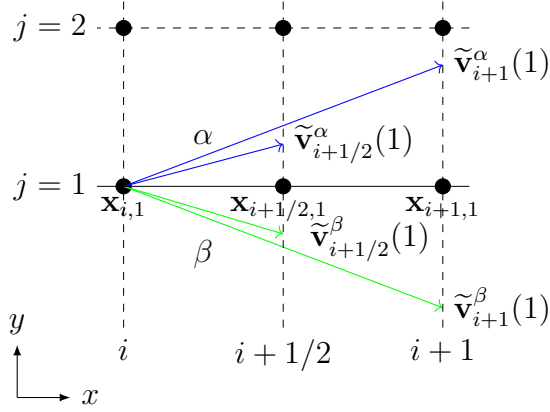
\begin{figure}[H]
        \centering
        \captionsetup{width=0.85\linewidth}
        \begin{tikzpicture}[scale=0.7]
        \coordinate (Origin)   at (0,0);
        \coordinate (XAxisMin) at (-2,-2);
        \coordinate (XAxisMax) at (-1,-2);
        \coordinate (YAxisMin) at (-2,-2);
        \coordinate (YAxisMax) at (-2,-1);
        \draw [thin,-latex] (XAxisMin) -- (XAxisMax) node [right] {$x$};
        \draw [thin,-latex] (YAxisMin) -- (YAxisMax) node [above] {$y$};
        
        
        \node[draw,circle,inner sep=2pt,fill] at (0,2) {};                    
        \node[draw,circle,inner sep=2pt,fill] at (0,5) {};
        \node[draw,circle,inner sep=2pt,fill] at (3,2) {};                    
        \node[draw,circle,inner sep=2pt,fill] at (3,5) {};
        \node[draw,circle,inner sep=2pt,fill] at (6,2) {};
        \node[draw,circle,inner sep=2pt,fill] at (6,5) {};
        
        \node [below] at (0,2) {$\mathbf{x}_{i,1}$};
        \node [below] at (2.92,2) {$\mathbf{x}_{i+1/2,1}$};
        \node [below] at (6,2) {$\mathbf{x}_{i+1,1}$};
        
        \draw[->, blue] (0,2) -- (3, 2.8) node [right, text=black] {$\widetilde{\mathbf{v}}^\alpha_{i+1/2}(1)$};
        \draw[->, green] (0,2) -- (3, 1.1) node [right, text=black, xshift=6pt, yshift=-2pt] {$\widetilde{\mathbf{v}}^\beta_{i+1/2}(1)$};
        
        \draw[->, blue] (0,2) -- (6, 4.3) node [right, text=black] {$\widetilde{\mathbf{v}}^\alpha_{i+1}(1)$};
        \draw[->, green] (0,2) -- (6, -0.3) node [right, text=black] {$\widetilde{\mathbf{v}}^\beta_{i+1}(1)$};
        
        \node[above] at (1.5, 2.6) {$\alpha$};
        \node[below] at (1.5, 1.2) {$\beta$};
        
        \draw[dashed, line width=0mm] (0, 5.5) -- (0,-0.7)node [below] {$i$};
        \draw[dashed, line width=0mm] (3, 5.5) -- (3,-0.7) node [below] {$i+1/2$};
        \draw[dashed, line width=0mm] (6, 5.5) -- (6,-0.7) node [below] {$i+1$};
        
        \draw[dashed, line width=0mm] (6.5, 5) -- (-0.5,5) node [left] {$j=2$};
        \draw[line width=0mm] (6.5, 2) -- (-0.5,2) node [left] {$j=1$};
        
        
        \end{tikzpicture}        
        \caption{ Modified Euler based method near the lower boundary of the domain.}
        \label{fig:multistageExtrapolation}
    \end{figure}

    In the case that we require extrapolation at a point $t_e < t_0$ or $t_e > t_{m-1}$, we simply estimate $g(t_e) \approx P(t_e)$.     
    This need for extrapolation does occur for both the modified Euler and Runge-Kutta based methods, as shown in Figure~\ref{fig:multistageExtrapolation} for the modified Euler scheme, due to the predictor lying outside of $\Omega$ and missing either the value $a$, or $b$, as it cannot be determined along the characteristic. Let without loss of generality $a > b$, then in order to approximate $b(x_{i+\frac{1}{2}}, \tilde{y}^\beta_{i+\frac{1}{2}}(1))$ which is needed to calculated $\widetilde{\mathbf{v}}^\beta_{i+1}(1)$ according to~\eqref{eqn:modifiedEulerFullStep}, we extrapolate using a spline based on $\tilde{\mathbf{v}}^\alpha_{i+1/2}(1)$ for $j=1,\dots,m+1$, known at the $y$-values $\tilde{y}^\alpha_{i+\frac{1}{2}}(1)$ inside the domain, to $\tilde{y}^\beta_{i+\frac{1}{2}}(1)$ outside the domain.
    
    \section[title]{Generating Solutions to the Hyperbolic \hfill \\ \MAE{}}
    \label{app:generateSolutionPair}
        Finding solutions to the \MAe{} can be problematic due to the hyperbolic and nonlinear nature of $u_{xx} u_{yy} - u_{xy}^2 = -f^2 < 0$ with $f \neq 0$. Therefore we introduce a method based on complex functions to quickly obtain a pair $(u, f)$ which solves the \MAe{}. To this end let $w$ be a complex analytical function and let $u(x,y) = \text{Re}(w(x + i y))$. Differentiation then yields
        \begin{equation*}
        \begin{aligned}
        u_{xx} & = \text{Re}(w^{\prime\prime}), \\
        u_{xy} & = \text{Re}(i w^{\prime\prime}) = - \text{Im}(w^{\prime\prime}), \\
        u_{yy} & = -\text{Re}(w^{\prime\prime}).
        \end{aligned} 
        \end{equation*}
        For $f$ it then follows
        \begin{align*}
        f^2 = -u_{xx} u_{yy} + u_{xy}^2 = \Big(\text{Re}(w^{\prime\prime})\Big)^2 + \Big(\text{Im}(w^{\prime\prime})\Big)^2 = |w^{\prime\prime}|^2.
        \end{align*}
        The implications of this are twofold. First, given $w$ we can construct $u(x,y) = \text{Re}(w(x + i y))$ and $f(x,y) = |w^{\prime\prime}(x + i y)|$, which form a solution. Secondly, if for given $f$ there exists an analytical function $w$ such that $f^2(x,y) = |w^{\prime\prime}|^2$, then $u(x,y) = \text{Re}(w(x + i y))$ solves the \MAe{}.

        Alternatively, considering $u(x,y) = -\text{Re}(w(x + i y))$ or $u(x,y) = \pm\text{Im}(w(x + i y))$ instead, yields the same conclusions as above.
        
    \section*{Acknowledgments}
        Special thanks goes to J. de Graaf, for contributing App.~\ref{app:generateSolutionPair}, which has greatly simplified generating test examples for the \MAe{}. Contact: CASA, Department of Mathematics and Computer Science, Eindhoven University of Technology, PO Box 513, 5600MB Eindhoven, The Netherlands\\
        
        \noindent This work is part of the research programme NWO-TTW Perspectief with project number P15-36, which is (partly) financed by the Netherlands Organisation for Scientific Research (NWO).

    \end{appendices}

    \end{document}